\documentclass[labelcite,aap,MSNbibl,seceqn,citesort,dvips]{arximspdf}

\doi{10.1214/11-AAP822} 
\volume{22}
\issue{5}
\pubyear{2012}
\firstpage{1989}
\lastpage{2047}

\makeatletter

\newtheorem{theorem}{Theorem}[section]
\newtheorem{lemma}[theorem]{Lemma}
\newtheorem{proposition}[theorem]{Proposition}

\newproclaim{example}[theorem]{Example}
\newproclaim{remark}[theorem]{Remark}

\newcommand{\ntoo}{{n\to\infty}}

\newcommand{\downto}{\searrow}

\newcommand{\tend}{\longrightarrow}
\newcommand{\dto}{\stackrel{\mathrm{d}}{\tend}}
\newcommand{\pto}{\stackrel{\mathrm{p}}{\tend}}

\newcommand{\psim}{\stackrel{\mathrm{p}}{\sim}}
\newcommand{\eqd}{\stackrel{\mathrm{d}}{=}}

\newcommand{\Op}{O_{\mathrm p}}

\newcommand{\bbR}{\mathbb R}
\newcommand{\bbN}{\mathbb N} 
\newcommand{\E}{\operatorname{\mathbb E{}}}

\newcommand{\Var}{\operatorname{Var}}
\newcommand{\Cov}{\operatorname{Cov}}
\newcommand{\Po}{\operatorname{Po}}
\newcommand{\Bin}{\operatorname{Bin}}
\newcommand{\Mul}{\operatorname{Mul}}
\newcommand{\Be}{\operatorname{Be}}
\newcommand{\NBi}{\operatorname{NegBin}}

\newcommand{\ga}{\alpha}
\newcommand{\gb}{\beta}
\newcommand{\gD}{\Delta}
\newcommand{\gf}{\varphi}
\newcommand{\gam}{\gamma}
\newcommand{\gl}{\lambda}
\newcommand{\go}{\omega}
\newcommand{\gs}{\sigma}
\newcommand{\gss}{\sigma^2}

\newcommand{\gth}{\theta}
\newcommand{\gthx}{\vartheta}
\newcommand{\eps}{\varepsilon}

\newcommand{\cA}{\mathcal A}
\newcommand{\cB}{\mathcal B}
\newcommand{\cF}{\mathcal F}
\newcommand{\cG}{\mathcal G}
\newcommand{\cP}{\mathcal P}
\newcommand{\cQ}{\mathcal Q}
\newcommand{\cZ}{{\mathcal Z}}

\newcommand{\hU}{{\widehat U}}
\newcommand{\tX}{{\tilde X}}

\newcommand{\tpi}{{\tilde\pi}}

\newcommand{\ett}{{\mathbf1}}

\newcommand{\qqqq}{^{1/4}}
\newcommand{\qqq}{^{1/3}}
\newcommand{\qqw}{^{-1/2}}
\newcommand{\qw}{^{-1}}
\newcommand{\qww}{^{-2}}
\newcommand{\qwr}{^{-1/r}}

\newcommand{\intoooo}{\int_{-\infty}^\infty}
\newcommand{\oi}{[0,1]}
\newcommand{\ooo}{[0,\infty)}
\newcommand{\oooo}{(-\infty,\infty)}

\newcommand{\dtv}{d_{\mathrm{TV}}}

\newcommand{\dd}{\,\mathrm{d}}
\newcommand{\ddq}{\mathrm{d}}
\newcommand{\ddt}{\frac{\ddq}{\ddq t}}

\newcommand{\gnm}{G(n,m)}

\newcommand{\iis}{I_i(s)}
\newcommand{\iit}{I_i(t)}
\newcommand{\cao}{\cA(0)}
\newcommand{\cDA}{\Delta\cA}
\newcommand{\Ax}{A^*}
\newcommand{\yix}{Y_1}
\newcommand{\mix}{M_1}
\newcommand{\cx}{_\mathsf{c}}
\newcommand{\cxx}{_\mathsf{c}^*}
\newcommand{\tc}{t\cx}
\newcommand{\ac}{a\cx}
\newcommand{\bc}{b\cx}
\newcommand{\bcq}{\bc'}
\newcommand{\bco}{b}
\newcommand{\bcx}{b^*}
\newcommand{\pc}{p\cx}
\newcommand{\tcx}{t\cxx}
\newcommand{\tcxx}{t\cx^{**}}
\newcommand{\acx}{a\cxx}
\newcommand{\pcx}{p\cxx}
\newcommand{\ccr}{c\cx}
\newcommand{\gthc}{\gth\cx}
\newcommand{\ff}{\bar f}

\newcommand{\ijx}{_{ij}}
\newcommand{\hh}{\gamma}
\newcommand{\aox}{A_0}
\newcommand{\aoxx}{A_0'}
\newcommand{\AsN}{\mathrm{AsN}}
\newcommand{\tf}{T^F}

\newcommand{\nie}{\tau'_\eps}

\newcommand{\tauw}{\tau'}
\newcommand{\bd}{c_\delta}
\newcommand{\gox}{\go'}
\newcommand{\goxx}{\go''}
\newcommand{\smu}{S^\mu}
\newcommand{\xo}{x_0}
\newcommand{\xox}{x_1}
\newcommand{\ddx}{\frac{\partial}{\partial x}}

\newcommand{\gthcq}{\gth\cx^-}
\newcommand{\gthcc}{\gth^*\cx}
\newcommand{\tT}{\tilde T}
\newcommand{\tsk}{\tilde S_k}
\newcommand{\xp}{p'}
\newcommand{\gx}{g}
\newcommand{\gy}{\tilde g}

\newcommand{\tS}{\tilde S}

\makeatother

\begin{document}
\begin{frontmatter}

\title{Bootstrap percolation on the random graph $G_{n,p}$}
\runtitle{Bootstrap percolation on $G_{n,p}$}

\begin{aug}
\author[A]{\fnms{Svante} \snm{Janson}\corref{}\ead[label=e1]{svante@math.uu.se}},
\author[B]{\fnms{Tomasz} \snm{{\L}uczak}\thanksref{t2}\ead[label=e2]{tomasz@amu.edu.pl}},
\author[C]{\fnms{Tatyana} \snm{Turova}\thanksref{t3}\ead[label=e3]{tatyana@maths.lth.se}}\\
\and
\author[D]{\fnms{Thomas} \snm{Vallier}\ead[label=e4]{vallierthomas@gmail.com}}
\runauthor{Janson, {\L}uczak, Turova and Vallier}
\affiliation{Uppsala University, Adam Mickiewicz University,
Lund University and~Helsinki~University}
\address[A]{S. Janson\\
Matematiska institutionen\\
Uppsala Universitet\\
Box 480\\
751 06 Uppsala\\
Sweden\\
\printead{e1}}
\address[B]{T. {\L}uczak\\
Collegium Mathematicum\\
Adam Mickiewicz University\hspace*{1.49pt}\\
Umultowska 87\\
61-614 Poznan\\
Poland\\
\printead{e2}}
\address[C]{T. Turova\\
Matematikcentrum\\
Lunds Universitet\\
Box 117\\
221 00 Lund\\
Sweden\\
\printead{e3}}
\address[D]{T. Vallier\\
Department of Mathematics\\
\quad and Statistics\\
University of Helsinki\\
P.O. Box 68\\
Gustaf H\"allstr\"omin katu 2b\\
FI-00014 Helsinki\\
Finland\\
\printead{e4}} 
\end{aug}

\thankstext{t2}{Supported in part by the Foundation for Polish Science.}

\thankstext{t3}{Supported in part by the Swedish Science Foundation.}

\received{\smonth{12} \syear{2010}}
\revised{\smonth{10} \syear{2011}}

%
\begin{abstract}
Bootstrap percolation on the random graph $G_{n,p}$ is a process of
spread of ``activation'' on a given realization of the graph with a
given number of initially active nodes. At each step those vertices
which have not been active but have at least $r\geq2$ active neighbors
become active as well.

We study the size $A^*$ of the final active set. The parameters of the
model are, besides $r$ (fixed) and $n$ (tending to $\infty$), the size
$a=a(n)$ of the initially active set and the probability $p=p(n)$ of
the edges in the graph. We show that the model exhibits a sharp phase
transition: depending on the parameters of the model, the final size of
activation with a high probability is either $n-o(n)$ or it is $o(n)$.
We provide a complete description of the phase diagram on the space of
the parameters of the model. In particular, we find the phase
transition and compute the asymptotics (in probability) for $A^*$; we
also prove a~central limit theorem for $A^*$ in some ranges.
Furthermore, we provide the asymptotics for the number of steps until
the process stops.
\end{abstract}

%

\end{frontmatter}

\section{Introduction}\label{sec1}

Bootstrap percolation on a graph $G$ is defined as the spread of
\textit{activation} or \textit{infection} according to
the following rule, with a given threshold $r\ge2$:
We start with a set $\cao\subseteq V(G)$ of
\textit{active} vertices.
Each inactive vertex that has at least $r$ active neighbors
becomes active. This is repeated until no more vertices become active,
that is, when no inactive vertex has $r$ or more active
neighbors.
Active vertices never become inactive, so the set of active vertices
grows monotonously.

To avoid confusion, we will use the terminology that each active vertex
\textit{infects} all its neighbors, so that a vertex that is infected (at
least) $r$ times becomes active.

We are mainly interested in the final size $\Ax$ of the active set,
and in
particular whether eventually all vertices will be active or not.
If they are, we say that the initial set $\cao$
\textit{percolates} (completely). We will study a~sequence of graphs
of order
$n\to\infty$; we then also say that (a sequence of)
$\cao$ \textit{almost percolates} if the number of
vertices that remain inactive is $o(n)$, that is, if $\Ax=n-o(n)$.

Bootstrap percolation on a lattice
(which is a~special example of a~cellular automata)
was introduced in 1979 by Chalupa, Leath and Reich \cite{ChalupaLR}
as a~simplified model of some magnetic systems.
Since then bootstrap percolation has been studied on various graphs, both
deterministic and random. One can study either a~random initial set or
the deterministic problem of choosing an initial set that is optimal
in some sense.
A simple example of the latter is the classical folklore problem to
find the
minimal
percolating set in a two-dimensional grid
(i.e., a finite square $[n]^2$ in the square lattice); see
Balogh and Pete \cite{BaloghPete} and
Bollob{\'a}s \cite{BollArt}.
(These references also treat higher-dimensional grids $[n]^d$.)
Another extremal problem is studied by Morris \cite{Morris}.
The problem with a random initial set was introduced by
Chalupa, Leath and Reich \cite{ChalupaLR} (lattices and regular
infinite tree),
and further studied on lattices by Schonmann \cite{Schonmann};
it has, in particular, been studied on finite
grids (in two dimensions or more),
see
Aizenman and Lebowitz \cite{AizenmanL}, 
Balogh and Pete~\cite{BaloghPete}, 
Cerf and Cirillo \cite{CerfC}, 
Cerf and Manzo \cite{CerfManzo}, 
Holroyd \cite{Holroyd}, 
Balogh, Bollob{\'a}s and Morris \cite{BaloghBM}, 
Gravner, Holroyd and Morris
\cite{GravnerHM}. 
In a recent paper, Balogh et al. \cite{BaloghBD-CM} derived a sharp
asymptotic for the critical density (i.e., the critical size of a
random initial set) for bootstrap percolation on grids of any
dimension, generalizing results of Balogh, Bollob{\'a}s, and Morris
\cite{BaloghBM2}. Grids with a different edge set where studied by
Holroyd, Liggett and Romik \cite{HolroydLR}. The study of bootstrap
percolation on lattices is partly explained by its origin in
statistical physics, and the bootstrap process is being successfully
used in studies of the Ising model; see
\cite{CerfManzo1,CerfManzo2,FontesSS,Morris1}. Lately bootstrap
percolation has also been studied on varieties of graphs different from
lattices and grids; see, for example, Balogh and Bollob{\'a}s
\cite{BaloghBollobas} (hypercube); Balogh, Peres and Pete
\cite{BaloghPeresPete} (infinite trees); Balogh and Pittel
\cite{BaloghPittel}, Janson \cite{SJ215} (random regular graphs); an
extension where the threshold may vary between the vertices is studied
by Amini \cite{Amini}. An anisotropic bootstrap percolation was studied
by Duminil-Copin and van Enter \cite{D-CvE}. Further, a graph bootstrap
percolation model introduced by Bollob{\'a}s \cite{Bo68} already in
1968, where edges are infected instead of vertices, was analyzed
recently by Balogh, Bollob{\'a}s and Morris~\cite{BaloghBM3} and
Balogh et al. \cite{BaloghBMO}.

In the present paper, we study bootstrap percolation on the
Erd\"os--R\'enyi random graph $G_{n,p}$ with an initial set $\cao$
consisting of a given number $a$ of vertices chosen at random. (By
symmetry, we obtain the same results for any deterministic set of $a$
vertices.) Recall that $G_{n,p}$ is the random graph on the set of
vertices $V_n = \{ 1,\ldots,n\}$ where all possible edges between pairs
of different vertices are present independently and with the same
probability~$p$. As usual, we let $p=p(n)$ depend on $n$.

A problem equivalent to bootstrap percolation on $G_{n,p}$ in the case
$p=\gl/n$ was studied by Scalia-Tomba \cite{Scalia-Tomba}, although he
used a different formulation as an epidemic. (Ball and Britton
\cite{BallBritton,BallBritton09} study a more general model with
different degrees of severity of infection.) Otherwise, bootstrap
percolation on $G_{n,p}$ was first studied by Vallier
\cite{VallierPhD}; we here use a simple method (the same as
\cite{Scalia-Tomba}) that allows us to both simplify the proofs and
improve the results. We will state the results for a general fixed
$r\ge2$ (the case $r=1$ is much different; see Remark~\ref{R1}); the
reader may for simplicity consider the case $r=2$ only, since there are
no essential differences for higher $r$.

We will see that there is a threshold phenomenon: typically, either the
final size $\Ax$ is small (at most twice the initial size $a$), or it
is large [sometimes exactly $n$, but if $p$ is so small that there are
vertices of degree less than $r$, these can never become active except
initially so eventually at most $n-o(n)$ will become infected]. We can
study the threshold in two ways: in the first version, we keep $n$ and
$p$ fixed and vary $a$. In the second version, we fix $n$ and~$a$ and
vary $p$. We will state some results for both versions and some for the
former version only; these too can easily be translated to the second
version. We will also study dynamical versions, where we add new
external infections or activations or new edges until we reach the
threshold; see Section~\ref{Sdyn}.

Apart from the final size $\Ax$, we will also study the time $\tau$ the
bootstrap process takes until completion. We count the time in
\textit{generations}: generation~0 is $\cao$, generation 1 is the set of
other vertices that have at least $r$ neighbors in generation~0, and so
on. The process stops as soon as there is an empty generation, and we
let $\tau$ be the number of (nonempty) generations. Thus, if we let
$\cG_k$ be the set of vertices activated in generation $k$, then
%
%
\begin{equation}\label{tau}
\tau:=\max\{k\ge0\dvtx\cG_k\neq\varnothing\}=\min\{k\ge1\dvtx
\cG_k=\varnothing\}-1.
\end{equation}
%

\begin{remark}
Bootstrap percolation does not seem to be a good model for usual infectious
diseases; see, however, Ball and Britton \cite{BallBritton}.
It might be a better
model for the spread of rumors or other
ideas or beliefs; cf. the well-known rule, ``What I tell you three
times is true'' in
Carroll \cite{snark}.

Bootstrap percolation can be also viewed as a simplified model for
propagation of activity in a neural network. Although related neuronal
models are too involved for a rigorous analysis (see, e.g.,
\cite{Ko,Tl,TV}) they
inspired study of bootstrap percolation on $G_{n,p}$ by Vallier
\cite{VallierPhD}. There is a further discussion on the
application of bootstrap percolation on $G_{n,p}$ to neuromodelling in
\cite{T}.
\end{remark}
%
%
\begin{remark}
Instead of $G_{n,p}$,
one might consider the random graph
$\gnm$, with a given number $m=m(n)$ of edges. It is easy to obtain
a~result for $\gnm$ from
our results for $G_{n,p}$, using monotonicity, but we usually leave
this to the reader. [In the dynamical model in Section~\ref{SSdynm}, we consider
$\gnm$, however.]
\end{remark}
%
%
\begin{remark}
An alternative to starting with an initial active set of fixed size $a$
is to let each vertex be initially activated with probability $\ga=\ga(n)>0$,
with different vertices activated independently. Note that
this is the same as taking
the initial size $a$ random with $a\in\Bin(n,\ga)$. For most results the
resulting random variation in $a$ in negligible, and we obtain the same
results as for $a=n\ga$, but for the Gaussian limit in Theorems
\ref{Tac}\ref{Tac0} and~\ref{TDG2}, the asymptotic variances are
changed by
constant factors. We leave the details to the reader.
\end{remark}

Some open problems arise from our study.
In \cite{BaloghBM}, Balogh,
Bollob{\'a}s and Morris determine the critical probability for bootstrap
percolation on grids when the dimension $d=d(n) \to\infty$. A similar idea
translated to the
$G(n,p)$ graph would be to study what happens when $r=r(n) \to\infty
$. This
problem is not treated here although our methods
might be useful also for such problems. The problem of majority
percolation where a vertex becomes activated if at least half of its
neighbors are active [$r(v) = d(v)/2$] has been studied on the
hypercube by
Balogh, Bollob{\'a}s and Morris \cite{BaloghBM1}. On the $d$-dimensional
grid $d(v)/2 = d$ but on the $G(n,p)$ graph, this problem is completely
different and still open.
(We thank the referee for these suggestions.)

The method is described in Section~\ref{Ssetup}. The main results are
stated in Section~\ref{Smain}, with further results in Sections
\ref{Sdyn} and~\ref{Sbound}. Proofs are given in Sections~\ref{Soverview}--\ref{Slast}.

\subsection{Notation}\label{sec11}

All unspecified limits are as $\ntoo$.
We use $\dto$ for convergence in distribution and
$\pto$ for convergence in probability of random variables;
we further use $\Op$ and $o_{\mathrm p}$ in the standard sense (see, e.g.,
\cite{SJN6} and \cite{JLR}), and we use w.h.p. (with high
probability) for events with probability tending to 1 as $\ntoo$.
Note that, for example, ``$=o(1)$ w.h.p.'' is equivalent to
``$=o_p(1)$''
and to ``$\pto0$,''
and that ``$\sim a_n$ w.h.p.'' is equivalent to ``$=(1+o_p(1))a_n$;''
see \cite{SJN6}.
A statement of the type ``when $\cP$, then w.h.p. $\cQ$'' (or similar
wording), where $\cP$ and $\cQ$ are two events, means that
${\mathbb P}(\cP\mbox{ and (not }\cQ))\to0$, that is, that w.h.p.
``(not $\cP$) or $\cQ$'' holds.
(See, e.g., Theorem~\ref{T2} and Proposition~\ref{PG3}.)
If ${\mathbb P}(\cP)$ is bounded away from 0, this is equivalent to
``conditioned on
$\cP$, $\cQ$ holds w.h.p.''

If $X_n$ is a sequence of random variables, and $\mu_n$ and $\gss_n$
are sequences of real numbers, with $\gss_n>0$, we say that
$X_n\in\AsN(\mu_n,\gss_n)$ if $(X_n-\mu_n)/\gs_n\dto N(0,1)$.

Occasionally we use the subsubsequence principle (\cite{JLR}, page 12),
which says that to prove a limit result (e.g., for real numbers, or
for random variables in probability or in distribution), it is sufficient
to show that every subsequence has a subsubsequence where the result holds.
We may thus, without loss of generality, select convenient subsequences
in a
proof, for example, such that another given sequence either converges or
tends to $\infty$.

\section{A useful reformulation}\label{sec2}\label{Ssetup}

In order to analyze the bootstrap percolation process on
$G_{n,p}$, we change the time scale; we forget the generations and consider
at each time step
the infections from one vertex only.
Choose $u_1\in\cao$ and give each of its neighbors a \textit{mark};
we then say
that $u_1$ is \textit{used}, and let
$\cZ(1):=\{u_1\}$ be the set of used vertices
at time 1.
We continue recursively: at time $t$, choose a vertex
$u_{t}\in\cA(t-1)\setminus\cZ(t-1)$.
We give each neighbor of $u_{t}$ a new mark. Let
$\cDA(t)$ be the set of inactive vertices with $r$
marks; these now become active, and we let
$\cA(t)=\cA(t-1)\cup\cDA(t)$ be the set of
active vertices at time $t$. We finally set
$\cZ(t)=\cZ(t-1)\cup\{u_{t}\}=\{u_s\dvtx s\le t\}$,
the set of used vertices. [We start with
$\cZ(0)=\varnothing$, and note that necessarily
$\cDA(t)=\varnothing$ for $t<r$.]

The process stops when
$\cA(t)\setminus\cZ(t)=\varnothing$, that is, when
all active vertices are used. We denote this time by $T$,
%
%
\begin{equation}
\label{t1}
T:=\min\{t\ge0\dvtx\cA(t)\setminus\cZ(t)=\varnothing\}.
\end{equation}
Clearly, $T\le n$; in particular, $T$ is finite. The final active set
is $\cA(T)$; it is clear that this is the same set as the one produced
by the bootstrap percolation process defined in the
\hyperref[sec1]{Introduction}; only the time development differs. Hence
we may as well study the version just described. [This is true for any
choice of the vertices~$u_t$. For definiteness, we may assume that we
keep the unused, active vertices in a queue and choose $u_{t}$ as the
first vertex in the queue, and that the vertices in $\cDA(t)$ are added
at the end of the queue in order of their labels. Thus $u_{t}$ will
always be one of the oldest unused, active vertices, which will enable
us to recover the generations; see further Section~\ref{Sgen}. In
Section~\ref{Sdyn}, we consider other ways of choosing $u_t$.] This
reformulation was used already by Scalia-Tomba \cite{Scalia-Tomba} (for
a more general model). It is related to the (continuous-time)
construction by Sellke \cite{Sellke} for an epidemic process.

Let $A(t):=|\cA(t)|$, the number of active vertices at time $t$.
Since $|\cZ(t)|=t$ and $\cZ(t)\subseteq\cA(t)$
for $t=0,\ldots,T$, we also have
%
%
\begin{equation}
\label{t2}
T=\min\{t\ge0\dvtx A(t)=t\}
=
\min\{t\ge0\dvtx A(t)\le t\}.
\end{equation}
Moreover,
since the final active set is $\cA(T)=\cZ(T)$,
its size $\Ax$ is
%
%
\begin{equation}\label{at}
\Ax:=A(T)=|\cA(T)|=|\cZ(T)|=T.
\end{equation}
Hence, the set $\cao$ percolates if
and only if $T=n$, and $\cao$ almost percolates if and only if $T=n-o(n)$.

We analyze this process by the standard method of revealing the edges
of the graph $G_{n,p}$ only on a need-to-know basis. We thus
begin by choosing~$u_1$ as above and then reveal its neighbors;
we then find $u_2$ and reveal its neighbors, and so on. Let,
for $i\notin\cZ(s)$, $\iis$ be the indicator
that there is an edge between the vertices $u_s$ and
$i$.
This is also the indicator that $i$ gets a mark at time
$s$, so if $M_i(t)$ is the number of marks $i$ has
at time $t$, then
%
%
\begin{equation}
\label{mi}
M_i(t)=\sum_{s=1}^t \iis,
\end{equation}
at least until $i$ is activated (and what happens later does not matter).
Note that if $i\notin\cao$, then, for every $t\le T$,
$i\in\cA(t)$ if and only if $M_i(t)\ge r$.
The crucial feature of this description of the process, which
makes the analysis simple, is that the random variables
$\iis$ are i.i.d. $\Be(p)$.

We have defined
$\iis$ only for $s\le T$ and
$i\notin\cZ(s)$, but it is convenient to add further (redundant)
variables so that $\iis$ are defined, and i.i.d., for all $i\in
V_n$ and all $s\ge1$.
One way to do this formally is to reverse the procedure above. We
start with i.i.d. $\iis\in\Be(p)$,
for $i\in V_n$ and $s\ge1$, and a set
$\cao\subseteq V_n$. We let $\cZ(0):=\varnothing$ and start with an empty
graph on~$V_n$.
We then, as above, for $t=1,\ldots,n$ select
$u_t\in\cA(t-1)\setminus\cZ(t-1)$ if this set is
nonempty; otherwise we select $u_t\in
V_n\setminus\cZ(t-1)$ (taking, e.g., the smallest such
vertex).
We define $M_i(t)$ by (\ref{mi}) for all $i\in
V_n$ and $t\ge0$, and update
$\cA(t):=\cao\cup\{i\dvtx M_i(t)\ge r\}$ and
\mbox{$\cZ(t):=\cZ(t-1)\cup\{u_t\}=\{u_s\dvtx s\le t\}$}.
Furthermore, add an edge $u_ti$ to the graph for each vertex
$i\notin\cZ(t)$ such that $\iit=1$.
Finally, define $T$ by (\ref{t1}) or
(\ref{t2}).

It is easy to see that this constructs a random graph $G_{n,p}$
and that $\cA(t)$, $t\le T$, is as above for this
graph, so the final active set of the bootstrap percolation on the
graph is
$\cA(T)$.

Define also, for $i\in V_n\setminus\cao$,
%
%
\begin{equation}
\label{yi}
Y_i:= \min\{t\dvtx M_i(t) \ge r\}.
\end{equation}
If $Y_i\le T$, then $Y_i$ is the time vertex $i$ becomes active,
but if $Y_i>T$, then $i$ never becomes active. Thus, for $t\le T$,
%
%
\begin{equation}
\label{at2}
\cA(t)=\cao\cup\{i\notin\cao\dvtx Y_i\le t\}.
\end{equation}

By (\ref{mi}), each $M_i(t)$ has a binomial distribution
$\Bin(t,p)$. Further, by~(\ref{mi}) and~(\ref{yi}), each $Y_i$ has a
negative binomial distribution $\NBi(r,p)$,
%
%
\begin{equation}\label{yik}
{\mathbb P}(Y_i=k)={\mathbb P}\bigl(M_i(k-1)=r-1, I_i(k)=1\bigr)
=\pmatrix{k-1\cr r-1}p^r(1-p)^{k-r};\hspace*{-28pt}
\end{equation}
moreover, these random variables $Y_i$ are i.i.d.

We let, for $t=0,1,2,\ldots,$
%
%
\begin{equation}\label{st}
S(t):=|\{i\notin\cao\dvtx Y_i\le t\}|
= \sum_{i\notin\cao} \ett\{Y_i\le t\},\vadjust{\goodbreak}
\end{equation}
so, by (\ref{at2}), and our notation $A(0)=a$,
%
%
\begin{equation}\label{as}
A(t)=A(0)+S(t)=S(t)+a.
\end{equation}

By (\ref{as}), (\ref{t2}) and (\ref{at}), it suffices to study the
stochastic process $S(t)$.
Note that $S(t)$ is a sum of $n-a$ i.i.d.
processes $\ett\{t\ge Y_i\}$, each of which is
$0/1$-valued and jumps from 0 to 1 at time $Y_i$, where $Y_i$ has
the distribution $\NBi(r,p)$ in~(\ref{yik}).
We write $S(t)=S_{n-a}(t)$ when we want to emphasize the
number of summands in $S(t)$; more generally we define
$S_m(t):=\sum_{i=1}^m\ett\{Y_i\le t\}$ for any $m\le n$ [assuming for
consistency that
$\cao=\{n-a+1,\ldots,n\}$].

The fact that $S(t)$, and thus $A(t)$, is a sum of
i.i.d. processes makes the analysis easy; in particular, for
any given $t$,
%
%
\begin{equation}\label{sbin}
S(t)\in\Bin\bigl(n-a,\pi(t)\bigr),
\end{equation}
where
%
%
\begin{equation} \label{pi}
\pi(t)
:={\mathbb P}(\yix\le t)
={\mathbb P}\bigl(\mix(t)\ge r\bigr)
={\mathbb P}\bigl(\Bin(t,p)\ge r\bigr).
\end{equation}
In particular, we have
%
%
\begin{eqnarray}
\label{es}
\E S(t)&=&(n-a)\pi(t),
\\
\label{svar}
\Var S(t)&=&(n-a)\pi(t)\bigl(1-\pi(t)\bigr) \le\E S(t) \le n\pi(t).
\end{eqnarray}

To avoid rounding to integers sometimes below, we define
$S(t):=S(\lfloor t\rfloor)$ and $\pi(t):=\pi(\lfloor t\rfloor)$
for all real $t\ge0$. We also sometimes (when it is obviously harmless)
ignore rounding to simplify notation.

\section{Main results}\label{sec3}\label{Smain}

\subsection{Limits in probability}\label{sec31}\label{SSprob}

For given $r$, $n$, and $p$ we define,
for reasons that will be seen later,
%
%
\begin{eqnarray}
\label{tc}
\tc&:=&\biggl(\frac{(r-1)!}{np^r}\biggr)^{1/(r-1)},
\\
\label{ac}
\ac&:=&\biggl(1-\frac1r\biggr)\tc,
\\
\label{bc}
\bc&:=&n\frac{(pn)^{r-1}}{(r-1)!} e^{-pn}.
\end{eqnarray}
In particular, for $r=2$, $\tc:=1/(np^2)$ and $\ac:=1/(2np^2)$.
For future use, note also that (\ref{tc}) can be written
%
%
\begin{equation}\label{tc1}
n\frac{(p\tc)^r}{r!} = \frac{\tc}{r}.
\end{equation}

Our standard assumptions
$p\ll n^{-1/r}$ and
$p\gg n\qw$
imply that
%
%
\begin{eqnarray} \label{tccond}
\tc&\to&\infty,\qquad p\tc\to0,\qquad \tc/n\to0,\nonumber\\[-8pt]\\[-8pt]
\ac&\to&\infty,\qquad \ac/n\to0,\qquad
\bc/n\to0,\qquad
p\bc\to0,\nonumber
\end{eqnarray}
and further
%
%
\begin{eqnarray}
\label{tc2}
\E S_n(\tc) &=& n\pi(\tc) \sim n\frac{(p\tc)^r}{r!} = \frac{\tc}{r},
\\
\label{bc2a}
n-\E S_n(n)
&=&n\bigl(1-\pi(n)\bigr)=
n{\mathbb P}\bigl(\Bin(n,p)\le r-1\bigr)
\nonumber\\
&\sim& n{\mathbb P}\bigl(\Bin(n,p)= r-1\bigr)
\\
&\sim&
\bcq:=
n\frac{(pn)^{r-1}}{(r-1)!} (1-p)^n.
\nonumber
\end{eqnarray}
If $p\ll n\qqw$, then $(1-p)^n\sim e^{-np}$ and (\ref{bc2a}) yields
%
%
\begin{equation}\label{bc2b}
n-\E S_n(n) =n\bigl(1-\pi(n)\bigr) \sim\bcq \sim\bc;
\end{equation}
if $p$ is larger [$p=\Omega(n\qqw)$, i.e., $n\qqw=O(p)$],
this is not quite true, but in this
case both $\bcq$ and $\bc$ decrease to 0 very fast; in all cases
%
%
\begin{equation}\label{bc2c}
n-\E S_n(n) =n\bigl(1-\pi(n)\bigr) = \bc+o(\bc+1).
\end{equation}

Recall that our main interest is in $S(t)=S_{n-a}(t)$ rather than $S_n(t)$;
see~(\ref{sbin}); for $S(t)$ we obviously have similar results, with
additional error terms depending on $a$; see (\ref{es}) and, for example,
(\ref{esapp}).

Note further that by (\ref{bc}), for any $\beta\in(-\infty,\infty)$,
\[
np-\bigl(\log n+(r-1)\log(np)\bigr) \to
\cases{-\infty, \cr
\beta, \cr
\infty,}
\quad\iff\quad
\bc\to
\cases{\infty, \cr
(r-1)!\qw e^{-\beta}, \cr
0,}
\]
which by simple calculations yields,
provided $p\ge n\qw$,
%
%
\begin{eqnarray}\label{bclim}
&&
np-\bigl(\log n+(r-1)\log\log n\bigr) \to
\cases{-\infty, \cr
\beta, \cr
\infty,}\nonumber\\[-8pt]\\[-8pt]
&&\quad\iff\quad
\bc\to
\cases{ \infty, \cr
(r-1)!\qw e^{-\beta}, \cr
0.}\nonumber
\end{eqnarray}

Our first result, to be refined later, shows that the threshold for
almost percolation is $a=\ac$.
The proof of the theorems in this section are given later
(Sections~\ref{Sthreshold}--\ref{Sgen}).
Let us recall that $A^*$ is the final size of the active set, and that
$A^*=T=A(T)=a+S_{n-a}(T)$.
%
%
\begin{theorem}\label{T1}
Suppose that $r\ge2$ and $n\qw\ll p\ll n^{-1/r}$.
{
\renewcommand\thelonglist{(\roman{longlist})}
\renewcommand\labellonglist{\thelonglist}
\begin{longlist}
\item\label{T1sub}
If $a/\ac\to\ga<1$, then
$\Ax=(\gf(\ga)+o_{\mathrm p}(1))\tc$, where
$\gf(\ga)$ is the unique root in $\oi$ of
%
%
\begin{equation}\label{ika2}
r\gf(\ga)-\gf(\ga)^r=(r-1)\ga.
\end{equation}
[For $r=2$, $\gf(\ga)=1-\sqrt{1-\ga}$.]

Further, $\Ax/a\pto\gf_1(\ga):=\frac{r}{r-1}\gf(\ga)/\ga$, with
$\gf_1(0):=1$.
\item\label{T1super}
If $a/\ac\ge1+\delta$, for some $\delta>0$, then
$\Ax=n-o_{\mathrm p}(n)$;
in other words, we have w.h.p. almost percolation.
More precisely, $\Ax=n-\Op(\bc)$.
\item\label{T1complete}
In case~\ref{T1super}, if further $a\le n/2$, say,
we further have complete percolation, that is,
$\Ax=n$ w.h.p., if and only if $\bc\to0$,
that is, if and only if
$np-(\log n+(r-1)\log\log n) \to\infty$.
\end{longlist}
}
\end{theorem}

It is easily verified that $\gf_1$ is a continuous, strictly
increasing function $[0,1]\to[1,r/(r-1)]$. In particular, in
the subcritical case
\ref{T1sub},
we thus have w.h.p. $\Ax<(r/(r-1))a\le2a$, so
the activation will not spread to many more than the originally
active nodes.

In the supercritical case~\ref{T1super}, we have the
following more detailed result.
%
%
\begin{theorem}\label{T2}
Suppose that $r\ge2$, $n\qw\ll p\ll n^{-1/r}$ and $a=o(n)$, and that
$\Ax=n-o_{\mathrm p}(n)$ as, for example, in Theorem~\ref{T1}\ref{T1super}.
Then:
{
\renewcommand\thelonglist{(\roman{longlist})}
\renewcommand\labellonglist{\thelonglist}
\begin{longlist}
\item\label{T2oo}
If
$np-(\log n+(r-1)\log\log n) \to-\infty$, so
$\bc\to\infty$ by (\ref{bclim}), then
$\Ax=n-\bc(1+o_{\mathrm p}(1))$. In particular, w.h.p.
we do not have complete percolation.
\item\label{T20}
If
$np-(\log n+(r-1)\log\log n) \to\infty$, so
$\bc\to0$ by (\ref{bclim}), then w.h.p.
$\Ax=n$, so we have complete percolation.
\item\label{T2b}
If
$np-(\log n+(r-1)\log\log n) \to\beta\!\in\!(-\infty,\infty)$, so
$\bc\to\bco>0$ by~(\ref{bclim}), then
$n-\Ax\dto\Po(\bco)$; in particular, ${\mathbb P}(\Ax=n)\to\exp
(-\bco
)\in(0,1)$.
\end{longlist}
}

More generally, even if we do not have almost percolation w.h.p.,
the result holds w.h.p. provided $\Ax\ge3\tc$.
\end{theorem}

By the last statement we mean that
${\mathbb P}(\mbox{the result fails and }\Ax\ge3\tc)\to0$.
In particular,
it holds w.h.p. conditioned on
$\Ax\ge3\tc$, provided we have
$\liminf{\mathbb P}(\Ax\ge3\tc)>0$.
%
%
\begin{remark}\label{Rdegrees}
Let $\cB$ be the set of vertices in $G_{n,p}$ with
degrees less than~$r$. These are never activated unless they happen to
be in the initially active set $\cao$, and for each of the vertices,
this has probability $a/n\to0$ if $a=o(n)$; hence trivially
$\Ax\le n-|\cB|(1-o_{\mathrm p}(1))$. We have [cf. (\ref{bc2a}) and~(\ref{bc2c})]
\[
\E|\cB|=n{\mathbb P}\bigl(\Bin(n-1,p)\le r-1\bigr)\sim\bc+o(\bc+1)
\]
with concentration of $|\cB|$ around its mean if $\bc\to\infty$ and
a limiting Poisson distribution
if $\bc\to b<\infty$; see
\cite{JLR}, Sections 6.2 and 6.3, and \cite{SJI}.
Comparing this with Theorem~\ref{T2} we see that in the
supercritical case, and with $a=o(n)$, the final inactive set
$V_n\setminus\cA(T)$ differs from $\cB$ by
$o_{\mathrm p}(|\cB|)$ vertices only, and in the case
$\bc=O(1)$ [combining cases~\ref{T20} and
\ref{T2b} in Theorem~\ref{T2}], w.h.p.
$V_n\setminus\cA(T)=\cB$.
In other words, when we get a large active set,
the vertices that remain inactive are mainly the ones with
degrees less than $r$, and if further $\bc=O(1)$, they are w.h.p. exactly
the vertices with degrees less than $r$.
\end{remark}

We can, as discussed earlier, also consider thresholds for~$p$
for a given $a$.
%
%
\begin{theorem}\label{Tpc}
Suppose that $r\ge2$ and that $a\to\infty$ with
$a=o(n)$. Then the threshold for $p$ for almost
percolation is
%
%
\begin{equation}\label{pc}
\pc:=\biggl(\frac{(r-1)^{r-1}(r-1)!}{r^{r-1}}\biggr)^{1/r}(n a^{r-1})^{-1/r}
\end{equation}
in the sense that if, for some $\delta>0$,
$p\le(1-\delta)\pc$, then $\Ax\le2a=o(n)$
w.h.p., while if $p\ge(1+\delta)\pc$, then
$\Ax=n-o(n)$ w.h.p. In the latter case, further $\Ax=n$ w.h.p.
if and only if $p= (\log n+(r-1)\log\log n+\go(n))/n$ for some
$\go(n)\to\infty$.
\end{theorem}

Note that $n\qw\ll\pc\ll n^{-1/r}$.
Equation (\ref{pc}) is the inverse to (\ref{ac}) in the sense that
the functions $a\mapsto\pc$ and $p\mapsto\ac$
that they define are the inverses
of each other.
For $r=2$, (\ref{pc}) simplifies to $\pc=(2na)\qqw$.
%
%
\begin{remark}
Note that the thresholds for complete and
almost percolation are different only for
large $a$. Indeed, for such a case
the threshold~$\pc$
for almost percolation
can be so small that the graph $G(n,\pc)$
may not be even connected. Then, besides
$\pc$, we have the second threshold
for the complete percolation; for example, if
$a=n/\log n$ and
$r=2$, there are two thresholds:
$\pc=\Theta(\sqrt{\log n}/n)$
for almost percolation,
and $\Theta(\log(n)/n)$ for complete
percolation.
If $a$ is small enough
so that $G(n,\pc)$ is dense enough
(e.g., if $a\le0.49 n/\log^2n$ when $r=2$),
these two thresholds coincide.
\end{remark}

\subsection{Gaussian limits}\label{sec32}\label{SSGauss}

To study the threshold at $\ac$ more precisely, we approximate
$\pi(t)$ in (\ref{pi}) by the corresponding Poisson probability,
%
%
\begin{equation}
\label{tpi}
\tpi(t):={\mathbb P}\bigl(\Po(tp)\ge r\bigr)
=\psi(tp)
:=\sum_{j=r}^\infty\frac{(pt)^j}{j!}e^{-pt}.
\end{equation}
Note that $\psi$ is a differentiable, increasing function on
$(0,\infty)$, and that
%
%
\begin{equation}\label{tpii}
\ddt\tpi(t)=p\psi'(pt)
=p\frac{(pt)^{r-1}}{(r-1)!} e^{-pt}
=\frac{p^rt^{r-1}}{(r-1)!} e^{-pt}.
\end{equation}
By a standard estimate for Poisson approximation of a binomial
distribution (see, e.g., \cite{SJI}, Theorem 2.M),
%
%
\begin{equation}\label{dtv}
|\pi(t)-\tpi(t)|\le\dtv(\Bin(t,p),\Po(tp))
<p,
\end{equation}
where $\dtv$ denotes the total variation distance. A sharper estimate for
small~$t$ will be given in Lemma~\ref{Ltpi}.

We define, for given $n$ and $p$,
%
%
\begin{equation}
\label{acx}
\acx:=-\min_{t\le3\tc}\frac{n\tpi(t)-t}{1-\tpi(t)},
\end{equation}
and let $\tcx\in[0,3\tc]$ be the point where the minimum is attained.
Under our standard assumptions $n\qw\ll p \ll n^{-1/r}$, for
$t\le3\tc$, when $pt\to0$ by (\ref{tccond}), we have, by
(\ref{tpi}) and (\ref{tc1}),
%
%
\begin{equation}
\label{tpi2}
n\tpi(t)\sim n\frac{(pt)^r}{r!}=\biggl(\frac{t}{\tc}\biggr)^{r}\frac{\tc}r
\end{equation}
and thus $\tpi(t)\to0$
and $1-\tpi(t)\sim1$; it
follows easily that $\acx\sim\ac$ and
$\tcx\sim\tc$.
More precise estimates are given in Lemma~\ref{LD},
where it also is shown that~$\tcx$ is unique (for large~$n$, at least).
Furthermore, by Lemma~\ref{LS1} below, for large~$n$, the
minimum in (\ref{acx}) could as well be taken over $t\le n/2$, say,
since $n\tpi(t)-t\ge0$ for $t\in[3\tc,n/2]$.

The following theorem shows that the precise threshold for $a$ is
$\acx\pm O(\sqrt{\ac})$, with a width of the threshold of the order
$\sqrt{\ac}\sim\sqrt{\acx}$. $\Phi$ denotes the standard normal
distribution function.
Note that Theorem~\ref{T2} applies, provided $a=o(n)$, and provides
more detailed information on $\Ax$ when $\Ax$ is large
[i.e., in (ii) and in (iii) conditioned on, say, $\Ax\ge3\tc$].
%
%
\begin{theorem}\label{Tac}
Suppose that $r\ge2$ and $n\qw\ll p\ll n^{-1/r}$.
{
\renewcommand\thelonglist{(\roman{longlist})}
\renewcommand\labellonglist{\thelonglist}
\begin{longlist}
\item\label{Tac-}
If $(a-\acx)/\sqrt{\ac}\to-\infty$, then
for every $\eps>0$, w.h.p. \mbox{$\Ax\le\tcx\le\tc(1+\eps)$}.
If further $a/\acx\to1$, then $\Ax=(1+o_{\mathrm p}(1))\tc$.
\item\label{Tac+}
If $(a-\acx)/\sqrt{\ac}\to+\infty$, then
$\Ax=n-\Op(\bc)$.
\item\label{Tac0}
If $(a-\acx)/\sqrt{\ac}\to y\in(-\infty,\infty)$, then
for every $\eps>0$ and every \mbox{$\bcx\gg\bc$} with $\bcx=o(n)$,
\begin{eqnarray*}
{\mathbb P}(\Ax>n-\bcx) &\to& \Phi\bigl((r-1)^{1/2} y\bigr),
\\
{\mathbb P}\bigl(\Ax\in[(1-\eps)\tc,(1+\eps)\tc]\bigr) &\to& 1-\Phi
\bigl((r-1)^{1/2} y\bigr).
\end{eqnarray*}
\end{longlist}
}
\end{theorem}

For the corresponding result when we keep $a$ fixed and change $p$, we
define, for given $n$ and $a$,
%
%
\begin{equation}\label{hp}
\hh(p):=
\inf_{t\le n/2}\{(n-a)\tpi(t)-t\}
.
\end{equation}
Since $\tpi$ is an increasing function of $p$, $\hh(p)$ is increasing,
with $\hh(0)=-n/2$ and, provided, for example, $a=o(n)$,
$\hh(1)=o(1)$ [attained at $t=o(1)$]. Given
$a=a(n)\to\infty$ with $a=o(n)$,\vadjust{\goodbreak} there is thus (for large $n$) a unique
$\pcx$ such that
%
%
\begin{equation}\label{pcx}
\hh(\pcx)=-a.
\end{equation}
We will see in Lemma~\ref{Lpcx} that $\pcx\sim\pc$.
It is easily verified that, for large $n$ at least,
$\acx=a\iff\hh(p)=-a$, and thus
$p\mapsto\acx$
and $a\mapsto\pcx$ are the inverses of each other.
%
%
\begin{theorem} \label{Tpcxx}
Suppose $r\ge2$ and $a\to\infty$ with $a=o(n)$.
{
\renewcommand\thelonglist{(\roman{longlist})}
\renewcommand\labellonglist{\thelonglist}
\begin{longlist}
\item
If $(p/\pcx-1)a^{1/2}\to-\infty$,
then
$\Ax\le((r/(r-1)+o_{\mathrm p}(1))a$.
If further $p/\pcx\to1$, then
$\Ax= ((r/(r-1)+o_{\mathrm p}(1))a$.

\item
If $(p/\pcx-1)a^{1/2}\to+\infty$,
then $\Ax=n-o_{\mathrm p}(n)$; if further $np-(\log n+(r-1)\log\log
n)\to\infty$,
then $\Ax=n$ w.h.p.
\item
If $(p/\pcx-1)a^{1/2}\to\gl\in(-\infty,\infty)$, then for every
$\eps>0$,
%
%
\begin{eqnarray}
\label{sam}
{\mathbb P}\bigl(\Ax>(1-\eps)n\bigr) &\to& \Phi\bigl(r(r-1)\qqw\gl\bigr),
\\
\qquad{\mathbb
P}\biggl(\Ax\in\biggl[\biggl(\frac{r}{r-1}-\eps\biggr)a,\biggl(\frac{r}{r-1}+\eps\biggr)a\biggr]\biggr)
&\to& 1-
\Phi\bigl(r(r-1)\qqw\gl\bigr).
\end{eqnarray}
If further $np-(\log n+(r-1)\log\log n)\to\infty$, then (\ref{sam})
can be replaced by
$ {\mathbb P}(\Ax=n) \to\Phi(r(r-1)\qqw\gl)$.
\end{longlist}
}
\end{theorem}

In the subcritical cases in Theorems~\ref{T1}\ref{T1sub} and
\ref{Tac}\ref{Tac-}, we also obtain a~Gaussian limit for the size
of the final active set.
%
%
\begin{theorem}
\label{TGsub}
Suppose $r\ge2$ and $n\qw\ll p\ll n\qwr$.
Let $t_*$ be the smallest positive root of
%
%
\begin{equation}
\label{tx}
(n-a)\tpi(t_*)+a-t_*=0.
\end{equation}

{
\renewcommand\thelonglist{(\roman{longlist})}
\renewcommand\labellonglist{\thelonglist}
\begin{longlist}
\item\label{TGsub1}
If $a/\ac\to\ga\in(0,1)$, then
$t_*\sim\gf(\ga)\tc$
with $\gf(\ga)\in(0,1)$ given by (\ref{ika2}),
and
$\Ax\in\AsN(t_*,\gf_2(\ga)\tc)$,
where
$\gf_2(\ga):=\gf(\ga)^r(1-\gf(\ga)^{r-1})\qww/r$.

\item\label{TGsub2}
If $a/\ac\to1$ and also $(a-\acx)/\sqrt{\ac}\to-\infty$, then
$t_*\sim\tc$, more precisely
%
%
\begin{equation}
\label{tx2}
t_*=\tcx-\bigl(1+o(1)\bigr)\sqrt{\frac{2\tc}{r-1}(\acx-a)}
\end{equation}
and
\[
\Ax\in\AsN\biggl(t_*,\frac{\tc}{2(r-1)^2(1-a/\acx) }\biggr).
\]
\end{longlist}
}
\end{theorem}

%
\begin{remark} \label{RTGsub}
It follows from the proof that
in both cases, for large $n$ at least, $t_*$ is the unique root of
(\ref{tx}) in $[0,\tcx]$. In~\ref{TGsub1}, also $t_*<\tc$, so $t_*$
is the
unique root in $[0,\tc]$.
In~\ref{TGsub2}, this is not always true. By Lemma~\ref{LD}, still
for large
$n$, $\tcx>\tc$ and
$\tcx-\tc\sim p\tc^2/(r-1)$. If, for example, $r=2$ and $p=\log
n/n$, then
$\tcx-\tc\sim n/\log^3n$, while $a=\acx-\sqrt n$ yields
$\tcx-t_*\sim\sqrt2 n^{3/4}/\log n \ll\tcx-\tc$.\vadjust{\goodbreak}
\end{remark}

\subsection{The number of generations}\label{sec33}\label{SSGen}

In the supercritical case
$a-\acx\gg\sqrt{\ac}$, when $\cA(0)$ w.h.p. almost percolates by
Theorem~\ref{Tac}, we have the following asymptotic formula for the
number of
generations until the bootstrap percolation process stops.
%
%
\begin{theorem}
\label{TGensuper}
Suppose that $r\ge2$, $n\qw\ll p \ll n\qwr$ and $a=o(n)$.
Assume than $a-\acx\gg\sqrt{\ac}$ [so that $\cA(0)$ w.h.p. almost
percolates].
Then, w.h.p.,
%
%
\begin{eqnarray}\label{tg}
\tau
&\sim&\frac{\pi\sqrt2}{\sqrt{r-1}}\biggl(\frac{\tc}{a-\acx}\biggr)^{1/2}
+
\frac{1}{\log r}
\biggl(\log\log(np) - \log_+\log\frac a{\ac}\biggr)
\nonumber\\[-8pt]\\[-8pt]
&&{}+
\frac{\log n}{np}
+\Op(1).
\nonumber
\end{eqnarray}
\end{theorem}

This theorem is an immediate consequence of Propositions
\ref{PGensuper0},~\ref{PGensuper+},~\ref{PGen2} and~\ref{PG3} in
Section~\ref{Sgen}.
Moreover, these propositions show that the three terms [excepting the error
term $\Op(1)$] in the formula (\ref{tg})
are the numbers of generations required for three
distinct phases of the evolution: the beginning including (possibly) a
bottleneck when the size is about $\tc$;
a period of doubly exponential
growth; and a final phase where the last vertices are activated.
Note that each of the three terms may be the largest one.
%
%
\begin{example}\label{Egen}
Let $p=n^{-\ga}$, with $1/r<\ga<1$, and suppose $a=O(\ac)$.
Then the third term in (\ref{tg}) is $O(1)$ and can be
ignored while the second term is $\log\log n/\log r+O(1)$.
If we are safely supercritical, say $a=2\ac$, then the first term too is
$O(1)$ and the result is $\tau\sim\log\log n/\log r$ w.h.p., dominated
by the
second term.

If instead the process is only barely supercritical, with $a=\acx+\ac
^\gb$
say, with $1/2<\gb<1$, then the first term in (\ref{tg})
is $C n^\gam$ with $C>0$ and the
exponent $\gam=\frac{1-\gb}2\cdot\frac{r\ga-1}{r-1}$, which dominates
the other
terms. Note that the exponent here can be any positive number in $(0,1/4)$
(with $\gam\approx1/4$ if $\ga\approx1$ and $\gb\approx1/2$, so the
graph is
very sparse and the initial set is minimal).

Finally, if $p=\log\log n/n$, say, so the graph is very sparse, and
$a=2\ac$,
then again the first term in (\ref{tg}) is $O(1)$, the second is
$O(\log\log\log\log n)$, while the third is $\log n/\log\log n$,
which thus
dominates the sum.
\end{example}

Note that the second term is $O(\log\log n)$, and the third is $o(\log n)$
[and in many cases $O(1)$ so it can be ignored],
while the first term may be as
large as $n^{1/4-o(1)}$ [although it too in many cases is $O(1)$].
%
%
\begin{remark}
In the subcritical case, one could presumably obtain similar results
for the number of generations until the process stops, but we have not
pursued this topic here.\vadjust{\goodbreak}
\end{remark}

\section{Dynamical models}\label{sec4}\label{Sdyn}

We usually assume, as above, that $a$ and $p$ are given, but
we can also consider dynamical models where one of them grows
with time.

\subsection{Adding external activations}\label{sec41}\label{SSdyna}
In the first dynamical model, we let $n$ and $p$ be given, and
consider a realization of $G_{n,p}$. We start with all vertices
inactive (and
completely uninfected).
We then activate the vertices (from the outside) one by one, in random
order. After each external activation, the bootstrap percolation mechanism
works
as before, activating all vertices that have at least $r$ active
neighbors until no such vertices remain; this is done instantaneously
(or very rapidly) so that this is completed before the next external
activation.
Let $\aox$ be the number of externally activated vertices the first
time that the active set $\cA$ is ``big'' in some sense. For example,
for definiteness, we may define ``big'' as $|\cA|>n/2$.
[It follows from Theorem~\ref{T1} that any threshold $|\cA|>cn$ for a constant
$c\in(0,1)$ will give the same asymptotic results, as well as
thresholds tending to 0 or $n$ sufficiently slowly. If $np-(\log
n+(r-1)\log\log n)\to\infty$, we may also choose the condition
$|\cA|=n$, that is, complete percolation $\cA=V_n$.]
Then $\aox$ is a random variable (depending both on the realization of
$G_{n,p}$ and on the order of external activations). In this formulation,
the threshold result in Theorem~\ref{T1} may be stated as follows.
%
%
\begin{theorem}
\label{TD1}
Suppose that $r\ge2$ and $n\qw\ll p\ll n^{-1/r}$.
Then $\aox/\break\ac\pto1$.
\end{theorem}
\begin{pf}
The active set after $a$ external activations is the same as
the final active set $\cA(T)$ in the static model considered in the
rest of
this paper
with these vertices chosen to be active initially. Hence,
for any given $a$, $\aox\le a$ if and only if bootstrap percolation
with $a$ initially active yields a~big final active set.
In particular, if $\delta>0$,
then Theorem~\ref{T1}\ref{T1sub} implies that ${\mathbb P}(\aox\le
(1-\delta
)\ac
)\to0$,
while
Theorem~\ref{T1}\ref{T1super} and~\ref{T1complete} imply that
${\mathbb P}(\aox\le(1+\delta)\ac)\to1$.
\end{pf}

More precisely, Theorem~\ref{Tac} yields a Gaussian limit.
%
%
\begin{theorem}
\label{TDG1}
Suppose that $r\ge2$ and $n\qw\ll p\ll n^{-1/r}$.
Then $\aox\in\AsN(\acx,\ac/(r-1))$.
\end{theorem}
\begin{pf}
Let $x\in\oooo$. Then, arguing as in the proof of Theorem~\ref{TD1}
but now
using Theorem~\ref{Tac}\ref{Tac0} (with $y=x/\sqrt{r-1}$), we find
\[
{\mathbb P}\biggl(\frac{\aox-\acx}{\sqrt{\ac/(r-1)}}\le x\biggr)
=
{\mathbb P}\bigl(\aox\le\acx+ x\sqrt{\ac/(r-1)} \bigr)
\to\Phi(x).
\]
\upqed
\end{pf}

We have here for simplicity assumed that the external activations are
done by sampling without replacement, but otherwise independently of
whether the vertices already are (internally) activated.
A natural variation is to only activate vertices that are inactive.
Let $\aoxx$ be
the number of externally activated vertices when the active set
becomes big in this version.
Since a new activation of an already active vertex does not matter at
all, $\aoxx$ equals in the version above the number of externally
active vertices among the first~$\aox$ that are not already
internally activated. Thus $\aox-\aoxx$ is the number of external
activations that hit an already active vertex. It is easily verified
that this is $o_{\mathrm p}(\ac)$, and thus Theorem~\ref{TD1} holds
for $\aoxx$ as
well; we omit the details.
It seems likely that it is possible to
derive a version of the Gaussian limit in Theorem~\ref{TDG1} for
$\aoxx
$ too,
but that would require a more careful estimate of $\aox-\aoxx$ (and in
particular its variance), which we have not done, so we leave this
possibility as an open problem.
%
%
\begin{remark}\label{Rdyna}
One way to think about this dynamical model,
where we add new active vertices successively and may
think of these as being initially active,
is to see it as a sequence of
bootstrap percolation processes, one for each $a=0,1,\ldots,n$; the processes
live on the same graph $G_{n,p}$ but have different numbers of
initially active
vertices, and they are coupled in a~natural way.
In order to really have the same realization of $G_{n,p}$ for different~$a$, we
have to be careful in the choice of the order in which we explore the vertex
neighborhoods, that is, the choice of $u_t$. [Recall that $G_{n,p}$ is
constructed from the indicators $I_i(s)$ and the sequence $(u_t)$; see
Section~\ref{Ssetup}.] We can achieve this by first making a list $L$
of all
vertices in the (random) order in which they are externally activated.
We then
at each time $t$ choose~$u_t$ as an unused internally activated vertex
(e.g., the most recent one) if there is any such vertex, and otherwise as
the next unused vertex in the list~$L$.

This model makes it possible to pose now questions about the
bootstrap percolation. For example, we may consider the critical process
starting with exactly $\aox$ initially active vertices (i.e., the
first process
that grows beyond the bottleneck and becomes big)
and ask for the
number of generations until the process dies out.
Alternatively, we may consider the process
starting with exactly $\aox-1$ initially active vertices (i.e., the last
process that does not become big) and ask for its final
size.

Such questions will not be treated in the present paper,
but we mention that it
is easily seen that the final size with $\aox-1$ initially active vertices
is $\tc(1+o_{\mathrm p}(1))$ so that the
final size jumps from about $\tc$ to about $n$ with the addition of a single
additional initial vertex.
Furthermore, we conjecture that, under suitable conditions, the
number of generations for the
process with
$\aox$ initially active vertices is of order
$\ac\qqq$
(which is much larger than the number of generations for any fixed $a$; see
Section~\ref{SSGen}).
\end{remark}

\subsection{Adding external infections}\label{sec42}\label{SSdyni}
An alternative to external activations is external infections, where
we again start with all vertices inactive and uninfected, and infect
vertices one by one from the outside, choosing the infected vertices at
random (independently and with replacement); as before, $r$ infections
(external or internal) are
needed for activation, and active vertices infect their neighbors.
Let $J_0$ be the number of external infections when the active set first
becomes ``big'' (as in Section~\ref{SSdyna}). (Thus, $J_0$ is a~random
variable.)

In the original model, each initially active vertex infects about $np$
other vertices so the total number of initial infections is about $npa$;
it is thus easy to guess that $J_0\approx np\ac$.
Indeed, this is the case as is shown by the next theorem.
We cannot
(as far as we know) directly derive this from our previous results,
since the dependencies between infections in the two versions are slightly
different, but it follows by a minor variation of our method; see
Section~\ref{Sthreshold}.
We believe that the result could be sharpened to a Gaussian limit as in
Theorem~\ref{TDG1}, but we leave this to the reader.
%
%
\begin{theorem}
\label{TJ}
Suppose that $r\ge2$ and $n\qw\ll p\ll n^{-1/r}$.
Then $J_0/(np\times\ac)\pto1$.
\end{theorem}

In particular, for $r=2$, we thus have $J_0\psim np\ac=1/(2p)$.

\subsection{Adding edges}\label{sec43}\label{SSdynm}
In the second dynamical model, $n$ and $a$ are given; we start
with $n$ vertices of which $a$ are active, but no edges.
We then add the edges of the complete graph $K_n$ one by one, in
random order. As in the previous dynamical model, bootstrap percolation
takes place instantaneously after each new edge is added.

It is convenient to use the standard method of adding the edges at
random times (as in, e.g., \cite{SJII}).
Thus, each edge $e$ in $K_n$ is added at a time~$U_e$, where~$U_e$ are
independent and uniformly distributed on $\oi$. Then, at a~time
$u\in\oi$, the resulting graph is $G_{n,p}$ with $p=u$.
(We use $u$ to denote this time variable, in order not to confuse it
with the time $t$ used to describe the bootstrap percolation process.)

Let the random variable $M$ be the number of edges required to obtain
a~big active set $\cA$, where ``big'' is defined as in Section~\ref{SSdyna}.
%
%
\begin{theorem} \label{TDG2}
Suppose $r\ge2$ and $a\to\infty$ with $a=o(n)$.
Then
\[
M=\pmatrix{n\cr2}\pc\bigl(1+o_{\mathrm p}(1)\bigr)=\frac12n^2\pc
\bigl(1+o_{\mathrm p}(1)\bigr).
\]
More precisely,
\[
M\in\AsN\biggl(\pmatrix{n\cr2}\pcx,\frac{r-1}{4r^2}\frac{(n^2\pc)^2}{a}\biggr).
\]
\end{theorem}

The proof is given in Section~\ref{SGauss}.
%
%
\begin{remark}[(Coupling different $p$)]\label{Rdynpp}
The proof of Theorem~\ref{TDG2} is based on using our earlier results
for a single
$p$.
We might also want to study the bootstrap percolation process for all $p$
at once [or equivalently, in $\gnm$ for all $m$ at once], that is,
with a
coupling of the models for different $p$, for given $n$ and $a$.
As in Remark~\ref{Rdyna}, this requires a careful choice of the order in
which the
vertices are inspected.
We can achieve this by modifying the formulation in Section \ref
{Ssetup} as follows:

When we have chosen a vertex $u_t$, we reveal the times $U_{tj}$ that
the edges from it appear; this tells us the neighborhood of $u_t$ at
any time $u$. We begin by choosing $u_1,\ldots,u_a$ as the initially
active vertices. We then, after each choice of $u_t$, $t\ge a$,
calculate for each of the remaining $n-t$ vertices the time when it
acquires the $r$th edge to $\{u_1,\ldots,u_t\}$, and let $u_{t+1}$ be
the vertex such that this time is minimal. Then, fixing any time
$u=p$, the chosen vertices $u_t$ will all be active until the first
time that
no unused active vertices remain, and the process stops.
In this manner, we have found a choice of $u_1,u_2,\ldots$ that
satisfies the description in Section~\ref{Ssetup} for all $p\in\oi$
simultaneously.

As in Remark~\ref{Rdyna}, we can use this model to study, for example,
the last
``small'' or the first ``big'' bootstrap percolation process, when we
add edges one by one with a given set of initially active vertices.
Again, we will not consider such questions in the present paper.
\end{remark}

\section{Boundary cases}\label{sec5}\label{Sbound}

We have above assumed $r\ge2$ and $n\qw\ll p\ll n^{-1/r}$.
In this section we treat the cases when these assumptions do not
hold.
Proofs are given in Section~\ref{Spf+}.

We begin with the sparse case $p(n) \sim{c/n}$ when $\tc$ and $\ac$
defined by~(\ref{tc}) and (\ref{ac}) are of order $n$. (The exact values are
no longer relevant, since they are based on approximations no longer valid.)
This suggests that the interesting case is when $a\asymp n$, that is,
when a
positive fraction of all vertices are initially active.
Indeed, Theorem~\ref{pn=cn} shows that, w.h.p.,
if we start
with a~positive fraction of the graph, then the activation spreads to a~larger part
of the graph but does not reach almost all vertices; if, on the
contrary, the
size of the original set of activated vertices is negligible with
respect to
the size of the graph, then the activation does not spread to a
positive fraction of the graph.
Provided $c$ is large enough, there is,
as found by Scalia-Tomba \cite{Scalia-Tomba},
also in this case a dichotomy, or ``phase transition,''
similar to Theorem~\ref{T1}, with a sudden jump from
a ``small'' to a ``large'' final active set, although in this case all
sets are
of order $n$ so the jump is not as dramatic as for larger~$p$.

Define, for $x\ge0$, $c\ge0$ and $\gth\in\oi$,
%
%
\begin{eqnarray}
\label{f1a}
f(x,c,\gth) :\!&=& (1 - \gth){\mathbb P}\bigl(\Po(cx) \ge r\bigr) +
\gth- x \nonumber\\[-8pt]\\[-8pt]
&=& (1 - \gth) \sum_{j=r}^\infty\frac
{(cx)^j}{j!}e^{-cx} - x + \gth\nonumber
\\
\label{f1b1}
&=&
1-x-(1-\gth){\mathbb P}\bigl(\Po(cx)\le r-1\bigr)
\\
\label{f1b2}
&=&
1-x-(1-\gth)\sum_{j=0}^{r-1} \frac{(cx)^j}{j!}e^{-cx},
\end{eqnarray}
and let $x_0(\theta)$ be the smallest root $x\ge0$ of
%
%
\begin{equation}
\label{f10}
f(x,c,\gth)=0;
\end{equation}
similarly, let
$\xox(\gth)$ be the largest root in $\oi$ of this equation.

Since $f(0,c,\gth)=\gth\ge0$ and
$f(1,c,\gth)=-(1-\gth){\mathbb P}(\Po(c)\le r-1)\le0$, there is always
at least one
root in $\oi$, and $0\le\xo(\gth)\le\xox(\gth)\le1$; further
$0< \xo(\gth)\le\xox(\gth)<1$ when $0<\gth<1$ while $\xo(0)=0$ and
$\xo(1)=\xox(1)=1$.
We also define
%
%
\begin{equation}
\label{ccr}\quad
\ccr=\ccr(r)
:=r+\frac{{\mathbb P}(\Po(r-1)\le r-2)}{{\mathbb P}(\Po(r-1)= r-1)}
=r+\frac{\sum_{j=0}^{r-2}(r-1)^j/j!}{(r-1)^{r-1}/(r-1)!}.
\end{equation}
Thus $\ccr(2)=3$, $\ccr(3)=9/2$, $\ccr(4)=53/9$.
%
%
\begin{lemma}\label{Lf}
\textup{(i)}
If $0\le c\le\ccr$, then (\ref{f10}) has a unique root $x=x_0(\theta
)\in
\oi$ for
every $\gth\in\oi$, and $x_0(\theta)$ is a continuous strictly increasing
function of $\gth$.

\textup{(ii)}
If $c>\ccr$, then there exists $\gthcq=\gthcq(c)$
and $\gthc=\gthc(c)$ with $0\le\gthcq<\gthc<1$ such that (\ref{f10})
has three roots in $\oi$ when $\gth\in(\gthcq,\gthc)$,
but a unique root when $\gth\in[0,\gthcq)$ or $\gth\in(\gthc,1]$;
if $\gth=\gthcq>0$
or $\gth=\gthc$, there are two roots, one of them double.
The smallest root
$x_0(\theta)$ is strictly increasing and continuous on $\oi$ except at
$\gthc$
where it has a jump from $\xo(\gthc)$ to $\xox(\gthc)>\xo(\gthc
)$, where
$\xox(\gthc)=\xo(\gthc+):=\lim_{\gth\downto\gthc}\xo(\gth)$
is the other
root for $\gth=\gthc$.
Furthermore, if
$\gth=\gthc$, then
$f(x,c,\gth)\ge0$ for
$x\in[0,\xox(\gth)]$, and $x_0(\theta)$ is a double root.
\end{lemma}
%
%
\begin{theorem}\label{pn=cn}
Suppose that $r\ge2$, $p\sim c/n$ and $a\sim\theta n$ for some constants
$c\ge0$ and $\theta\ge0$.

{
\renewcommand\thelonglist{(\roman{longlist})}
\renewcommand\labellonglist{\thelonglist}
\begin{longlist}
\item\label{p=cnaon}
If $\gth=0$, that is, if $a=o(n)$, then $\Ax/a\pto1$.
\item\label{pn=cn=0}
If $c=0$, that is, if $p=o(1/n)$, then $\Ax/a\pto1$.
\item\label{pn=cnsub}
If $0\le c\le\ccr$, then $\Ax/n\pto x_0(\theta)$, where $x_0(\theta)$
is the unique
nonnegative root of (\ref{f10}).
\item\label{pn=cnsuper}
If $c>\ccr$ and $\gth\neq\gthc(c)$ given by Lemma~\ref{Lf},
then $\Ax/n\pto x_0(\theta)$, where $x_0(\theta)$ is the smallest nonnegative
root of (\ref{f10}).
\end{longlist}
}
\end{theorem}

There is thus a jump in the final size at $a=\gthc n$.
Remark~\ref{Rgthc} shows how to find~$\gthc$.\vadjust{\goodbreak}
%
%
\begin{remark}\label{Rgthcc}
$\gthc(c)$ and $\gthcq(c)$ are decreasing functions of $c$.
[$\gthc(c)$ is strictly decreasing, while $\gthcq(c)$ is constant $0$
for large
$c$.]
Hence their largest value is, by the calculation in Remark~\ref{Rgthc},
\begin{eqnarray*}
\gthcc
&=&\gthcc(r)
:= \gthc(\ccr)
=\gthcq(\ccr)
\\
&=&1-\frac{1}{r{\mathbb P}(\Po(r-1)=r-1)+{\mathbb P}(\Po(r-1)\le r-2)}.
\end{eqnarray*}
Thus $\gthcc(2)=1-e/3$, $\gthcc(3)=1-e^2/9$, $\gthcc(4)=1-2e^3/53$.

The threshold for $\gthcq=0$ can be calculated too.
For $r=2$, $\gthcq(c)=0$ for
$c\ge e^y/y$, where $e^y=1+y+y^2$; numerically, this is
$c\ge3.35091\ldots.$
\end{remark}
%
%
\begin{remark}
We have here considered a given $p\sim c/n$ and varied $a\sim\gth n$.
If we instead, as in Theorem~\ref{Tpc}, take a given $a\sim\gth n$ for
a fixed
$\gth$ and vary $c=pn$, we have a similar phenomenon. Lemma~\ref{Lf} and
Theorem~\ref{pn=cn} apply
for every combination of $\gth$ and $c$, and by considering the set of
$(c,\gth)\in\bbR_+^2$ such that (\ref{f10}) has two or three roots,
it
follows from Remark~\ref{Rgthcc} that
if $\gth\ge\gthcc$, then $\Ax/n\pto\xo(c)$, where $\xo(c)$ is
the unique
root of (\ref{f10}) and thus a continuous function of~$c$, while if
$\gth<\gthcc$, then there is a range of $c$ where~(\ref{f10}) has three
roots, and one value of $c$ where the limit value $\xo(c)$ jumps from a
``small'' to a ``large'' value. Thus there is, again, a kind of phase
transition.
\end{remark}

The following theorem shows that if we for simplicity take $p=c/n$,
then the precise threshold for $a$ in Theorem~\ref{pn=cn}(iv) is $
\gthc n \pm O(\sqrt{n})$, with a width of the threshold of the order
$\sqrt{n}$.
%
%
\begin{theorem}\label{cngauss}
Suppose that $r\ge2$ and $p = \frac{c}{n}$ with $c>\ccr$ fixed. Let
$\gthc=\gthc(c)$, $x_0=x_0(\gthc)$ and $\xox=\xox(\gthc)$ be as in
Lemma~\ref{Lf};
thus $x_0$ and $\xox$ are the two roots in $\oi$ of $f(x,c,\gthc
)=0$, with
$x_0<\xox$.

{
\renewcommand\thelonglist{(\roman{longlist})}
\renewcommand\labellonglist{\thelonglist}
\begin{longlist}
\item\label{cngauss-}
If $(a-\gthc n)/\sqrt{n}\to-\infty$, then for any $\eps>0$, w.h.p.
$\Ax
\le(1+\eps)x_0n$. If further $a\sim\gthc n$, then $\Ax
=(1+o_{\mathrm p}(1))x_0n$.
\item\label{cngauss+}
If $(a-\gthc n)/\sqrt{n}\to+\infty$, then for any $\eps>0$, w.h.p.
$\Ax\geq(1-\eps)x_1n$.
If further $a\sim\gthc n$, then $\Ax=(1+o_{\mathrm p}(1))\xox n$.
\item\label{cngauss0}
If $(a-\gthc n)/\sqrt{n}\to y \in(-\infty,\infty)$,
then there exists a sequence $\eps_n\to0$ such that
\begin{eqnarray*}
{\mathbb P}\bigl(\Ax\in[(1-\eps_n)x_1n,(1+\eps_n)x_1n]\bigr) &\to&
\Phi(y/\gs),
\\
{\mathbb P}\bigl(\Ax\in[(1-\eps_n)x_0n,(1+\eps_n)x_0n]\bigr) &\to&
1-\Phi(y/\gs),
\end{eqnarray*}
where $\gss=(1-\gthc)\psi(c\xo)/(1-\psi(c\xo))>0$.
\end{longlist}
}
\end{theorem}

At the other, dense, endpoint of our range we have
$p(n) \sim c n^{-1/r}$. Then $\tc$ and $\ac$ in (\ref{tc}) and (\ref{ac})
are of order constant.
(Again the exact values are irrelevant.)
This suggests, and the following theorem makes more precise,
that the process\vadjust{\goodbreak} will either die out or grow very quickly,
with the outcome determined by the first few steps,
and that the activation can spread from a set of constant size to the
entire graph with a positive probability, which, however, is bounded
away from
$1$.
%
%
\begin{theorem}\label{sqrtn}
Suppose $r \geq2$ and $p\sim c n^{-1/r}$ for a constant $c>0$.
{
\renewcommand\thelonglist{(\roman{longlist})}
\renewcommand\labellonglist{\thelonglist}
\begin{longlist}
\item\label{sqrtn1}
If $a$ is fixed with $a \geq r$, then
%
%
\begin{equation}\label{eqsqrtn}
{\mathbb P} ( \Ax= n ) \rightarrow\zeta(a,c)
\end{equation}
for some $\zeta(a,c)\in(0,1)$.
Furthermore, there exist numbers $\zeta(a,c,k)>0$ for $k\ge a$ such that
${\mathbb P}(\Ax=k)\to\zeta(a,c,k)$ for each fixed $k\ge a$, and
$\sum_{k=a}^\infty\zeta(a,c,\break k)+\zeta(a,c)=1$.
\item
If $a \to\infty$, then ${\mathbb P}(\Ax=n)\to1$, that is, $\Ax=n$ w.h.p.
\end{longlist}
}
\end{theorem}
%
%
\begin{remark}
\label{RMS}
The limiting probabilities in Theorem~\ref{sqrtn} can be expressed as hitting
probabilities of an inhomogeneous random walk.
Let $\xi_k\in\Po\bigl({k-1\choose r-1}c^r\bigr)$, $k\ge1$, be
independent, and
let
$\tsk:=\sum_{j=1}^k(\xi_j-1)$ and
$\tT:=\min\{k\dvtx a+\tsk=0\}\in\bbN\cup\{\infty\}$.
Then
%
%
\begin{equation}
\zeta(a,c)={\mathbb P}(\tT=\infty)={\mathbb P}(a+\tsk\ge1 \mbox{ for
all $k\ge1$})
\end{equation}
and $\zeta(a,c,k)={\mathbb P}(\tT=k)$.
Consequently, Theorem~\ref{sqrtn}\ref{sqrtn1} can also be written as
\[
\dtv(\Ax,\min(\tT,n))\to0,
\]
where $\dtv$ is the total variation distance.
\end{remark}

If the probability of connections $p$ is even larger,
$p \gg n^{-1/r}$,
then the initial set percolates as long as $a \geq r$.
%
%
\begin{theorem}\label{sqrtnn}
Let $r\ge2$. If $p \gg n^{-1/r}$ and $a \geq r$, then $\Ax=n$ w.h.p.
\end{theorem}
%
%
\begin{remark}
\label{R1}
The case $r=1$ is different. In this case,
infection is equivalent to activation, and spreads to
every vertex connected to an active vertex.
Thus the final active set is the union of the components of the graph
that contain at least one initially active vertex.
It is well known that this is equivalent to the Reed--Frost model for
epidemics, where each infected person infects everyone else with probability
$p$, all infections being independent. (This equivalence is easily seen by
the argument in Section~\ref{Ssetup}.)
The Reed--Frost model
has been much
studied; see, for example, von Bahr and Martin-L\"of \cite{vBahrML},
Martin-L\"of
\cite{ML86,ML98}.
We state some known result for comparison with our results for $r\ge2$;
proofs can be found in \cite{vBahrML,ML86}, where also further details are
given
(including central limit theorems as in Section~\ref{SSGauss}),
or by modifying the proofs of the results above.
Many results follow also
easily from known results on the component structure of
$G_{n,p}$.\vadjust{\goodbreak}

If $p=\log n/n + \go(n)/n$, with $\go(n)\to\infty$, then w.h.p. $G_{n,p}$
is connected and thus $\Ax=n$ as soon as $a\neq0$.
More generally, if $p\gg n\qw$ and $a\ge1$, then w.h.p. $\Ax=n-o(n)$;
cf. Theorem~\ref{sqrtnn}.

The case $p=c/n$ is perhaps more interesting. There are many (w.h.p.
$\ge
c'n$) isolated
vertices, so we cannot have percolation or almost percolation unless
$a/n\to1$. If $c>1$, there is a single giant component of size
$\rho n+o_{\mathrm p}(n)$, with $\rho=\rho(c)>0$, and thus, if $a\ge1$
is fixed, then there is a dichotomy, with either $\Ax=o(n)$
or $\Ax=\rho n+o(n)$ w.h.p., with probabilities converging to the positive
$(1-\rho)^a$ and $1-(1-\rho)^a$, respectively; cf. Theorem~\ref{sqrtn}.

If $c\le1$ and $a$ is fixed, then, by the same argument,
$\Ax$ converges to the total size of a Galton--Watson process with
$\Po(c)$
offspring distribution and $a$ initial individuals
(a Borel--Tanner distribution). Thus $\Ax/a$ is stochastically bounded but
does not converge in probability to a constant; cf. Theorem~\ref{pn=cn}\ref{p=cnaon}.

If $c<1$ and $a\to\infty$ but $a=o(n)$, then
$\Ax/a\pto1/(1-c)$; cf. Theorem~\ref{pn=cn}\ref{p=cnaon}.

If $p\sim c/n$ with any $c>0$ and $a\sim\gth n$ with $\gth>0$, then
$\Ax/n\pto\xo(\gth)$, where $\xo$ is the unique positive root of
(\ref{f10}), where now $f(x,c,\gth)=1-x-(1-\gth)e^{-cx}$ by (\ref
{f1b2}) with
$r=1$. In other words, Theorem~\ref{pn=cn}\ref{pn=cnsub} holds for
$r=1$, too, with
$\ccr=\infty$,
and there is no threshold. [For $\gth=0$ there is the well-known threshold
at $c=1$, but note that also in this case, $\xo$ is continuous in both $c$
and $\gth$ and there is no jump as in Theorem
\ref{pn=cn}\ref{pn=cnsuper}.]
\end{remark}

\section{Overview of the proofs}\label{sec6}\label{Soverview}

By (\ref{at}), (\ref{t2}) and (\ref{as}), for
$u=1,2,\ldots$
%
%
\begin{eqnarray}
\label{a1}
\Ax\ge u \quad&\iff&\quad T\ge u \quad\iff\quad\min_{t<u}\bigl(A(t)-t\bigr)
>0\nonumber\\[-8pt]\\[-8pt]
\quad&\iff&\quad a+\min_{t<u}\bigl(S(t)-t\bigr) >0.\nonumber
\end{eqnarray}
Hence, $\Ax=T$ is completely determined by the stochastic
process\break
$\min_{t<u}(S(t)-t)$, $u\ge1$.
In particular, $\cao$ percolates if and only if
$a>-\min_{t<n}(S(t)-t)=\max_{t<n}(t-S(t))$.

Note that (\ref{a1}) is an exact representation of
$\Ax$; we have not yet made any approximations. To obtain
asymptotic results, we introduce some simple approximations. We give
an informal overview of the argument here; details will follow in
later sections.

First, $S(t)\approx\E S(t)$ by the law of large numbers. A
simple calculation will show that $f(t):=\E S(t)-t$ starts
at 0 for $t=0$, then decreases to a~minimum at
$t\approx\tc$ given in (\ref{tc}), and then
increases until $\E S(t)\approx n$ and thus
$f(t)\approx n-t$; then $f(t)\approx n-t$ holds
until $t=n$, so $f(t)$ decreases again in this range to a
final value $f(n)=\E S(n)-n\approx0$.

There are thus two candidates for the minimum point of $S(t)-t$:
either $t\approx\tc$ or $t\approx n$. What
happens at $t\approx n$ makes the difference between almost
percolation and complete percolation; we will study this too in detail
later, but for the moment we ignore it and concentrate on whether we
have almost percolation or not, and we see that, roughly, this is
determined by whether $a>-(S(\tc)-\tc)$ or not, which can
be approximated by \mbox{$a>-(\E S(\tc)-\tc)$}.
A simple calculation yields
$\E S(\tc)-\tc\approx-\ac$, which establishes the threshold at
$\ac$.

This argument also gives the following picture of the course of the
activation $\cA(t)$ in the critical case
$a\approx\ac$. (We leave the modifications in the
subcritical and supercritical cases to the reader.)
We start with $A(0)=a$. At first, there are very few new
vertices that reach the threshold of $r$ infections,
and the number $A(t)-t$ of unused vertices goes down, and
approaches 0 as~$t$ approaches $\tc$. However, the rate
of activation of new vertices is increasing, because a pool of
vertices with $r-1$ infections has accumulated, and as
$t\to\tc$, new vertices are activated at about the same
rate as they are used. There are now two possibilities: either the
activation dies out at this point, with a total size
about $\tc= r/(r-1)a_c$, or it survives this bottleneck, and it then
rapidly grows after time $\tc$ until almost all vertices are
active. In the latter case there are again two possibilities: either
all remaining vertices are finally active (complete percolation), or
a few are not.

\section{Approximation of $S(t)$ by its mean}\label{sec7}

For notational convenience, we assume that
$V_n\setminus\cao=\{1,\ldots,n-a\}$.
Note first that $(n-a)\qw
S(t)=(n-a)\qw\sum_{i=1}^{n-a}\ett\{Y_i\le t\}$ is the empirical distribution
function of $\{Y_i\}_1^{n-a}$.
By the law of large numbers for the binomial distribution
(\ref{sbin}), for every
$t=t(n)$,
$S(t)=(n-a)\pi(t)+o_{\mathrm p}(n)$.
Moreover,
by the Glivenko--Cantelli theorem (\cite{Kallenberg},
Proposition 4.24),
the following holds uniformly for all $t$:
%
%
\begin{lemma}
\label{L1}
$\sup_{t\ge0}|S(t)-\E S(t)|=o_{\mathrm p}(n)$.
\end{lemma}
\begin{pf}
If $n-a\ge\sqrt{n}$, say, this is a weaker version of
\cite{Kallenberg}, Proposition~4.24. For smaller
$n-a$, the result is trivial, since $0\le S(t)\le n-a$.~%
\end{pf}

For small $t$, the uniform error bound in Lemma~\ref{L1} is not good
enough. [It can be improved to $\Op(n^{1/2})$, see Lemma~\ref{L2},
but this too is too large for our purposes.]
For each $t$, (\ref{svar}) gives a bound
$\Op((n\pi(t))^{1/2})$. We extend this to a uniform bound
for a range of $t$ by a martingale argument.
We begin by introducing a pair of well-known martingales for empirical
distribution functions. (See \cite{SJII}, Lemma 2.1, for a continuous
time version.)
%
%
\begin{lemma}
\label{Lmartin}
The stochastic process
%
%
\begin{equation}
\label{lmartina}
\frac{S(t)-\E S(t)}{1-\pi(t)},\qquad t=0,1,\ldots,
\end{equation}
is a martingale, and
%
%
\begin{equation}
\label{lmartinb}
\frac{S(t)-\E S(t)}{\pi(t)},\qquad t=r,r+1,\ldots,
\end{equation}
is a reverse martingale.
\end{lemma}
\begin{pf}
Since $S(t)$ is a sum of $n-a$ i.i.d. processes
$\ett\{Y_i\le t\}$, it suffices to treat each of
these separately, that is, for the first part to show that, for each $i$,
\[
X(t)=X_i(t):=\frac{\ett\{Y_i\le t\}-{\mathbb P}(Y_i\le t)}{1-{\mathbb
P}(Y_i\le t)}
=1-\frac{\ett\{Y_i> t\}}{{\mathbb P}(Y_i> t)}
\]
is a martingale. This is elementary: if $Y_i\le t$, then
$X(t)=X(t+1)=1$.
If $Y_i>t$, then
$X(t)=-\pi(t)/(1-\pi(t))$ either jumps to $X(t+1)=1$ or decreases
to $X(t+1)=-\pi(t+1)/(1-\pi(t+1))$, and the conditional probability of
these events are $(\pi(t+1)-\pi(t))/(1-\pi(t))$ and
$(1-\pi(t+1))/(1-\pi(t))$, respectively, so a simple
calculation yields $\E(X(t+1)\mid Y_i>t)=-\pi(t)/(1-\pi(t))$.
[Alternatively, this follows from the case $X(t)=1$ and the fact
that $\E X(t+1)=\E X(t)$.]
Hence, $\E(X(t+1)\mid X(1),\ldots,X(t))=X(t)$.

For the second part, we similarly find that
$\tX(t):=\ett\{Y_i\le t\}/\pi(t)$ is a~reverse martingale,
that is, that $\E(\tX(t)\mid\tX(t+1),\ldots)=\tX(t+1)$.
\end{pf}
%
%
\begin{lemma}
\label{L2}
For any $t_0$,
%
%
\begin{eqnarray}
\label{l2a}
\E\Bigl(\sup_{t\le t_0} |S(t)-\E S(t)|\Bigr)^2 &\le& 16 n \pi(t_0),
\\
\label{l2b}
\E\Bigl(\sup_{t\ge t_0} |S(t)-\E S(t)|\Bigr)^2 &\le& 16 n \bigl(1-\pi(t_0)\bigr).
\end{eqnarray}
\end{lemma}
\begin{pf}
Assume first $\pi(t_0)\le1/2$.
Let $\xi(t):=(S(t)-\E S(t))/(1-\pi(t))$. This is a
martingale by Lemma~\ref{Lmartin}, and Doob's inequality
(\cite{Kallenberg}, Proposition~7.16) yields, using (\ref{svar}),
%
%
\begin{eqnarray}\label{er8a}
\E\Bigl(\sup_{t\le t_0} |S(t)-\E S(t)|\Bigr)^2
&\le&
\E\sup_{t\le t_0} |\xi(t)|^2
\nonumber\\
&\le& 4 \E|\xi(t_0)|^2
=4\frac{\Var S(t_0)}{(1-\pi(t_0))^2}\\
&\le&
8n\pi(t_0),
\nonumber
\end{eqnarray}
which proves (\ref{l2a}) in this case. Similarly, if $\pi(t_0)\ge1/2$,
then we obtain, using the reverse
martingale (\ref{lmartinb}),
%
%
\begin{equation}\label{er8b}
\E\Bigl(\sup_{t\ge t_0} |S(t)-\E S(t)|\Bigr)^2
\le4\frac{\Var S(t_0)}{\pi(t_0)^2}
\le8n\bigl(1-\pi(t_0)\bigr).
\end{equation}

Now, let $t_1$ be the largest integer such that
$\pi(t_1)\le1/2$. We can apply~(\ref{er8a}) with
$t_0=t_1$ and (\ref{er8b}) with $t_0=t_1+1$,
and thus
\begin{eqnarray*}
\E\Bigl(\sup_{t\ge0} |S(t)-\E S(t)|\Bigr)^2
&\le&
\E\Bigl(\sup_{t\le t_1} |S(t)-\E S(t)|\Bigr)^2\\
&&{}+
\E\Bigl(\sup_{t\ge t_1+1} |S(t)-\E S(t)|\Bigr)^2
\\
&\le& 8n.
\end{eqnarray*}
This immediately implies (\ref{l2a}) for
$\pi(t_0)>1/2$ and (\ref{l2b}) for $\pi(t_0)<1/2$.
\end{pf}

\section{\texorpdfstring{Approximation of $\E S(t)$ and proofs of Theorems \protect\ref{T1}--\protect\ref{Tpc}}
{Approximation of $\E S(t)$ and proofs of Theorems 3.1--3.4}}\label{sec8}
\label{Sthreshold}

For (real) $t>0$ and $pt\le1$, say,
by (\ref{pi}),
%
%
\begin{eqnarray}\label{pia}
\pi(t)&=&\sum_{j=r}^{\lfloor t\rfloor}\pmatrix{{\lfloor t\rfloor}\cr j}
p^j(1-p)^{\lfloor t\rfloor-j}
=\pmatrix{\lfloor t\rfloor\cr r} p^r\bigl(1+O(pt)\bigr)\nonumber\\[-8pt]\\[-8pt]
&=&\frac{t^rp^r}{r!}\bigl(1+O(pt+t\qw)\bigr)\nonumber
\end{eqnarray}
[cf. (\ref{tpi2})],
and thus, by (\ref{es}),
%
%
\begin{equation}\label{esapp}
\E S(t)-n\frac{t^rp^r}{r!}
=O\bigl(nt^rp^r(pt+t\qw+a/n)\bigr).
\end{equation}
It thus\vspace*{1pt} makes sense to approximate $f(t):=\E S(t)-t$ by
$\ff(t):=n(tp)^r/r!-t$. An elementary calculation shows
that $\ff$ has, on $[0,\infty)$ a unique,\vspace*{1pt} global minimum at $\tc$
given by
(\ref{tc}), and that the minimum value is $\ff(\tc)=-\ac$.
We obtain, for example, the following estimate.
%
%
\begin{lemma}\label{L10r}
Suppose that $r\ge2$, $n\qw\ll p\ll n^{-1/r}$ and $a=o(n)$. Then
\[
\sup_{0\le x\le10r} \biggl|S(x\tc)-\frac1r x^r\tc\biggr|
=o_{\mathrm p}(\tc).
\]
\end{lemma}
\begin{pf}
First, (\ref{esapp}) and (\ref{tc1}) yield, recalling (\ref{tccond}),
uniformly for \mbox{$x\le10r$},
\[
\E S(x\tc)=n x^r\frac{\tc^rp^r}{r!}\bigl(1+o(1/x)\bigr)
=x^r\frac{\tc}r\bigl(1+o(1/x)\bigr)
=\frac{x^r}{r}{\tc}+o(\tc).
\]
Further, Lemma~\ref{L2} yields by (\ref{pia}) and (\ref{tc1}),
\[
\sup_{0\le x\le10r} |S(x\tc)-\E S(x\tc)|^2 = \Op(n\pi(10r\tc)) =
\Op(n\tc^rp^r) = \Op(\tc) = o_{\mathrm p}(\tc^2),
\]
and the result follows.
\end{pf}

We shall use Lemma~\ref{L10r} to prove now that in the subcritical case
($a\sim\alpha a_c$ with $\alpha<1$) there exists $t< t_c$ such that
w.h.p. $A(t)\leq t$
and then determine the precise value of $A(T) = T$.\vadjust{\goodbreak}
\begin{pf*}{Proof of Theorem~\ref{T1}\ref{T1sub}}
The assumption on $a$ may be written by~(\ref{ac}),
%
%
\begin{equation}
\label{t1aa}
a=\bigl(\ga+o(1)\bigr)\ac=\bigl(\ga(1-r\qw)+o(1)\bigr)\tc.
\end{equation}
Hence, (\ref{as}) and Lemma~\ref{L10r}, taking $x=1$, yield
\begin{eqnarray*}
A(\tc)-\tc&=& S(\tc)+a-\tc=\tc/r+ o_{\mathrm p}(\tc)+a-\tc\\
&=& \tc
\bigl(r\qw+\ga(1-r\qw)-1+o_{\mathrm p}(1)\bigr).
\end{eqnarray*}
Since $\ga(1-r\qw)<1-r\qw$, w.h.p. $A(\tc)-\tc<0$, and
thus, by (\ref{t2}), $T<\tc$.

We apply Lemma~\ref{L10r} again, now taking $x=T/\tc$, and see that
$S(T)=(T/\tc)^r\tc/r+o_{\mathrm p}(\tc)$.
Since $S(T)=A(T)-a=T-a$, we find, using (\ref{t1aa}), that
\[
T-\ga(1-r\qw)\tc=S(T)+o(\tc)
=\biggl(\frac{T}{\tc}\biggr)^r\frac{\tc}{r}+o_{\mathrm p}(\tc)
\]
and thus
%
%
\begin{equation}\label{ika}
r\frac{T}{\tc}-(r-1)\ga=\biggl(\frac{T}{\tc}\biggr)^r+o_{\mathrm p}(1).
\end{equation}
Since the function $h(x):=rx-x^r$ is strictly increasing from
0 to $r-1$ on $[0,1]$, (\ref{ika})
implies (using the fact just shown that $T/\tc<1$ w.h.p.) that
$T/\tc\pto y$, where $y$ is the unique root in
$\oi$ of $h(y)=(r-1)\ga$, that is, $y=\gf(\ga)$ given by (\ref{ika2}).

This proves the first assertion, and if $\ga>0$, the second follows.
If $\ga=0$, then $a=o(\tc)$, and
(\ref{esapp}) implies, for every fixed $\gl>0$,
$\E S(\gl a)=O(na^rp^r)=o(an\tc^{r-1}p^r)=o(a)$. Hence,
for every fixed $\gl>1$,
$A(\gl a)=S(\gl a)+a=a+o_{\mathrm p}(a)$, so w.h.p.
$A(\gl a)<\gl a$, and thus $a\le T<\gl
a$. Consequently, when $\ga=0$, $T/a\pto1$.
\end{pf*}

We turn to the proof of the supercritical case in Theorem~\ref{T1}. The
following
lemma shows that if the process of activation can escape the bottleneck at
$t_c$, then the process continues until (almost) percolation. The idea
is to
split the time interval $[3t_c,n]$ into different intervals. Then in the
proof of Theorem~\ref{T1}\ref{T1super} and~\ref{T1complete}, it remains to
show that
if $a$ is supercritical, then $A(t)>t$ for $t< 3 t_c$.

Let $\bcx:=\bc\go(n)$, where
$\go(n)\to\infty$ slowly but is otherwise arbitrary.
%
%
\begin{lemma}
\label{Lbulk}
Suppose that $r\ge2$ and $n\qw\ll p\ll n^{-1/r}$.
Then, for any $a$,
w.h.p. $A(t) > t$ for all $t\in[3\tc,n-\bcx]$.
\end{lemma}
\begin{pf}
By (\ref{as}), $A(t)=S_{n-a}(t)+a\ge
S_n(t)$, so it suffices to show that $S_n(t)>t$ (or
equivalently, to take $a=0$).
We separate the proof into a number of different cases for different ranges
of $t$.
We assume at some places, without further mention, that $n$ is large enough.\vadjust{\goodbreak}

\textit{Case} 1: $t\in[3\tc,8r\tc]$.
By Lemma~\ref{L10r}, w.h.p. for all such $t$,
\[
S_n(t)\ge\frac1r\biggl(\frac{t}{\tc}\biggr)^r\tc-\tc
\ge\frac{3^{r-1}}r t-\tc\ge\frac32 t-\tc> t.
\]

\textit{Case} 2: $t\in[8r\tc,p\qw]$.
Let $t_j:=2^jr\tc$, $j\ge1$, and let $J:=\min\{j\ge1\dvtx\break pt_j\ge1\}$.
For $\tc\le t\le p\qw$,
using (\ref{tc}),
\begin{eqnarray*}
\pi(t)&\ge&\pmatrix{t\cr r}p^r(1-p)^{t-r}
= \frac{t^r}{r!}p^re^{-tp}\bigl(1+o(1)\bigr)
\\[-2pt]
&\ge& \frac13\frac{t^rp^r}{r!}
= \frac13t\biggl(\frac{t}{\tc}\biggr)^{r-1}\frac{\tc^{r-1}p^r}{r!}
= \frac{t}{3rn}\biggl(\frac{t}{\tc}\biggr)^{r-1}.
\end{eqnarray*}
Hence, for $3\le j\le J-1$,
$\E S_n(t_j)=n\pi(t_j)\ge\frac{2^j}3 t_j\ge\frac83 t_j$,
and thus, using Chebyshev's inequality and (\ref{svar}),
\[
{\mathbb P}\bigl(S_n(t_j)\le2t_j\bigr)
\le
{\mathbb P}\biggl(S_n(t_j)\le\frac34\E S_n(t_j)\biggr)
\le
\frac{\Var S_n(t_j)}{((1/4)\E S_n(t_j))^2}
\le\frac{16}{n\pi(t_j)}\le\frac{6}{t_j}.
\]

Hence,
\begin{eqnarray*}
{\mathbb P}\bigl(S_n(t)\le t \mbox{ for some $t\in[8r\tc,t_J]$}\bigr)
&\le& \sum_{j=3}^{J-1} {\mathbb P}\bigl(S_n(t_j)\le2t_j\bigr)
\\[-2pt]
&\le& \sum_{j=3}^{J-1} \frac6{t_j}
<\frac{12}{t_3}
<\frac2{r\tc}=o(1).
\end{eqnarray*}

\textit{Case} 3:
$t\in[p\qw,c_1 n]$ for a suitable small $c_{1}>0$.
Let $t_1':=\lceil p\qw\rceil$. Then
\[
\pi(t'_1)={\mathbb P}\bigl(\Bin(t'_1,p)\ge r\bigr)
={\mathbb P}\bigl(\Po(t'_1p)\ge r\bigr)+O(p) \ge2c_{1}
\]
for some small $c_{1}$. Hence w.h.p.
$S_n(t'_1)>c_{1} n$
and consequently
$S_n(t)\ge S_n(t'_1)>c_1n\ge t$.

\textit{Case} 4:
$t\in[c_{1} n,n-p\qw]$.
Let $t_2':=\lfloor c_{1} n\rfloor$ and $t'_3:=n-p\qw$.
Then
\begin{eqnarray*}
1-\pi(t'_2)&=&{\mathbb P}\bigl(\Bin(t'_2,p)<r\bigr)
\\[-1pt]
&=& O((t'_2p)^{r-1}e^{-t'_2p})
=O((np)^{r-1}e^{-c_{1} np})\\[-1pt]
&=&o((np)\qw).
\end{eqnarray*}
Thus,
$\E(n-S_n(t'_2))=n(1-\pi(t'_2))=o(p\qw)$, and
w.h.p.,
$n-S_n(t'_2)<p\qw$, that is, $S_n(t'_2)>n-p\qw=t'_3$.

\textit{Case} 5:
$t\in[n-p\qw,n-\bcx]$.
We have $t'_3:=n-p\qw$.
Then
\begin{eqnarray*}
1-\pi(t'_3) &=& {\mathbb P}\bigl(\Bin(\lfloor t'_3\rfloor,p)<r\bigr)
\\[-1pt]
&=& O((t'_3p)^{r-1}e^{-t'_3p})
=O((np)^{r-1}e^{-np})\\[-1pt]
&=&O(\bc/n).
\end{eqnarray*}
Hence,
$\E(n-S_n(t'_3))=n(1-\pi(t'_3))=O(\bc)=o(\bcx)$, and
thus w.h.p.
$n-S_n(t'_3)<\bcx$, that is, $S_n(t'_3)>n-\bcx$.\vadjust{\goodbreak}
\end{pf}
%
%
\begin{remark}
\label{Rbulk}
The proof shows that once we reach at least $1/p$
active vertices, the active set will w.h.p. grow to at least
$n-\bcx$ in at most 3 generations.
(Hence, the size then is $n-\Op(\bc)$; see \cite{SJN6}, Lemma 3.)
\end{remark}
%
%
\begin{lemma}
\label{Lmin}
\[
\min_{x\ge0}\biggl(\frac{x^r}{r}-x\biggr) =\frac1r-1,
\]
attained at $x=1$ only.
\end{lemma}
\begin{pf}
Elementary calculus.
\end{pf}
\begin{pf*}{Proof of Theorem~\ref{T1}\ref{T1super} and~\ref{T1complete}}
For $0\le t\le3\tc$, we may assume $a\le3\tc$
since otherwise $A(t)>t$ trivially.
In this case, Lemmas~\ref{L10r} and~\ref{Lmin} (with $x=t/\tc$)
show that w.h.p., uniformly in
$t\le3\tc$,
\begin{eqnarray*}
A(t)&=&a+S(t)
\ge (1+\delta)(1-r\qw)\tc+\frac1r\biggl(\frac{t}{\tc}\biggr)^r\tc-o(\tc)
\\
&\ge& \delta(1-r\qw)\tc+\frac{t}{\tc}\tc-o(\tc)
>t.
\end{eqnarray*}
This and Lemma~\ref{Lbulk} show that w.h.p. $A(t)>t$ for all
$t\le n-\bcx$, and thus $\Ax>n-\bcx$.

Hence $n-\Ax<\bcx=\bc\go(n)$ w.h.p., for any choice
of $\go(n)\to\infty$, which is equivalent to
$n-\Ax=\Op(\bc)$; see, for example, \cite{SJN6}, Lemma 3.
This proves~\ref{T1super}.

If $\bc\to0$, we may choose $\bcx=1$; then
w.h.p. $n-\Ax<1$, so $\Ax=n$.
Conversely, if $\bc\not\to0$, then, at least for a
subsequence, there exists with probability at least $c>0$ a
vertex with degree $\le r-1$, and with probability $1-a/n$,
this vertex will never be activated so $\Ax<n$; see Remark~\ref{Rdegrees}.
This proves~\ref{T1complete}.
\end{pf*}
\begin{pf*}{Proof of Theorem~\ref{T2}}
Choose $\bcx:=np\bc\gg\bc$. By (\ref{bc}),
$\bcx p=(np)^{r+1}e^{-np}/(r-1)!\to0$.
Hence, $(n-\bcx)p=np+o(1)\to\infty$ and
\begin{eqnarray*}
1-\pi(n-\bcx)
&=& {\mathbb P}\bigl(\Bin(\lfloor n-\bcx\rfloor,p)\le r-1\bigr)
\\
&\sim& \frac{(n-\bcx)^{r-1}p^{r-1}}{(r-1)!}(1-p)^{n-\bcx}\\
&\sim&1-\pi(n).
\end{eqnarray*}
Consequently [see (\ref{bc2a})],
%
%
\begin{eqnarray}\label{magn}
\E\bigl(A(n)-A(n-\bcx)\bigr)
&=&
\E\bigl(S(n)-S(n-\bcx)\bigr)
\nonumber\\
&\le& n\bigl(\pi(n)-\pi(n-\bcx)\bigr)\\
&=&o\bigl(n\bigl(1-\pi(n)\bigr)\bigr)
=o(\bcq).
\nonumber
\end{eqnarray}

By assumption and Lemma~\ref{Lbulk}, w.h.p. $T>n-\bcx$,
and thus $A(n-\bcx)\le A(T)\le A(n)$.
Hence (\ref{magn}) implies
%
%
\begin{equation}\label{fri}
\Ax=T=A(T)=A(n)+o_{\mathrm p}(\bcq).
\end{equation}
Further,
%
%
\begin{equation}
\label{fri2}
n-A(n)=n-a-S(n)\in\Bin\bigl(n-a,1-\pi(n)\bigr)
\end{equation}
with
mean $(n-a)(1-\pi(n))\sim\bcq$; see
(\ref{bc2a}).

If $\bc\to\infty$, then
$\bcq\sim\bc$; thus
(\ref{fri2}) implies
$n-A(n)=\bc+o_{\mathrm p}(\bc)$, and (\ref{fri}) yields~\ref{T2oo}.

In $\bc\to b<\infty$, then
$\bcq=\bc+o(1)\to b$; thus
(\ref{fri}) yields
$\Ax=A(n)+o_{\mathrm p}(1)$, and hence (since the variables are
integer valued) $\Ax=A(n)$ w.h.p.
Further, in this case (\ref{fri2}) implies
$n-A(n)\dto\Po(b)$,
and~\ref{T20} and~\ref{T2b} follow.
\end{pf*}
\begin{pf*}{Proof of Theorem~\ref{Tpc}}
An easy consequence of Theorem~\ref{T1}.
\end{pf*}

We end this section with a proof of Theorem~\ref{TJ}, where we start
with a number
of external infections (but no initially active vertices).
As said in Section~\ref{SSdyni},
we do this by a minor variation of our method. We include this proof to show
the flexibility of the method, but we omit parts that are identical or
almost identical
to the proofs above.
\begin{pf*}{Proof of Theorem~\ref{TJ}}
In order to preserve independence between vertices, we consider the model
with a Poisson
number $W\in\Po(\mu)$ of external infections (independent of everything
else). Then each vertex $i$ receives
$W_i\in\Po(\mu/n)$ external infections, and these random variables are
independent. The analysis in Section~\ref{Ssetup} becomes slightly
modified: the
number of infections (marks) at time $t$ now is $M^\mu
_i(t):=W_i+M_i(t)$, so
$Y_i$ is
replaced by $Y^\mu_i:=\min\{t\dvtx M^\mu_i(t)\ge r\}$ and $S(t)$ is
replaced by
$\smu(t):=\sum_{i=1}^n \ett\{Y^\mu_i\le t\}$. We now have
$A(t)=\smu
(t)$, so
$\Ax=T=\min\{t\ge0\dvtx\smu(t)=t\}$.

We take $\mu=ynp\ac$ for a fixed $y>0$ and claim that if $y<1$, then w.h.p.
$\Ax<\tc/r$ and thus $J_0>W$, while if $y>1$, then w.h.p. $\Ax
=n-o_{\mathrm p}
(n)$ and
thus $J_0\le W$. The result then follows by taking $y=1\pm\eps/2$ for small
$\eps>0$.

To prove these claims, we first note that
$ \E\smu(t)=n{\mathbb P}(M^\mu_i(t)\ge r)$ with, for such $\mu$ and
$t=O(\tc)$,
%
%
\begin{eqnarray}\label{tjo}\quad
{\mathbb P}\bigl(M^\mu_i(t)\ge r\bigr)
&=& {\mathbb P}\bigl(W_i+M_i(t)\ge r\bigr)
\nonumber\\
&=& \sum_{j=0}^{r-1}{\mathbb P}(W_i=j){\mathbb P}\bigl(M_i(t)\ge
r-j\bigr)+{\mathbb
P}(W_i\ge r)
\\
&\sim& \sum_{j=0}^{r}\frac{(\mu/n)^j}{j!}\cdot\frac{(tp)^{r-j}}{(r-j)!}
=\frac{(tp+\mu/n)^r}{r!}
=\frac{p^r(t+y\ac)^r}{r!}.
\nonumber
\end{eqnarray}
We obtain as in Lemma~\ref{L10r}, using versions of Lemmas
\ref{Lmartin} and~\ref{L2} for~$\smu(t)$,
%
%
\begin{equation}\label{smu}
\sup_{0\le x\le10r} \biggl|\smu(x\tc)-\frac1r (x+y\ac/\tc)^r\tc\biggr|
=o_{\mathrm p}(\tc).
\end{equation}
Recall that $\ac/\tc=1-1/r$ by (\ref{ac}). If $y<1$, then (\ref
{smu}) with
$x=1/r$
implies that w.h.p.
\[
A\biggl(\frac{\tc}{r}\biggr)=\smu\biggl(\frac{\tc}{r}\biggr)
<\frac1r\biggl(\frac1r+\biggl(1-\frac1r\biggr)\biggr)^r\tc=\frac\tc r
\]
and thus $\Ax=T<\tc/r$ as claimed.

Conversely, if $y>1$, then Lemma~\ref{Lmin} shows that
\[
\frac{(x+y\ac/\tc)^r} r \ge x+y\frac{\ac}{\tc}+\biggl(\frac1r-1\biggr)
= x+(y-1)\biggl(1-\frac1r\biggr).
\]
Hence, (\ref{smu}) shows that w.h.p.
$A(x\tc)=\smu(x\tc)>x\tc$ for $x\le10r$, and
thus $\Ax=T>10r\tc$. Further, since $\smu(t)\ge S_n(t)$, Lemma~\ref{Lbulk}
implies that w.h.p. $A(t)>t$ for all $t\in[3\tc,n-\bcx]$, and thus
$\Ax\ge n-\bcx$ w.h.p., which proves the second claim and completes
the proof.
\end{pf*}

Note that (\ref{tjo}) and (\ref{smu}) show that $A(t)=\smu(t)$ is,
to the
first order,
$\E S_n(t)$ shifted horizontally by $\mu/(np)=y\ac$, while in our standard
model
$A(t)$ is $S_t$ shifted vertically by $a$. Since we study the hitting time
of the linear barrier $A(t)=t$, these are essentially equivalent.

\section{\texorpdfstring{Proofs of Theorems \protect\ref{Tac}--\protect\ref{TGsub} and \protect\ref{TDG2}}
{Proofs of Theorems 3.6--3.8 and 4.5}}\label{sec9}
\label{SGauss}

We begin with an estimate of~$\tpi(t)$ defined in (\ref{tpi}).
%
%
\begin{lemma}
\label{LS1}
Suppose that $r\ge2$ and $n\qw\ll p\ll n^{-1/r}$.
Then, for large~$n$,
$n\tpi(t)\ge1.4 t$ for $t\in[3\tc,n/2]$.
\end{lemma}
\begin{pf}
Assume not. Then we can find, for a subsequence $n=n_k\to\infty$,
$t=t_k\in[3\tc,n/2]$ such that $n\tpi(t)<1.4 t$.
Selecting a subsequence, we may further assume that $pt\to
z\in[0,\infty]$.
We consider three cases separately.

{
\renewcommand\thelonglist{(\roman{longlist})}
\renewcommand\labellonglist{\thelonglist}
\begin{longlist}
\item
$z=0$, that is, $pt\to0$. Then, from (\ref{tpi}) and (\ref{tc1}),
\[
n\tpi(t)\sim n\frac{(pt)^r}{r!}
=\frac{t^r}{r\tc^{r-1}}
\ge\frac{3^{r-1}}rt\ge\frac32t.
\]

\item$pt\to z\in(0,\infty)$.
Then $n\tpi(t)\sim n\psi(z)$ with $\psi(z)>0$, and
$t=O(1/p)=o(n)\ll n\tpi(t)$.

\item
$z=\infty$, that is, $pt\to\infty$. Then $n\tpi(t)\sim n\ge2t$.
\end{longlist}
}

In all cases we have for large $n$ a contradiction to $n\tpi(t)<1.4 t$.
\end{pf}
%
%
\begin{lemma}
\label{Lpcx}
Suppose $r\ge2$ and $a\to\infty$ with $a=o(n)$. Then
$\pc$ and $\pcx$
defined by (\ref{pc}) and (\ref{pcx}) satisfy
$\pcx\sim\pc$.
In particular, $n\qw\ll\pcx\ll n^{-1/r}$.\vadjust{\goodbreak}
\end{lemma}
\begin{pf}
Let $p=y\pc$ for some fixed $y>0$, and define $\tc$ and $\ac$ by
(\ref{tc}) and~(\ref{ac}). Then $\tc=y^{-(r-1)/r}(r/(r-1))a$ and
$\ac=a y^{-(r-1)/r}$. Further,
$n\qw\ll p \ll n\qwr$ and, by (\ref{tccond}), $p\tc\to0$ and
$\tc=o(n)$. Hence, if $x=O(1)$, and $t=x\tc$, then $pt=o(1)$,
$t=o(n)$ and, uniformly in bounded $x$,
by (\ref{tpi}) and (\ref{tc1}),
\[
\tpi(t)=\frac{(pt)^r}{r!} + O((pt)^{r+1})
=\frac{x^r\tc}{nr}\bigl(1+o(1)\bigr).
\]
Hence, uniformly in $x\le3$,
%
%
\begin{equation}\label{jul}
(n-a)\tpi(t)-t = \biggl(\frac{x^r}{r}+o(1)-x\biggr)\tc.
\end{equation}
By Lemma~\ref{LS1}, for large~$n$, $(n-a)\tpi(t)-t\ge0$ for
$t\in[3\tc,n/2]$, and thus, by (\ref{jul}) and Lemma~\ref{Lmin},
%
%
\begin{eqnarray}\label{jull}
\hh(p)
&=&
\inf_{t\le3\tc}\{(n-a)\tpi(t)-t\}
=\biggl(\inf_{x\le3}\biggl(\frac{x^r}r-x\biggr)+o(1)\biggr)\tc
\nonumber\\
&=& \biggl(\frac{1}{r}-1+o(1)\biggr)\tc
=-\bigl(1+o(1)\bigr)\ac
\\
&=& -\bigl(y^{-(r-1)/r}+o(1)\bigr)a.
\nonumber
\end{eqnarray}
Hence, if $y=1-\delta<1$, then $y^{-(r-1)/r}>1$ and thus, for large~$n$,
$\hh(p)<-a$
so $\pcx>p=(1-\delta)\pc$.
Conversely, if $y=1+\delta>1$, then (\ref{jull}) yields, for large~$n$,
$\hh(p)>-a$
so $\pcx<p=(1+\delta)\pc$.

Consequently, $\pcx/\pc\to1$.
\end{pf}

We also need more precise estimates of $S(t)$.
The following Gaussian process limit is fundamental.
$D[0,B]$ denotes the space of right-continuous functions on $[0,B]$,
with the Skorohod topology; see, for example, \cite{Bill} (for $B=1$; the
general case is similar by a change of variables)
or \cite{Kallenberg}, Chapter 16.
%
%
\begin{lemma}
\label{LG}
Suppose $r\ge2$ and $a\to\infty$ with $a=o(n)$.
Then
%
%
\begin{equation}
\label{lg}
Z(x):=
\frac{S(x\tc)-\E S(x\tc)}{\sqrt{\tc}}
\dto
W(x^r/r)
\end{equation}
in $D[0,B]$ for any fixed~$B$, where $W$ is a standard Brownian motion.
\end{lemma}

The conclusion, convergence in $D[0,B]$ for every fixed~$B$, can also
be expressed as convergence in $D\ooo$.
\begin{pf*}{Proof of Lemma~\ref{LG}}
This is a result on convergence of empirical distribution functions
(of $\{Y_i\}$); cf. \cite{Bill}, Theorem 16.4; we get here a
Brownian motion
instead of a Brownian bridge as in \cite{Bill} because we consider
for each~$B$ only a small initial part of the distribution of $Y_i$.\vadjust{\goodbreak}

For every fixed $x>0$, by (\ref{pia}) and (\ref{tccond}),
$\pi(x\tc)\sim(x\tc p)^r/r!\to0$,
and thus by (\ref{svar}) and (\ref{tc1})
\[
\Var S(x\tc)
\sim n\pi(x\tc)
\sim\frac{np^rx^r\tc^r}{r!}
= \frac{x^r\tc}{r}
\to\infty.
\]
Hence (\ref{sbin}) and the central limit theorem yield
$Z(x)\dto N(0,x^r/r)$ for every $x>0$, which proves $Z(x)\dto W(x^r/r)$ for
each fixed $x$.

This is easily extended to finite-dimensional convergence:
Suppose that $0<x_1<\cdots<x_\ell$ are fixed, and let
$I_{ij}:=\ett\{Y_i\in(x_{j-1}\tc,x_j\tc]\}$, with $x_0=0$.
Thus, $S(x_j\tc)-S(x_{j-1}\tc)=\sum_{i=1}^{n-a}I_{ij}$.
Then, for $1\le j\le\ell$ and $k\neq j$,
\begin{eqnarray*}
\E I_{ij}&=&\pi(x_j\tc)-\pi(x_{j-1}\tc), \\
\Var I_{ij}&=&\E I_{ij}(1-\E I_{ij})\sim\E I_{ij}\\
&=&\pi(x_j\tc)-\pi
(x_{j-1}\tc)
\sim\biggl(\frac{x_j^r}r-\frac{x_{j-1}^r}r\biggr)\frac{\tc}n,\\
\Cov(I_{ij},I_{ik})&=&-\E I_{ij}\E I_{ik}
=O(\pi(x_\ell\tc)^2)
=O\bigl((\tc/n)^2\bigr)
=o(\tc/n).
\end{eqnarray*}
Note that $(I_{ij})_{j=1}^\ell$, $i=1,2,\ldots,n$, are i.i.d.
random vectors. The multi-dimen\-sional central limit theorem with (e.g.)
the Lindeberg condition (which follows from the one-dimensional
version in, for example, \cite{Gut}, Theorem~7.2.4, or~\cite{Kallenberg}, Theorem~5.12,
by the Cram\'er--Wold device) thus shows that
$(Z(x_j)-Z(x_{j-1}))_{j=1}^\ell\dto(V_j)_{j=1}^\ell$
with $V_j$ jointly normal with $\E V_j=0$, $\Var
V_j=x_j^r/r-x_{j-1}^r/r$ and $\Cov(V_j,V_k)=0$ for $j\neq k$. Hence,
$(V_j)_{j=1}^\ell\eqd(W(x_j^r/r)-W(x_{j-1}^r/r)_{j=1})^\ell$, and thus
$(Z(x_j))_{j=1}^\ell\dto(W(x_j^r/r))_{j=1}^\ell$.

To show (\ref{lg}), it thus remains to show tightness of $Z(x)$.
We use \cite{Bill}, Theorem 15.6, with $\gamma=2$ and $\ga=1$
(an alternative would be to instead use Aldous's tightness criterion
(\cite{Kallenberg}, Theorem 16.11));
it thus suffices to prove
that, for every $x_1,x_2,x_3$ with $0\le x_1\le x_2\le x_3\le B$,
and some constant $C$ depending on $B$ but not on $n$ or $x_1,x_2,x_3$,
%
%
\begin{equation}\label{tight}
\E\{|Z(x_2)-Z(x_1)|^2|Z(x_3)-Z(x_2)|^2\} \le C (x_3-x_1)^2.
\end{equation}
With the notation above and $I'_{ij}:=I_{ij}-\E I\ijx$, the left-hand
side of (\ref{tight}) can be written
\[
\tc\qww\E\sum_{i,j,k,l=1}^{n-a}I'_{i2}I'_{j2}I'_{k3}I'_{l3}
=
\tc\qww\sum_{i,j,k,l=1}^{n-a}\E(I'_{i2}I'_{j2}I'_{k3}I'_{l3}).
\]

By independence, the only nonzero terms are those where $i,j,k,l$
either coincide in two pairs, or all four indices coincide, and it follows
easily that (for any~$i$)
%
%
\begin{equation}\label{emma}
\E\{|Z(x_2)-Z(x_1)|^2|Z(x_3)-Z(x_2)|^2\}
\le
3\tc\qww(n-a)^2\E I_{i2}\E I_{i3}.
\end{equation}

Further, since each $Y_i$ is integer-valued,
the left-hand side of (\ref{tight}) vanishes unless there is at least
one integer in each of the intervals $(x_1\tc,x_2\tc]$ and
$(x_2\tc,x_3\tc]$, which implies that $x_3\tc-x_1\tc>1$, so we only
have to consider this case. It follows from (\ref{yik}) that
for $m\le x_3\tc\le B\tc$,
\[
{\mathbb P}(Y_i=m)\le\frac{m^{r-1}}{(r-1)!}p^r
\le\frac{B^{r-1}\tc^{r-1}p^r}{(r-1)!}
=\frac{B^{r-1}}{n}
\]
and thus, assuming $x_3\tc-x_1\tc>1$,
\begin{eqnarray*}
\E I_{i2}+\E I_{i3}&
\le& (\lfloor x_3\tc\rfloor-\lfloor x_1\tc\rfloor)\frac{B^{r-1}}n
\le({x_3\tc}-{x_1\tc}+1)\frac{B^{r-1}}n
\\
&\le& 2({x_3\tc}-{x_1\tc})\frac{B^{r-1}}n.
\end{eqnarray*}
Consequently, (\ref{emma}) yields,
\begin{eqnarray*}
\E\{|Z(x_2)-Z(x_1)|^2|Z(x_3)-Z(x_2)|^2\}
&\le&
3\frac{n^2}{\tc^2} \biggl(2({x_3\tc}-{x_1\tc})\frac{B^{r-1}}n\biggr)^2
\\ &=& 12(x_3-x_1)^2B^{2(r-1)},
\end{eqnarray*}
which proves (\ref{tight}) with $C=12B^{2(r-1)}$.
The proof is complete.
\end{pf*}

We also need a more careful estimate of $\pi(t)$ than above, and we
use the corresponding Poisson probability $\tpi(t)$ defined in
(\ref{tpi}).
%
%
\begin{lemma}
\label{Ltpi}
Assume $n\qw\ll p \ll n\qwr$.
Uniformly for $t\ge1$, $\pi(t)=\tpi(t)(1+O(t\qw))$.
In particular, uniformly for $t\le3\tc$,
\[
\pi(t)=\tpi(t)+O\bigl((pt)^r/t\bigr)
=\tpi(t)+O\bigl(\tpi(\tc)/\tc\bigr)
=\tpi(t)+O(n\qw).
\]
\end{lemma}
\begin{pf}
Assume first $pt\le1$.
By (\ref{pi}),
\begin{eqnarray*}
\pi(t)
&=&\sum_{j=r}^t{\mathbb P}\bigl(\Bin(t,p)=j\bigr)
=\sum_{j=r}^t\frac{t^j}{j!}\biggl(1+O\biggl(\frac{j^2}t\biggr)\biggr)p^j(1-p)^{t+O(j)}
\\
&=&\sum_{j=r}^\infty\frac{(pt)^j}{j!}e^{-pt+O(tp^2)}\bigl(1+O(j^2/t+jp)\bigr)
\\
&=&\tpi(t)\bigl(1+O(tp^2+t\qw+p)\bigr)
=\tpi(t)\bigl(1+O(t\qw)\bigr).
\end{eqnarray*}
For $pt>1$, $\tpi(t)$ is bounded below, and the result follows from
(\ref{dtv}).

If $t\le3\tc$, then $t=O(\tc)=o(1/p)$ by (\ref{tccond}),
and thus, using (\ref{tpi}) and~(\ref{tc1}),
\[
\pi(t)-\tpi(t)
=O\bigl(\tpi(t)/t\bigr)
=O\bigl((pt)^r/t\bigr)
=O\bigl((p\tc)^r/\tc\bigr)
=O(1/n).
\]
\upqed
\end{pf}
\begin{pf*}{Proof of Theorem~\ref{Tac}}
It suffices to consider $a$ such that $a\sim\ac=(1-r\qw)\tc$.
It then follows by (\ref{as}) and Lemma~\ref{L10r} that, uniformly
for\vadjust{\goodbreak}
$x\le10r$,
%
%
\begin{eqnarray}\label{ax1}
A(x\tc)-x\tc& = &a+S(x\tc)-x\tc
=\ac+\frac1r x^r\tc-x\tc+o_{\mathrm p}(\tc)
\nonumber\\[-8pt]\\[-8pt]
&=&\biggl(1-r\qw+\frac1r x^r-x\biggr)\tc+o_{\mathrm p}(\tc).
\nonumber
\end{eqnarray}
By Lemma~\ref{Lmin},
the coefficient $1-r\qw+x^r/r-x$
equals $0$ at $x=1$ but is strictly positive for all other
$x\ge0$. It
follows that for every $\delta>0$,
w.h.p. $A(x\tc)-x\tc>0$ for all $x\in[0,1-\delta]\cup[1+\delta,10r]$.
By a simple standard argument, there thus exists a sequence
$\delta_n\to0$,
where we may further assume that $\delta_n>|\tcx/\tc-1|$,
such that
w.h.p. $A(x\tc)-x\tc>0$ for all $x\in[0,1-\delta_n]\cup[1+\delta_n,10r]$.

Hence, w.h.p. either $T\in[(1-\delta_n)\tc,(1+\delta_n)\tc]$, or $A(t)>t$
for all $t\le10r\tc$; in the latter case,
for any $\bcx\gg\bc$,
w.h.p. $A(t)>t$ for all $t\le n-\bcx$ by Lem\-ma~\ref{Lbulk},
so $T\ge n-\bcx$; hence $T=n-\Op(\bc)$ and, more precisely,
provided $a=o(n)$,
Theorem~\ref{T2} applies.

We thus only have to investigate the interval $[(1-\delta_n)\tc
,(1+\delta
_n)\tc]$
more closely.
By
the Skorohod coupling theorem (\cite{Kallenberg}, Theorem 4.30),
we may assume that the processes for different $n$ are coupled such
that the limit~(\ref{lg}) holds a.s., and not just in distribution.
Since convergence in $D[0,B]$ to a~continuous function is equivalent
to uniform convergence, this means that (a.s.) $Z(x)\to W(x^r/r)$
uniformly for $x\le B$; in particular,
uniformly for $x\in[1-\delta_n,1+\delta_n]$,
%
%
\begin{eqnarray}\label{ems}
S(x\tc)& = &(n-a)\pi(x\tc)+\tc^{1/2} Z(x)
\nonumber\\[-8pt]\\[-8pt]
& = &(n-a)\pi(x\tc)+\tc^{1/2} \bigl(W(1/r)+o(1)\bigr).
\nonumber
\end{eqnarray}
Let $\xi:=W(1/r)\in N(0,1/r)$. Then, by (\ref{ems}) and Lemma~\ref{Ltpi},
uniformly for $x\in[1-\delta_n,1+\delta_n]$,
%
%
\begin{eqnarray}
S(x\tc) & = &(n-a)\tpi(x\tc)+O(1)+\tc^{1/2} \bigl(\xi+o(1)\bigr)
\nonumber\\[-8pt]\\[-8pt]
& = &(n-a)\tpi(x\tc)+\tc^{1/2} \bigl(\xi+o(1)\bigr)
\nonumber
\end{eqnarray}
and thus, refining (\ref{ax1}),
%
%
\begin{eqnarray}\label{emm0}
A(x\tc)-x\tc&=& a+S(x\tc)-x\tc\nonumber\\[-8pt]\\[-8pt]
&=&a+(n-a)\tpi(x\tc)-x\tc+\tc^{1/2}\xi+o_{\mathrm
p}(\tc^{1/2}).\nonumber
\end{eqnarray}
Hence, recalling (\ref{acx}) and that the minimum there is attained at
$\tcx\in[(1-\delta_n)\tc,(1+\delta_n)\tc]$,
%
%
\begin{eqnarray}\label{emm00}
&&\min_{t\in[(1-\delta_n)\tc,(1+\delta_n)\tc]} \frac{ A(t) -
\lfloor t\rfloor} {1 -\tpi(t)}
\nonumber\\
&&\qquad= a + \min_{t\in[(1-\delta_n)\tc,(1+\delta_n)\tc]} \frac{n\tpi
(t) - t} {1-\tpi
(t)} + \tc^{1/2}\xi+ o_{\mathrm p}(\tc^{1/2})
\\
&&\qquad=a-\acx+\tc^{1/2}\xi+o_{\mathrm p}(\tc^{1/2}).
\nonumber
\end{eqnarray}
%
We have shown that w.h.p.
$\Ax=T\le(1+\delta_n)\tc$
if and only if this minimum is $\le0$, and otherwise $T=n-\Op(\bc)$,
and the results follow; for~\ref{Tac-} we also observe that (\ref
{emm0}) and
(\ref{emm00}) imply that w.h.p. $A(\tcx)-\tcx<0$ and thus $T<\tcx$.
For example, in~\ref{Tac0} we have
\begin{eqnarray*}
a-\acx+\tc^{1/2}\xi+o_{\mathrm p}(\tc^{1/2}) & = &y\ac^{1/2}+\tc
^{1/2}\xi+o_{\mathrm p}(\tc^{1/2})
\\
& = &\bigl((r-1)^{1/2} y+r^{1/2}\xi+o_{\mathrm p}(1)\bigr)(\tc/r)^{1/2},
\end{eqnarray*}
and the probability that this is positive tends to
\[
{\mathbb P}\bigl((r-1)^{1/2} y+r^{1/2}\xi>0\bigr)=\Phi\bigl((r-1)^{1/2} y\bigr),
\]
since $r^{1/2}\xi\in N(0,1)$.
\end{pf*}
\begin{pf*}{Proof of Theorem~\ref{Tpcxx}}
It suffices to consider $p\sim\pcx\sim\pc$, which implies that
$\ac=\ac(p)\sim\ac(\pc)=a$. Hence the arguments in the proof of
Theorem~\ref{Tac} apply. In particular, again it suffices to consider
$t\in J=J_n:=[(1-\delta_n)\tc,(1+\delta_n)\tc]$, where now
$\tc=\tc(\pc)=(r/(r-1))a$.
The infimum in~(\ref{hp}) is attained for some $t=\tcxx$, where by
Lemma~\ref{LS1} $\tcxx\le3\tc$ for large~$n$, and an argument as in
(\ref{ax1}) shows that $\tcxx\sim\tc$. We may assume that
$\delta_n$ is chosen such that $\tcxx\in J$.
Then, by (\ref{emm0}),
%
%
\begin{equation}\label{emm1}\quad
\min_{t\in J}\{A(t)-\lfloor t\rfloor\}
=a+\min_{t\in J}\{(n-a)\tpi(t)-t\}+\tc^{1/2}\xi+o_{\mathrm
p}(\tc^{1/2}),
\end{equation}
where, by (\ref{hp}) and the comments just made (for large $n$),
%
%
\begin{equation}\label{em2}
a+\min_{t\in J}\{(n-a)\tpi(t)-t\}
=a+\hh(p)
=-\hh(\pcx)+\hh(p).
\end{equation}
Further, writing (\ref{hp}) as $\hh(p):=\min_t\{F(tp)-t\}$, with
$F(x)=(n-a)\psi(x)$, we have
at the minimum point $t:=\tcxx$ the derivative
$pF'(p\tcxx)-1=0$.
Hence, uniformly for
$|p_1-\pcx|\le\eps\pcx$ and $|t-\tc|\le\eps\tc$, for any
$\eps=\eps_n\to0$,
using (\ref{tpii}),
$F'(p_1t)=(1+o(1))F'(p\tcxx)=(1+o(1))/\pcx$ and thus by the mean-value
theorem, for some
$p_1$ between $p$ and $\pcx$,
\[
F(tp)-F(t\pcx)= t(p-\pcx)F'(tp_1)= \tc\frac{p-\pcx}{\pcx}\bigl(1+o(1)\bigr).
\]
Since the minimum in (\ref{hp}) may be taken over such $t$ only, for
suitable $\eps_n$, this yields
\[
\hh(p)-\hh(\pcx)= \tc\frac{p-\pcx}{\pcx}\bigl(1+o(1)\bigr).
\]

Consequently, (\ref{emm1}) and (\ref{em2}) yield
\[
\min_{t\in J}\{A(t)-\lfloor t\rfloor\}
=
\tc\frac{p-\pcx}{\pcx}\bigl(1+o(1)\bigr)+\tc^{1/2}\xi+o_{\mathrm
p}(\tc^{1/2}).
\]
Hence,
\begin{eqnarray*}
{\mathbb P}\Bigl(\min_{t\in J}\{A(t)-\lfloor t\rfloor\}>0\Bigr)
& = &{\mathbb P}\biggl(\tc^{1/2}\frac{p-\pcx}{\pcx}+\xi>0\biggr)+o(1)
\\
& = &{\mathbb P}\biggl(-r^{1/2}\xi<(r\tc)^{1/2}\frac{p-\pcx}{\pcx}\biggr)+o(1),
\end{eqnarray*}
where $r^{1/2}\xi\in N(0,1)$ and $\tc=\frac{r}{r-1}a$, and the
different parts of the theorem
follow.
\end{pf*}
%
%
\begin{lemma}
\label{LD}
Suppose that $r\ge2$ and $n\qw\ll p\ll n^{-1/r}$.
Then, for large~$n$ at least, the minimum point $\tcx$ in (\ref{acx}) is
unique,
and $\tcx\sim\tc$, $\acx\sim\ac$; more precisely,
%
%
\begin{eqnarray}\label{ldt}
\tcx&=&\biggl(1+\frac{p\tc}{r-1}+o(p\tc)\biggr)\tc,
\\
\label{lda}
\acx&=&\biggl(1-\frac1r+\frac{p\tc}{r+1}+o(p\tc)\biggr)\tc
=\biggl(1+\frac{rp\tc}{r^2-1}+o(p\tc)\biggr)\ac.
\end{eqnarray}
\end{lemma}
\begin{pf}
Let
%
%
\begin{equation}\label{ldg}
\gx(t):=\frac{n\tpi(t)-t}{1-\tpi(t)}
=\frac{n-t}{1-\tpi(t)}-n.
\end{equation}
Then
\[
\gx'(t)
=\frac{-(1-\tpi(t))+(n-t)\tpi'(t)}{(1-\tpi(t))^2}
=\frac{(n-t)\tpi'(t)+\tpi(t)-1}{(1-\tpi(t))^2}.
\]
Let
%
%
\begin{equation}\label{ldh}
h(t):=\bigl(1-\tpi(t)\bigr)^2 \gx'(t) = (n-t)\tpi'(t)+\tpi(t)-1.
\end{equation}
Then
$h'(t)= (n-t)\tpi''(t)>0$ for $t<(r-1)/p$, and in particular (for large~$n$)
for $t\le3\tc$; see (\ref{tccond}). Further, $h(0)=-1$ and by (\ref{tpii})
and (\ref{tc}), for large~$n$,
\[
h(3\tc)= 3^{r-1}-1+o(1)>0;
\]
hence, there is a unique $\tcx\in[0,3\tc]$ such that $h(\tcx)=0$, or
equivalently
$\gx'(\tcx)=0$. Further, $\gx''(\tcx)=h'(\tcx)/(1-\tpi(\tcx
))^2>0$ so
$\tcx$ is
the unique minimum point of
$\gx(t)$ in $[0,3\tc]$, as we defined $\tcx$ after (\ref{acx}).

Let $x=\tcx/\tc\in[0,3]$. Then, by (\ref{tpii}), (\ref{tc}) and
(\ref{tpi2}),
\[
0=h(\tcx)=\biggl(1-\frac{\tcx}n\biggr)x^{r-1}e^{-p\tcx} + O\biggl(\frac{\tcx}{n}\biggr)-1
=x^{r-1}e^{-p\tcx}-1+ O\biggl(\frac{\tc}{n}\biggr).
\]
Hence, recalling $n\qw\ll p$ and $p\tc\to0$,
%
%
\begin{equation}\label{ldq}
x=e^{p\tcx/(r-1)}\biggl(1+O\biggl(\frac{\tc}{n}\biggr)\biggr)
=1+\frac{p\tcx}{r-1}+o(p\tc).
\end{equation}
In particular, $x=1+o(1)$, so $\tcx\sim\tc$, and (\ref{ldt})
follows from
(\ref{ldq}).
Finally, substituting (\ref{ldt}) in (\ref{acx}) yields, using (\ref{ldg})
together with $\tpi(\tcx)=O(\tc/n)=o(p\tc)$ by (\ref{tpi2}) and
$n\qw
\ll p$,
and also
(\ref{tpi}) and (\ref{tpi2}),
\begin{eqnarray*}
\frac{\acx}{\tc}
&=&-\frac{\gx(\tcx)}{\tc}
=\bigl(1+o(p\tc)\bigr){\frac{\tcx-n\tpi(\tcx)}{\tc}}
=\bigl(1+o(p\tc)\bigr)\biggl(x-\frac{n\tpi(\tcx)}{\tc}\biggr)
\\[-2pt]
&=&
\bigl(1+o(p\tc)\bigr)\biggl(x-\frac{n(p\tcx)^r}{r! \tc}e^{-p\tcx} \biggl(1+\frac{p\tcx
}{r+1}+o(p\tcx)\biggr)\biggr)
\\[-2pt]
&=& \bigl(1+o(p\tc)\bigr)\biggl(x-\frac{x^r}{r}e^{-p\tc} \biggl(1+\frac{p\tc
}{r+1}+o(p\tc)\biggr)\biggr)
\\[-2pt]
&=& x-\frac{x^r}{r}
\biggl(1-p\tc+\frac{p\tc}{r+1}\biggr)
+o(p\tc)
\\[-2pt]
&=& x-\frac{x^r}{r}
+\frac1r \biggl(p\tc-\frac{p\tc}{r+1}\biggr)
+o(p\tc),
\end{eqnarray*}
and (\ref{lda}) follows by (\ref{ldq}) and $x-x^r/r=1-1/r+O(x-1)^2$.
\end{pf}
\begin{pf*}{Proof of Theorem~\ref{TGsub}}
In case~\ref{TGsub1}, that is, when $\alpha<1$,
by Theorem~\ref{T1}\ref{T1sub},
%
%
\begin{equation}\label{5a}
T=\Ax=\bigl(\gf(\ga)+o_{\mathrm p}(1)\bigr)\tc.
\end{equation}
By Theorem~\ref{Tac}\ref{Tac-}, this holds as well in case~\ref{TGsub2},
that is, when $\alpha=1$ and,
correspondingly, $\gf(\ga)=\gf(1)=1$.
Thus, for any $\alpha\leq1$
there exist $\delta_n\to0$ such that w.h.p.
$T\in I_n:=[(\gf(\ga)-\delta_n)\tc,(\gf(\ga)+\delta_n)\tc]$.
As in the proof of Theorem~\ref{Tac}, we may, by
the Skorohod coupling theorem (\cite{Kallenberg}, Theorem~4.30)
assume that the limit in (\ref{lg}) holds a.s., uniformly in $x\le B$.
For $t\in I_n$, $t/\tc\to\gf(\ga)$,
and (\ref{lg}) then implies that, uniformly for $t\in I_n$,
%
%
\begin{eqnarray}\label{tib}
S(t)&=&\E S(t)+\tc^{1/2} Z(t/\tc)\nonumber\\[-8pt]\\[-8pt]
&=& (n-a)\pi(t)+\tc^{1/2} W\bigl(\gf(\ga)^r/r\bigr)+o_{\mathrm
p}(\tc^{1/2}).\nonumber
\end{eqnarray}
Let $\xi:=W(\gf(\ga)^r/r)\in N(0,\gf(\ga)^r/r)$.
Then, by (\ref{tib}) and Lemma~\ref{Ltpi}, for \mbox{$t\in I_n$},
\[
S(t)
= (n-a)\tpi(t)+\tc^{1/2} \bigl(\xi+o_{\mathrm p}(1)\bigr).
\]
Since w.h.p. $T\in I_n$, we may here substitute $t=T$, and obtain
%
%
\begin{eqnarray}\label{J16}
0&=&A(T)-T=a+S(T)-T\nonumber\\[-8pt]\\[-8pt]
&=& a+(n-a)\tpi(T)-T+\tc^{1/2} \bigl(\xi+o_{\mathrm
p}(1)\bigr).\nonumber\vadjust{\goodbreak}
\end{eqnarray}
Define the function $\gy(t)$ by
%
%
\begin{equation}\label{gt}
\gy(t):= a+(n-a)\tpi(t)-t;
\end{equation}
thus (\ref{tx}) is $\gy(t_*)=0$.
Then we have shown in (\ref{J16}),
%
%
\begin{equation}\label{aahm}
\gy(T)=-\tc^{1/2} \bigl(\xi+o_{\mathrm p}(1)\bigr).
\end{equation}

The function $\gy$ is continuous on $\ooo$ with $\gy(0)=a>0$. Consider
the two
cases separately.

\ref{TGsub1}:
When $\ga<1$ we have, by (\ref{gt}) and (\ref{tpi2}),
$\gy(\tc)=a+(1+o(1))\tc/r-\tc
=a-(1+o(1))\ac
<0
$ (for large $n$), since
$a\sim\ga\ac$.
Further, on $[0,\tc]$, using~(\ref{tpii}) and (\ref{tc}),
%
%
\begin{eqnarray}\label{g}
\gy'(t)&=&(n-a)\tpi'(t)-1
=\frac{n-a}n
\biggl(\frac{t}{\tc}\biggr)^{r-1}e^{-pt}-1\nonumber\\[-8pt]\\[-8pt]
&=& \biggl(\frac{t}{\tc}\biggr)^{r-1}-1+o(1);\nonumber
\end{eqnarray}
this is negative for $t<(1-\eps)\tc$ for any $\eps>0$ and large $n$,
and it follows that (for large $n$, at least), $\gy$ has a unique root
$t_*$ in $[0,\tc]$. It follows from (\ref{tpi2}) and (\ref{ika2}) that
$t_*/\tc\to\gf(\ga)$.

Since also $T/\tc\pto\gf(\ga)$, (\ref{g}) implies that
$\gy'(t)=-(1-\gf(\ga)^{r-1})+o_p(1)$ for all $t$ between $t_*$ and $T$,
and thus the mean value theorem yields
\[
\gy(T)=\gy(T)-\gy(t_*)=(T-t_*)\bigl(-\bigl(1-\gf(\ga)^{r-1}\bigr)+o_p(1)\bigr),
\]
which together with (\ref{aahm}) yields, recalling $\gf(\ga)<1$,
\[
T-t_*=-\bigl(\bigl(1-\gf(\ga)^{r-1}\bigr)\qw+o_p(1)\bigr)\gy(T)
=\bigl(\bigl(1-\gf(\ga)^{r-1}\bigr)\qw\xi+o_p(1)\bigr)\tc^{1/2}.
\]
The result in~\ref{TGsub1} follows.

\ref{TGsub2}:
Let $\gx(t):=\gy(t)/(1-\tpi(t))-a$ and $h(t)$ be as in the proof of
Lemma~\ref{LD}, (\ref{ldg}) and (\ref{ldh}).
We know that $\gx(\tcx)=-\acx$ and $\gx'(\tcx)=0$.
Further, for $t\sim\tc$, we have by
(\ref{tccond}),
(\ref{tpi}), (\ref{tpi2}) and
(\ref{tpii}),
\begin{eqnarray*}
\tpi(t)&\sim&\frac{(tp)^r}{r!}
\sim\frac{\tc}{rn}=o(1),
\\
\tpi'(t)&\sim&\frac{r}{t}\tpi(t)\sim\frac1n,
\\
\tpi''(t)&\sim&\frac{r-1}{t}\tpi'(t)\sim\frac{r-1}{n\tc}.
\end{eqnarray*}
Hence, by (\ref{ldh}),
$
h(t)=o(1)$,
$h'(t)=(n-t)\pi''(t)\sim{(r-1)/\tc}$ and
%
%
\begin{eqnarray}\label{ldgg}
\gx''(t)&=&\frac{h'(t)}{(1-\tpi(t))^2}+2\frac{h(t)\tpi
'(t)}{(1-\tpi(t))^3}
\nonumber\\[-8pt]\\[-8pt]
&=&\frac{r-1}{\tc}\bigl(1+o(1)\bigr)+o\biggl(\frac{1}{n}\biggr)
=\frac{r-1}{\tc}\bigl(1+o(1)\bigr).
\nonumber
\end{eqnarray}
Consequently, a Taylor expansion yields, for $t\sim\tc\sim\tcx$,
%
%
\begin{equation}\label{gtay}
\gx(t)=-\acx+\frac{r-1}{2\tc}(t-\tcx)^2\bigl(1+o(1)\bigr).
\end{equation}

We have $\gy(t_*)=0$ and thus $\gx(t_*)=-a$.
Further, (\ref{tpi2}) and (\ref{ika2}) again yield $t_*/\tc\to\gf(1)=1$.
Hence, (\ref{gtay}) yields (\ref{tx2}).

Since Theorem~\ref{Tac} yields $T=\tc(1+o_{\mathrm p}(1))$ and
$T<\tcx$ w.h.p.,
(\ref{aahm}) yields $\gx(T)=\gy(T)/(1-\tpi(T))-a =-a-\tc^{1/2}
(\xi+o_{\mathrm p}(1))$;
thus, similarly, (\ref{gtay}) yields, using $\acx-a\gg\tc^{1/2}$,
\begin{eqnarray*}
\tcx-T&=&\bigl(1+o_{\mathrm p}(1)\bigr)\sqrt{\frac{2\tc}{r-1}\bigl(\acx-a-\tc
^{1/2}\bigl(\xi+o_{\mathrm p}(1)\bigr)\bigr)}
\\
&=&\bigl(1+o_{\mathrm p}(1)\bigr)\sqrt{\frac{2\tc}{r-1}(\acx-a)}
=\bigl(1+o_{\mathrm p}(1)\bigr)(\tcx-t_*).
\end{eqnarray*}
Hence, w.h.p., every $t$ between $T$ and $t_*$ satisfies
$\tcx-t=(1+o(1))(\tcx-t_*)$, and then by (\ref{ldgg}),
\[
\gx'(t)
=\bigl(1+o(1)\bigr)\frac{r-1}{\tc}(t-\tcx)
=-\bigl(1+o(1)\bigr)\frac{r-1}{\tc}(\tcx-t_*).
\]
Finally, the mean value theorem yields, similarly to case~\ref{TGsub1},
\[
T-t_*=\frac{\gx(T)-\gx(t_*)}{-(1+o(1))(({r-1})/{\tc})(\tcx-t_*)}
=\frac{\tc^{1/2}(\xi+o_{\mathrm p}(1))}{(2{(r-1)}{\tc\qw}(\acx
-a))^{1/2}},
\]
and the result in~\ref{TGsub2} follows, since $\xi\in N(0,1/r)$
and $\acx\sim\ac=\frac{r-1}r\tc$.
\end{pf*}
\begin{pf*}{Proof of Theorem~\ref{TDG2}}
We use the version described in Section~\ref{SSdynm}
where edges are added at random times.
Let $\hU$ be the time the active set becomes big, that is, the time the
$M$th edge is added. For any given $p$, then $\hU\le p$ if and only
if at time $p$, the active set is big, which is the same as saying
that there is a big active set in $G_{n,p}$.
Fix $x\in\oooo$ and choose $p=\pcx+(r-1)^{1/2} r\qw xa\qqw\pc$. Then
Theorem~\ref{Tpcxx} [with $\gl=(r-1)^{1/2} r\qw x$] yields
${\mathbb P}(\hU\le p)\to\Phi(x)$. In other words,
$(\hU-\pcx)/((r-1)^{1/2} r\qw a\qqw\pc)\dto N(0,1)$, or
%
%
\begin{equation}
\label{hU}
\hU\in\AsN\bigl(\pcx,(r-1)r\qww a\qw\pc^2\bigr).
\end{equation}

Let $N(u)$ be the number of edges at time $u$.
Then $N(0)=0$ and, in
analogy with Lemma~\ref{Lmartin},
$\bigl(N(u)-{n\choose2}u\bigr)/(1-u)$ is a martingale on $[0,1)$. Thus,
Doob's inequality yields, as in the proof of Lemma~\ref{L2}, for any
$u_0\in
\oi$,
$\E\bigl(\sup_{u\le u_0}\bigl|N(u)-{n\choose2} u\bigr|^2\bigr)\le16 {n\choose2}
u_0=O(n^2u_0)$;
cf. \cite{SJII}, Lemma 3.2.
Hence,
\[
\sup_{u\le u_0}\biggl|N(u)- \pmatrix{n\cr2} u\biggr|= \Op(n u_0^{1/2}).
\]
Choosing $u_0=2\pc$, we thus obtain, since $\hU\le2\pc$ w.h.p.,
%
%
\begin{equation}\label{dym}
M=N(\hU)=\pmatrix{n\cr2}\hU+ \Op(n\pc^{1/2}).
\end{equation}

We have by (\ref{pc}), for some constant $c=c(r)$,
\[
\frac{n\pc^{1/2}}{n^2\pc/a^{1/2}}
=\frac{a^{1/2}}{n\pc^{1/2}}=
\frac{a^{1/2} c(na^{r-1})^{1/(2r)}}{n}
=c\biggl(\frac an\biggr)^{1-1/2r}
=o(1).
\]
Consequently, the error term in (\ref{dym}) is $o_{\mathrm p}(n^2\pc
/a^{1/2})$,
and the result follows from (\ref{dym}) and (\ref{hU}).
\end{pf*}

\section{The number of generations}\label{sec10}\label{Sgen}

Let $T_{0}:=0$ and define inductively
%
%
\begin{equation}
\label{tj}
T_{j+1}:=A(T_{j}),\qquad
j\ge0.
\end{equation}
Thus $A(T_0)=A(0)=|\cao|=|\cG_0|$, the size of generation 0 (the
initially active vertices). Further, by our choice
of $u_t$ as one of the oldest unused, active vertices,
$\cZ(T_1)=\cZ(A(0))=\cA(0)=\cG_0$ and
$\cZ(T_2)=\cZ(A(T_1))=\cA(T_1)=\cG_0\cup\cG_1$; in general, by induction,
all vertices in generation $k$ (and earlier)
have been found and declared active at time $T_k$, and they have
been used at time $T_{k+1}=A(T_k)$. In other words,
\[
\bigcup_{j=0}^{k}\cG_j = \cA(T_k)=\cZ(T_{k+1}),\qquad k\ge0.
\]
In particular, the size of generation $k$ equals
\[
|\cG_k|=|\cZ(T_{k+1})\setminus\cZ(T_{k})|=T_{k+1}-T_k,\qquad k\ge0,
\]
and the number of generations $\tau$ defined by
(\ref{tau}) is
\[
\tau=\max\{k\ge0\dvtx T_{k+1}>T_k\}
=\min\{k\ge1\dvtx T_{k+1}=T_k\}-1.
\]

We begin by considering the supercritical case.
We then consider the spread of activation in the bootstrap
percolation
process in three different stages in each of the following subsections.
We first consider the bottleneck when the size is close to $\tc$; we
know that this is where the activation will stop in the critical case,
and in the slightly supercritical case, the activation will grow slowly
here, and this will dominate the total time.
Then follows a~period of doubly exponential growth, and finally, when there
are only
$o(n)$ vertices remaining, it may take some time to sweep up the last
of them.
Recall that Example~\ref{Egen} shows that each of the three phases may
dominate the
two others.

We define, for any $m\le n$,
%
%
\begin{equation}
\tau(m):=\inf\{j\dvtx T_j\ge m\}
\end{equation}
with the interpretation that $\tau(m)=\infty$ if this set of $j$ is empty,
that is, if $m>\Ax=T$.


\subsection{The bottleneck}\label{sec101}\label{SSbottle}

We consider first $\tau(3\tc)$, that is,
the number of generations required to achieve at least $3\tc$
active vertices.
[The constant~3 is chosen for convenience; any constant $>1$
would
give the same result within $O(1)$ w.h.p.]
In the really supercritical case, this is achieved quickly.

%
\begin{proposition}
\label{PGensuper0}
Suppose that $r\ge2$ and $n\qw\ll p \ll n\qwr$. Assume
$a\ge(1+\delta)\ac$ for some $\delta>0$.
Then, w.h.p.
$\tau(3\tc)=O(1)$.
\end{proposition}
\begin{pf}
Lemmas~\ref{L10r} and~\ref{Lmin}
imply that uniformly for $0\le t\le3\tc$, with $x=t/\tc$,
\[
A(t)-t=S(t)-t+a=\biggl(\frac1r x^r -x\biggr)\tc+a +o_{\mathrm p}(\tc)
\ge-\ac+ (1+\delta)\ac+o_{\mathrm p}(\ac)
\]
and thus w.h.p.
\[
A(t)-t\ge\frac{\delta}2\ac
\ge\frac{\delta}4\tc.
\]
Hence, in this range, w.h.p. each generation has size at least $(\delta
/4)\tc$,
and the numbers of generations $\tau(3\tc)$ required to
reach $3\tc$ is thus w.h.p. bounded by $12/\delta$.
\end{pf}

In the slightly supercritical case when $a\sim\ac$, this part may be
a real
bottleneck, however.
We will approximate $A(t)$ by deterministic functions and begin with a
definition:
given a function $F\dvtx\ooo\to\ooo$, define the iterates
$\tf_{j+1}:=F(\tf_j)$ with $\tf_0:=0$.
Thus $T_j=T^A_j$.
%
%
\begin{lemma}
\label{Lgen1}
If $A\le F$, then $T_j\le\tf_j$ for every $j$.
If $A\ge F$, then \mbox{$T_j\ge\tf_j$} for every $j$.
\end{lemma}
\begin{pf}
By induction. Assume, for example, $A\le F$ and $T_j\le\tf_j$. Then,
since $A$ is (weakly) increasing,
\[
T_{j+1}=A(T_j)\le A(\tf_j) \le F(\tf_j)=\tf_{j+1} .
\]
\upqed
\end{pf}

We next prove a deterministic lemma.
%
%
\begin{lemma}\label{Lgenquad}
Let $a,b,t_0>0$, and let $F(t):=t+a+b(t-t_0)^2$.
Assume $a\le t_0$ and $bt_0\le1$.
Let $N$ be the smallest integer such that $\tf_N>2t_0$. Then
\[
N=
\bigl(1+O(bt_0)\bigr)\int_{-t_0}^{t_0}\frac1{a+bx^2}\dd x + O(1).
\]
\end{lemma}
\begin{pf}
Assume that $t\in[0,2t_0]$ and let $\gD:=F(t)-t$.
The assumptions on $a$ and $b$ imply $0< \gD\le a+bt_0^2\le2t_0$.
For $s\in[t,t+\gD]$ we have $|F'(s)-1|=|2b(s-t_0)|\le6bt_0$, and\vadjust{\goodbreak}
thus, by the mean-value theorem,
$|F(s)-s-(F(t)-t)|\le6bt_0 \gD=6bt_0(F(t)-t)$.
Thus, uniformly for such $s$,
$F(s)-s=(F(t)-t)(1+O(bt_0))$ and thus
$(F(t)-t)\qw=(F(s)-s)\qw(1+O(bt_0))$.
Consequently,
\begin{eqnarray*}
1&=&\int_t^{t+\gD}\frac1{F(t)-t}\dd s
=\bigl(1+O(bt_0)\bigr)\int_t^{t+\gD}\frac1{F(s)-s}\dd s
\\
&=&
\bigl(1+O(bt_0)\bigr)\int_t^{t+\gD}\frac1{a+b(s-t_0)^2}\dd s.
\end{eqnarray*}
If $t=\tf_j$, then $t+\gD=F(t)=\tf_{j+1}$. Summing for $j=0,\ldots
,N-1$ we
thus obtain
\[
N =\bigl(1+O(bt_0)\bigr)\int_0^{\tf_{N}}\frac1{a+b(s-t_0)^2}\dd s
\ge\bigl(1+O(bt_0)\bigr)\int_0^{2t_0}\frac1{a+b(s-t_0)^2}\dd s,
\]
and similarly, omitting $j=N-1$,
\[
N-1
\le\bigl(1+O(bt_0)\bigr)\int_0^{2t_0}\frac1{a+b(s-t_0)^2}\dd s.
\]
The result follows, using the change of variable $s=x+t_0$.
\end{pf}
%
%
\begin{proposition}
\label{PGensuper+}
Suppose that $r\ge2$ and $n\qw\ll p \ll n\qwr$. Assume
$a/\ac\to1$ and $a-\acx\gg\sqrt{\ac}$.
Then,
\[
\tau(3\tc)
=\frac{\pi\sqrt2+o_{\mathrm p}(1)}{\sqrt{r-1}}\biggl(\frac{\tc}{a-\acx
}\biggr)^{1/2}.
\]
\end{proposition}
\begin{pf}
By (\ref{pia}) and (\ref{tc1}),
$n\pi(3\tc)=O(\tc)$. Hence,
by Lemmas~\ref{L2} and~\ref{Ltpi}, for $t\le3\tc$,
%
%
\begin{equation}\label{gs}
S(t)=\E S(t)+\Op(\tc^{1/2})
=(n-a)\tpi(t)+\Op(\tc^{1/2}).
\end{equation}

Let $H(t):=a+(n-a)\tpi(t)-t$ and define
$h:=\inf_{t\le3\tc} H(t)$. Let the infimum be attained at $t_*$; it
follows from (\ref{tpi2}) and Lemma~\ref{Lmin} that $t_*\sim\tc$;
cf.~(\ref{jull}).
We have $H(t_*)=h$,
$H'(t_*)=0$ and, uniformly for $t\le3\tc$, using~(\ref{tpii}),
(\ref{tccond}) and (\ref{tc1}),
\begin{eqnarray*}
H''(t)&=&(n-a)\tpi''(t)
=(n-a)p^r\frac{t^{r-1}}{(r-1)!}\biggl(\frac{r-1}{t}-p\biggr)e^{-pt}
\\
&=&
np^r\frac{t^{r-2}}{(r-1)!}\bigl(r-1+o(1)\bigr)
=\biggl(\frac{t}{\tc}\biggr)^{r-2}\frac{r-1+o(1)}{\tc}
\\
&=&
\frac{r-1}{\tc}\biggl(1+o(1)+O\biggl(\frac{|t-\tc|}{\tc}\biggr)\biggr).
\end{eqnarray*}
Hence, by a Taylor expansion, for $0\le t\le3\tc$,
%
%
\begin{equation}\label{gH}
H(t)=h+\frac{r-1}{2\tc}(t-t_*)^2\biggl(1+o(1)+O\biggl(\frac{|t-\tc|}{\tc}\biggr)\biggr).\vadjust{\goodbreak}
\end{equation}
Notice that in the last two formulas, the term $o(1)$ tends to 0 as
$\ntoo$,
uniformly in $t\le3\tc$, and $O(\cdots)$ is uniform in $n$; these
uniformities allow us to combine the two terms in a meaningful way.

On the interval $[0,3\tc]$, $\tpi(t)=o(1)$ by (\ref{tpi2}) and~(\ref{tccond}), and thus
by (\ref{acx})
%
%
\begin{eqnarray}\label{hacx}
h&\sim&\inf_{t\le3\tc}\frac{ H(t)}{1-\tpi(t)}
=\inf_{t\le3\tc}\frac{
a+(n-a)\tpi(t)-t}{1-\tpi(t)}\nonumber\\[-8pt]\\[-8pt]
&=&a+\inf_{t\le3\tc}\frac{ n\tpi(t)-t}{1-\tpi(t)}
=a-\acx.\nonumber
\end{eqnarray}

In particular, by our assumption, $h\gg\ac^{1/2}$.
Consequently, by (\ref{gs}) and (\ref{gH}), for
any fixed small $\eps>0$ and $|t-\tc|\le2\eps\tc$, w.h.p.
%
%
\begin{eqnarray}\label{ga}
A(t)-t&=&a+S(t)-t
=H(t)+o_{\mathrm p}(h)\nonumber\\[-8pt]\\[-8pt]
&=&
\bigl(1+O(\eps)\bigr)\biggl(h+\frac{r-1}{2\tc}(t-t_*)^2\biggr).\nonumber
\end{eqnarray}

Let $t_1:=(1-\eps)t_*\sim(1-\eps)\tc$
and $t_2:=(1+\eps)t_*\sim(1+\eps)\tc$.
For $0\le t\le
t_1$ and $t_2 \le t \le3\tc$,
Lemmas~\ref{L10r} and~\ref{Lmin}
imply
that w.h.p. $A(t)-t\ge c\tc$,
for some constant $c=c(\eps)>0$. The numbers of generations required to
cover the intervals $[0,t_1]$ and $[t_2,3\tc]$
are thus $O(1/c(\eps))$, so $\tau(3\tc)=\nie+O(1/c(\eps))$,
where $\nie$ is the number of generations needed to increase the
size from at least $t_1$ to at least~$t_2$.
To find $\nie$, we may redefine $T_n$ by starting with
$T_0:=t_1$ and iterate as in (\ref{tj}) until we reach $t_2$.
(Note that since $A$ is increasing, if we start with a
larger $T_0$, then every $T_n$ will be larger. Hence, to start with
exactly $t_1$ can only affect $\nie$ by at most 1.)
By (\ref{ga}) and Lemma~\ref{Lgen1}, we may on the interval $[t_1,t_2]$
w.h.p.
obtain upper and lower bounds from $F_\pm(t)=t+(1\pm C\eps)(h+b(t-t_*)^2)$,
where $b:=(r-1)/(2\tc)>0$ and $C$ is some constant.
Let $t_0:=t_*-t_1=\eps t_*>0$.
We have $\acx\sim\ac$ and by assumption
$a\sim\ac$, so by (\ref{hacx}), $h=o(\ac)=o(\tc)$ and thus
$h<t_0/2$ for
large~$n$.
Furthermore, $bt_0=O(\eps t_*/\tc)=O(\eps)$. If $\eps$ is small
enough, we
thus have $bt_0\le1/2$ and,
by a translation $t\mapsto t-t_1$,
Lemma~\ref{Lgenquad}
applies to both~$F_+$ and $F_-$ and
yields,
w.h.p.,
using (\ref{hacx}),
\begin{eqnarray*}
\nie
&=&
\bigl(1+O(\eps)\bigr)\int_{-\eps t_*}^{\eps t_*}\frac{\ddq x}{h+bx^2} +
O(1)
\\
&=&
\bigl(1+O(\eps)\bigr)\intoooo\frac{\ddq x}{h+bx^2} +O\biggl(\frac{1}{b\eps t_*}\biggr)+
O(1)
\\
&=&
\bigl(1+O(\eps)\bigr)\intoooo\frac{\ddq x}{h+bx^2} +O\biggl(\frac{1}{\eps}\biggr)
\\
&=&\bigl(1+O(\eps)\bigr)\frac{\pi}{(hb)^{1/2}}
+O\biggl(\frac{1}{\eps}\biggr)
\\
&=&
\bigl(1+O(\eps)\bigr)\biggl(\frac{2\tc}{(r-1)(a-\acx)}\biggr)^{1/2}\pi
+ O(1/\eps).
\end{eqnarray*}
Since $\tc/(a-\acx)\to\infty$, it follows that for every $\eps>0$, w.h.p.,
with $c'(\eps):=\min(c(\eps),\eps)>0$,
\begin{eqnarray*}
\tau(3\tc)& = &
\nie+O\bigl(1/c(\eps)\bigr)
=\bigl(1+O(\eps)\bigr)\frac{\pi\sqrt2}{\sqrt{r-1}}\biggl(\frac{\tc}{a-\acx}\biggr)^{1/2}
+ O\bigl(1/c'(\eps)\bigr)
\\
&=&
\frac{\pi\sqrt2+O(\eps)}{\sqrt{r-1}}\biggl(\frac{\tc}{a-\acx}\biggr)^{1/2}.
\end{eqnarray*}
The result follows since $\eps>0$ is arbitrary.
\end{pf}
%
%
\begin{remark}
In the critical case $(a-\acx)/\sqrt{\ac}\to y\in\oooo$,
we can use a minor variation of the same argument, now using Lemma~\ref{LG},
where $h\sim a-\ac$ above is replaced by the random
\[
h'= a-\ac+{\tc^{1/2}}W(1/r)+o_{\mathrm p}(\tc^{1/2})
=r\qqw\tc^{1/2}\bigl(y\sqrt{r-1}+\xi+o_{\mathrm p}(1)\bigr),
\]
where $\xi\sim N(0,1)$.
We have $\tau(3\tc)<\infty\iff\Ax\ge3\tc\iff h'>0$; this is w.h.p.
equivalent to $y\sqrt{r-1}+\xi>0$.
[This\vspace*{1pt} thus happens with probability $\Phi((r-1)^{1/2} y)+o(1)$,
as stated in Theorem~\ref{Tac}\ref{Tac0}.]
The argument above then shows that
conditioned on $\tau(3\tc)<\infty$ (i.e., on $\Ax\ge3\tc$),
\[
\tau(3\tc)/\tc\qqqq\dto
\biggl(\frac{2^{1/2}\pi r\qqqq}{\sqrt{r-1}}\bigl(\xi+y\sqrt{r-1}\bigr)\qqw\Bigm
|\xi+y\sqrt{r-1}>0\biggr).
\]
In particular, then $\tau(3\tc)=\Theta_p(\tc^{1/4})$.

Note that in the supercritical case in Proposition~\ref{PGensuper+},
the time
$\tau(3\tc)$ is always smaller than $\tc\qqqq$,
but that it approaches the order $\tc\qqqq$ when $a-\acx$ grows
only a little faster than the critical value $\ac^{1/2}$.
Hence, we can say that
the worst possible number of generations to pass the bottleneck at
$\tc$ is of the order $\tc\qqqq$.
\end{remark}

\subsection{The doubly exponential growth}\label{sec102}

We next consider the growth from size $3\tc$ up to $1/p$. We will show that
in this range, the growth is doubly exponential. Again, we approximate
$A(t)$ by deterministic functions.

Define for any $\delta\in\bbR$ [cf. (\ref{tc1})],
%
%
\begin{equation}
\label{Fd}
F_{\delta}(t):= n \frac{(tp)^r}{r!} (1+\delta)
= \biggl(\frac{t}{\tc}\biggr)^{r-1}\frac{t}{r}(1+\delta).
\end{equation}
%

\begin{lemma}
\label{LGen2a}
For every $\delta>0$,
there are positive constants $\eps$ and $K$ such that
w.h.p.
$F_{-\delta}(t) \leq A(t) \leq F_{\delta}(t)$ \mbox{ for all }
$t\in[ K(\tc+ a),\eps/p ]$.
\end{lemma}
\begin{pf}
By (\ref{pia}) and (\ref{Fd}), for $K\tc\le t\le\eps/p$ (with
$\eps\le1$),
if $n$ is large enough so $\tc\ge1$,
\[
\pi(t) = \frac{(tp)^r}{r!} \bigl(1+O(\eps+K\qw)\bigr)=
\frac{1}{n} F_0(t)\bigl(1+O(\eps+K\qw)\bigr).
\]
We may thus choose $\eps$ and $K$ such that for all such $t$ (and
large $n$)
%
%
\begin{equation}\label{D3}
F_{-\delta/4}(t)\leq n\pi(t) \leq F_{\delta/4}(t).
\end{equation}

For $t\ge K(\tc+a)$, (\ref{Fd}) implies
\[
F_{0}(t)= \biggl(\frac{t}{\tc}\biggr)^{r-1}\frac{t}{r}
\ge K^{r} \frac ar,
\]
so choosing $K$ large enough, we have $a\le(\delta/4)F_0(t)$ for all
$t\in[K(\tc+a),\break\eps/n]$,
and thus by (\ref{D3})
\[
F_{-\delta/4}(t)-a\leq
\E S(t)=( n-a)\pi(t)
\leq F_{\delta/4}(t)
\le F_{\delta/2}(t)-a.
\]
Hence, by Chebyshev's inequality, using (\ref{svar}) and (\ref{D3}),
%
%
\begin{eqnarray}\label{D3b}
&&{\mathbb P}\{A(t)\notin[F_{-3\delta/4}(t),F_{3\delta/4}(t)]\}\nonumber\\
&&\qquad=
{\mathbb P}\{S(t)\notin[F_{-3\delta/4}(t)-a,F_{3\delta/4}(t)-a]\}
\\
&&\qquad\le
\frac{n\pi(t)}{(\delta F_0(t)/4)^2}
\le\frac{F_{\delta/4}(t)}{(\delta F_0(t)/4)^2}
=\frac{16(1+\delta/4)}{\delta^2 F_0(t)}.
\nonumber
\end{eqnarray}

Define
$t_j:=(1+\delta/5)^{j/r}K(\tc+a)$. Then,
(\ref{D3b}) and (\ref{Fd}) show that, assuming as we may $\delta\le1$,
\begin{eqnarray*}
\sum_{j\ge0\dvtx t_j\le\eps/p}
{\mathbb P}\{A(t_j)\notin[F_{-3\delta/4}(t_j),F_{3\delta/4}(t_j)]\}
&\le&
\sum_{j\ge0}\frac{20}{\delta^2 F_0(t_j)}
\\
& = &
\sum_{j\ge0}\frac{20}{\delta^2 F_0(t_0)}(1+\delta/5)^{-j}
\\
&=&
\frac{100(1+\delta/5)}{\delta^3 F_0(t_0)}
\to0,
\end{eqnarray*}
since, using (\ref{Fd}) again and (\ref{tccond}),
\[
F_0(t_0)=F_0\bigl(K(\tc+a)\bigr)\ge F_0(\tc)=\frac{\tc}r\to\infty.
\]
Consequently, w.h.p.
$A(t_j)\in[F_{-3\delta/4}(t_j),F_{3\delta/4}(t_j)]$
for all $j\ge0$ with $t_j\le\eps/p$.
However, if $t_j\le t\le t_{j+1}$, then $F_0(t_j)\le F_0(t)\le
F_0(t_{j+1})=(1+\delta/5)F_0(t_j)$, and it follows that, since both
$A(t)$ and
$F_0(t)$ are monotone, w.h.p.
\[
(1+\delta/5)\qw F_{-3\delta/4}(t)
\le A(t)
\le(1+\delta/5) F_{3\delta/4}(t)
\]
for all $t\in[K(\tc+a),(1+\delta/5)^{-1/r}\eps/p]$,
which,
provided $\delta$ is small and $\eps$ is replaced by $\eps/2$, say,
yields the result.
\end{pf}
%
%
\begin{proposition}\label{PGen2}
Suppose that $r\ge2$ and $n\qw\ll p \ll n\qwr$.
Then w.h.p., when $\Ax\ge3\tc$,
\[
\tau(1/p)-\tau(3\tc)
=\frac{1}{\log r}
\biggl(\log\log(np) - \log_+\log\frac a{\ac}\biggr)+O(1).
\]
\end{proposition}
\begin{pf}
Choose a fixed $0<\delta<1$, and choose $\eps$ and $K$ as in Lemma
\ref{LGen2a}.
(In this proof, we do not have to let $\delta\to0$, so we can take
$\delta=1/2$,
say.)
First, $\tau(K(\tc+a))-\tau(3\tc)$, the number of generations from
$3\tc$ to $K(\tc+a)$, is w.h.p. $O(1)$. Indeed, after $\tau(3\tc)$
generations
we have
at least $\max(3\tc,a)$ active vertices, and
in each of the following generations until well beyond $K(\tc+a)$,
the number is w.h.p. multiplied by at least $1.3$, say,
by the proof of Lemma~\ref{Lbulk} or by Lemmas~\ref{LS1},~\ref{Ltpi}
and~\ref{L2}.
Similarly, $\tau(1/p)-\tau(\eps/p)\le1$ w.h.p., arguing as in
Case 3
of the
proof of Lemma~\ref{Lbulk}.

Consequently it suffices to consider $\tau(\eps/p)-\tau(K(\tc+a))$.
We define iterates $T_{j}^{F_{\delta}}$ as in Section~\ref{SSbottle} by
$T_{j+1}^{F_{\delta}}:=F_\delta(T_{j}^{F_{\delta}})$, $j\ge0$, but
now starting with
$T_{0}^{F_{\delta}}:=K(\tc+a)$. Further, let
%
%
\begin{equation}
\label{Nd}
N_\delta:=\min\{j\ge0\dvtx T_{j}^{F_{\delta}}\ge\eps/p\}.
\end{equation}
By Lemma~\ref{LGen2a} we may assume that
$F_{-\delta}(t) \leq A(t) \leq F_{\delta}(t)$ for all
$t\in[ K(\tc+ a),\eps/p ]$, and then,
by induction as in\vspace*{1pt} Lemma~\ref{Lgen1},
$T_{j}^{F_{-\delta}}\le T_{j+\tau(K(\tc+a))} \le T_{j+1}^{F_{\delta
}}$ for all $j\ge
0$ with
$T_{j-1+\tau(K(\tc+a))}\le\eps/p$. Consequently, w.h.p.
%
%
\begin{equation}\label{t8}
N_{-\delta}\ge\tau(\eps/p)-\tau\bigl(K(\tc+a)\bigr)\ge N_\delta-1.
\end{equation}

To find $N_\delta$, rewrite (\ref{Fd}) as
\[
\frac{F_{\delta}(t)}{\bd\tc} =\biggl(\frac{t}{\bd\tc}\biggr)^r,
\]
where $\bd:=(r/(1+\delta))^{1/(r-1)}$. Iterating we see that, for
$j\ge0$,
\[
\frac{T_{j}^{F_{\delta}}}{\bd\tc} =\biggl(\frac{T_{0}^{F_{\delta
}}}{\bd\tc}\biggr)^{r^j}
=\biggl(\frac{K(\tc+a)}{\bd\tc}\biggr)^{r^j}
\]
and thus
\[
\log\biggl(\frac{T_{j}^{F_{\delta}}}{\bd\tc}\biggr) =r^j\log\biggl(\frac{K(\tc
+a)}{\bd\tc}\biggr)
\]
and
\[
j\log r
= \log\log\biggl(\frac{T_{j}^{F_{\delta}}}{\bd\tc}\biggr)
-\log\log\biggl(\frac{K(\tc+a)}{\bd\tc}\biggr).
\]
Consequently,
%
%
\begin{equation}\label{jb}
N_\delta
= \biggl\lceil\biggl(\log\log\biggl(\frac{\eps/p}{\bd\tc}\biggr) -\log\log\biggl(\frac
{K(\tc+a)}{\bd\tc}\biggr)\biggr)\Big/\log
r\biggr\rceil.\vadjust{\goodbreak}
\end{equation}
In order to simplify this, note that, using (\ref{tc}),
%
%
\begin{equation}
\log\biggl(\frac{\eps/p}{\bd\tc}\biggr)=
\log\biggl(\frac{1}{p\tc}\biggr)+O(1)
=\frac{1}{r-1}\log(np)+O(1)
\end{equation}
and thus
%
%
\begin{equation}\label{acu}
\log\log\biggl(\frac{\eps/p}{\bd\tc}\biggr)
=\log\log(np)+O(1).
\end{equation}
Further, we may assume that $a\ge\ac/2\ge\tc/4$, since otherwise
the process
is subcritical and $\Ax<3\tc$ w.h.p. by Theorem~\ref{T1}. Hence,
$\log(K(\tc+a))=\log a+O(1)$
and thus, since also $\log(\bd\tc)=\log\ac+O(1)$,
%
%
\begin{equation}\label{k1}
\log\biggl(\frac{K(\tc+a)}{\bd\tc}\biggr)=\log a-\log\ac+O(1)
=\log\frac a{\ac}+O(1).
\end{equation}
We may assume that $K\ge e\bd$, so $\log(K(\tc+a)/(\bd\tc))\ge
1$, and
then (\ref{k1}) yields
%
%
\begin{equation}\label{kll}
\log\log\biggl(\frac{K(\tc+a)}{\bd\tc}\biggr)
=\log_+\log\frac a{\ac}+O(1).
\end{equation}
Finally, (\ref{jb}), (\ref{acu}) and (\ref{kll}) yield
\[
N_\delta\log r
=\log\log(np) -\log_+\log\frac a{\ac}+O(1).
\]
Note that the right-hand side depends on $\delta$ only in the error
term $O(1)$.
Hence, we
have the same result for $N_{-\delta}$, and the result follows by
(\ref{t8})
and the comments at the beginning of the proof.
\end{pf}

\subsection{The final stage}\label{sec103}

We finally consider the evolution after $1/p$ vertices have become active.
We let, as in Section~\ref{Sthreshold}, $\bcx:=\bc\go(n)$ where
$\go
(n)\to
\infty$
slowly; we assume that $\bcx\ll1/p$ [which is possible since $p\bc
\to
0$ by
(\ref{tccond})].
By Remark~\ref{Rbulk}, $\tau(n-\bcx)\le\tau(1/p)+3$ w.h.p., so it
suffices to
consider the
evolution when less than $\bcx$ vertices remain.

Let $\cF_t:=\gs\{I_i(s)\dvtx1\le i\le n, 1\le s\le t\}$ be the
$\gs$-field
describing the evolution up to time $t$.
%
%
\begin{lemma}
\label{Lnassjo}
For any $t$ and $u$ with $0\le t\le t+u\le n$, the conditional
distribution of $A(t+u)-A(t)=S(t+u)-S(t)$ given $\cF_t$ is
$\Bin(n-A(t),\pi(t;u))$, where
%
%
\begin{equation}
\label{pitu}
\pi(t;u):=\frac{\pi(t+u)-\pi(t)}{1-\pi(t)}
.
\end{equation}

If further $n-\bcx\le t\le t+u\le n$, then,
uniformly in all such $t$ and $u$,
%
%
\begin{equation}
\label{pitu1}
\pi(t;u)=p u\bigl(1+o(1)\bigr).
\end{equation}
\end{lemma}
\begin{pf}
Conditioned on $\cF_t$, $A(t)$ is a given number, and of the $n-a$
summands in (\ref{st}), $n-a-S(t)=n-A(t)$ are zero. For any of these
terms, the probability that it changes from $0$ at time $t$ to 1 at time
$t+u$ is, by~(\ref{pi}),
\[
{\mathbb P}(Y_i\le t+u\mid Y_i>t)
=
\frac{{\mathbb P}(t<Y_i\le t+u)}{{\mathbb P}(Y_i>t)}
=\frac{\pi(t+u)-\pi(t)}{1-\pi(t)}
=\pi(t;u).
\]
Hence, the conditional distribution of $S(t+u)-S(t)$ is
$\Bin(n-A(t),\pi(t;u))$.

To see the approximation (\ref{pitu1}), note first that for $n-\bcx
\le
t\le n$,
since we assume $p\bcx\to0$, we have $\bcx\ll1/p\ll n$ so
$t\sim n$. Hence, using again $p\bcx\to0$ and recalling the notation
$\bcq$
from (\ref{bc2a}),
%
%
\begin{eqnarray}\label{pitt}
\pi(t+1)-\pi(t)&=&{\mathbb P}(\yix=t+1)
= \pmatrix{t\cr r-1}p^r(1-p)^{t+1-r}
\nonumber\\[-8pt]\\[-8pt]
& \sim&\frac{n^{r-1}}{(r-1)!}p^r(1-p)^n
=\frac{p\bcq}n.
\nonumber
\end{eqnarray}
Furthermore [cf. (\ref{bc2a})], still for $n-\bcx\le t\le n$,
%
%
\begin{eqnarray}\label{pit}
1-\pi(t)&=&{\mathbb P}\bigl(\Bin(t,p)\le r-1\bigr)
\sim{\mathbb P}\bigl(\Bin(t,p)= r-1\bigr)
\nonumber\\[-8pt]\\[-8pt]
&\sim& \frac{n^{r-1}}{(r-1)!}p^{r-1}(1-p)^n =\frac{\bcq}n.
\nonumber
\end{eqnarray}
Consequently, $\pi(t+u)-\pi(t)=(1+o(1)) up\bcq/n$ and
\[
\pi(t;u)
=\frac{\pi(t+u)-\pi(t)}{1-\pi(t)}
=\bigl(1+o(1)\bigr)\frac{up\bcq/n}{\bcq/n}
=\bigl(1+o(1)\bigr) up.
\]
\upqed
\end{pf}
%
%
\begin{lemma}
\label{LG3}
Suppose that $r\ge2$, $n\qw\ll p \ll n\qwr$ and $a=o(n)$.
If $\bc\to\infty$ and $n-\bcx\le t\le n$, then $A(t)=n-\bc
(1+o_{\mathrm p}(1))$; in
particular,
$n-A(t)<2\bc$ w.h.p.
\end{lemma}
\begin{pf}
We have, using (\ref{es}) and (\ref{pit}), since $\bc\to\infty$ implies
$\bcq\sim\bc$,
\[
\E\bigl(n-A(t)\bigr)
=n-a-\E S(t)
=(n-a)\bigl(1-\pi(t)\bigr)
\sim(n-a)\frac{\bcq}n\sim\bc
\]
and similarly, using (\ref{svar}),
\[
\Var\bigl(n-A(t)\bigr)=\Var S(t)\le
(n-a)\bigl(1-\pi(t)\bigr)
\sim\bc.
\]
Thus, by Chebyshev's inequality, since $\bc\to\infty$,
\[
n-A(t)=\bigl(1+o(1)\bigr)\bc+\Op(\bc^{1/2})
=\bigl(1+o_{\mathrm p}(1)\bigr)\bc.
\]
\upqed
\end{pf}
%
%
\begin{proposition}\label{PG3}
Suppose that $r\ge2$, $n\qw\ll p \ll n\qwr$ and $a=o(n)$.
Then, when $\Ax\ge3\tc$,
\[
\tau-\tau(1/p)
=\bigl(1+o(1)\bigr)\frac{\log n}{np}
+\Op(1).\vadjust{\goodbreak}
\]
In particular, if further $p\ge c\log(n)/n$ for some $n\ge0$, then
$\tau-\tau(1/p)=\Op(1)$.

Furthermore, when $\Ax=n$, w.h.p. $\tau-\tau(1/p)\le3$.
\end{proposition}
\begin{pf}
By Remark~\ref{Rbulk}, after $\tau(1/p)+3$ generations, the active
size is
$T_{\tau(1/p)+3}\ge n-\bcx$ w.h.p.

If $\bc\to0$, we can choose $\bcx=1/2$, so w.h.p. $T_{\tau
(1/p)+3}=n$ and
\mbox{$\tau\le\tau(1/p)+3$}.

More generally, if $\bc=O(1)$, we have
by (\ref{pitt}),
\begin{eqnarray*}
\E\bigl(S(n)-S(n-\bcx)\bigr)
&\le& n\bigl(\pi(n)-\pi(n-\bcx)\bigr)
\sim n\bcx\frac{p\bcq}{n}\\
&=& p\bcx\bc= O(p\bcx)=o(1).
\end{eqnarray*}
Hence, w.h.p. $S(n)=S(n-\bcx)$, which means that no further activations
occur after $n-\bcx$. Consequently, in this case too, w.h.p.
$\tau=\tau(n-\bcx)\le\tau(1/p)+3$.
In particular, this proves that $\tau\le\tau(1/p)+3$ w.h.p.
when $\Ax=n$, since w.h.p. $\Ax<n$
if $\bc\to\infty$ by Theorem~\ref{T2}.

Further, when $\bc=O(1)$, (\ref{bclim}) implies that $np\ge\log n$
for large
$n$, so $\log n/(np)\le1$, and the result holds in this case.

Now assume that $\bc\to\infty$.
For convenience, we modify the counting of generations and start at
$t=n-\bcx$, regarding the active but unused vertices at $n-\bcx$ as
``generation 0.'' (We may assume that $\bcx$ is an integer.)
Thus define, recursively,
\begin{eqnarray*}
T_0'&:=&n-\bcx,\\
T_{j+1}'&:=&A(T_j'),\qquad j\ge0,\\
\gD_j&:=&T'_{j+1}-T'_j=A(T'_j)-T'_j,\\
\tauw&:=&\max\{j\ge0\dvtx\gD_j>0\}.
\end{eqnarray*}
Since w.h.p. $T_{\tau(1/p)-1}\le\max(1/p,a)< n-\bcx\le T_{\tau(1/p)+3}$,
it follows by induction that $T_{\tau(1/p)-1+j}\le T'_j\le T_{\tau(1/p)+3+j}$,
$j\ge0$, and thus w.h.p.
%
%
\begin{equation}\label{cec}
\tauw+\tau(1/p)-1\le\tau\le\tauw+\tau(1/p)+3.
\end{equation}
Consequently, it suffices to estimate $\tauw$.

By Lemma~\ref{Lnassjo}, conditioned on $\cF_{T'_j}$ [i.e., on $T'_j$
and the
evolution up to~$T'_j$, which in particular specifies $A(T'_j)$],
for large $n$,
\[
\E(\gD_{j+1}\mid\cF_{T'_j})
=\bigl(n-A(T'_j)\bigr)\pi(T'_j;\gD_j)
\le
\bigl(n-A(T'_0)\bigr)2p\gD_j
\]
and thus, by induction, since $\gD_0\le n-T'_0=\bcx$,
%
%
\begin{equation}\label{lisa}
\E(\gD_{j}\mid\cF_{T'_0})
\le
\bigl(2\bigl(n-A(T'_0)\bigr)p\bigr)^j \gD_0
\le
\bigl(2\bigl(n-A(T'_0)\bigr)p\bigr)^j \bcx.
\end{equation}
Further, Lemma~\ref{LG3} yields
$ n-A(T'_0)=n-A(n-\bcx)<2\bc$ w.h.p.
Consequently, (\ref{lisa}) implies, w.h.p. for all $j\ge0$ (simultaneously),
%
%
\begin{equation}\label{jan}
\E(\gD_{j}\mid\cF_{T'_0})
\le(4p\bc)^j \bcx.\vadjust{\goodbreak}
\end{equation}
Recall that $p\bc\to0$ by (\ref{tccond}),
so we may assume $4p\bc<1$.
If $j$ is chosen such that $(4p\bc)^j\bcx\to0$,
then (\ref{jan}) implies that w.h.p.
$\gD_j=0$ and thus $\tauw<j$. Hence, for any $\gox=\gox(n)\to
\infty$, w.h.p.
\[
\tauw\le\frac{\log\bcx}{|{\log}(p\bc)+\log4|}+\gox(n),
\]
which is another way of saying \cite{SJN6}, Lemma 3,
%
%
\begin{equation}\label{g4}
\tauw\le\frac{\log\bcx}{|{\log}(p\bc)+\log4|}+\Op(1)
=\frac{\log\bcx}{|{\log}(p\bc)|}\bigl(1+o(1)\bigr)+\Op(1).\hspace*{-28pt}
\end{equation}

For a lower bound, fix $\eps$ with $0<\eps<1$, and
define the deterministic numbers $\gD_{j}^{-}$ by
%
%
\begin{equation}\label{siv}
\gD_{j}^{-}:=(1-\eps)^{j+1}(p\bc)^j\bcx.
\end{equation}
Let $\goxx:=1/(p\bc)\to\infty$. We claim that w.h.p.
%
%
\begin{equation}\label{manne}
\gD_j\ge\gD_{j}^{-} \qquad\mbox{for all $j\ge0$ such that } \gD
_{j}^{-}\ge
\goxx.
\end{equation}
By our assumption $4p\bc<1$, we have
$\gD_{j+1}^{-}/\gD_{j}^{-}<1/4$, so $\gD_{j}^{-}\to0$ geometrically fast.

By Lemma~\ref{LG3} and $\bc/\bcx\to0$, w.h.p.
\[
\gD_0=A(T'_0)-(n-\bcx)=A(T'_0)-n+\bcx\ge\bcx-2\bc\ge(1-\eps
)\bcx=\gD_{0}^{-} ,
\]
so (\ref{manne}) holds w.h.p. for $j=0$.

Say that $j\ge0$ is \textit{good} if $\gD_j\ge\gD_{j}^{-}$ and
\textit
{fat} if
$A(T'_j)>n-(1-\eps/4)\bc$.
Let $j\ge0$.
At time $T'_j$ we have $A(T'_j)-T'_j=\gD_j$ active but unused vertices.
Further, by Lemma~\ref{Lnassjo} we have, conditioned on $\cF_{T'_j}$ (which
specifies both~$T'_j$ and $\gD_j$),
\[
\gD_{j+1}=T'_{j+2}-T'_{j+1}
=A(T'_j+\gD_j)-A(T'_j)\in\Bin\bigl(n-A(T'_j),\pi(T'_j;\gD_j)\bigr).
\]
By Lemma~\ref{Lnassjo},
$\pi(T'_j;\gD_j)=p\gD_j(1+o(1))\ge p\gD_j(1-\eps/4)$ for $n$ large,
so if~$j$ is good but not fat,
\begin{eqnarray*}
\E(\gD_{j+1}\mid\cF_{T'_j})
&=&\bigl(n-A(T_j')\bigr)\pi(T'_j;\gD_j)
\ge(1-\eps/4)^2\bc p\gD_j
\\
&\ge& (1-\eps/2)\bc p\gD_{j}^{-}
\ge(1+\eps/2)\gD_{j+1}^{-}
\end{eqnarray*}
and Chebyshev's inequality yields, since $x\mapsto x/(x-a)^2$ is decreasing
for $x>a$,
\begin{eqnarray*}
{\mathbb P}(\gD_{j+1}<\gD_{j+1}^{-}\mid\cF_{T'_j})
&\le&\frac{\Var(\gD_{j+1}\mid\cF_{T'_j})}
{(\E(\gD_{j+1}\mid\cF_{T'_j})-\gD_{j+1}^{-})^2}
\\
&\le&
\frac{\E(\gD_{j+1}\mid\cF_{T'_j})}
{(\E(\gD_{j+1}\mid\cF_{T'_j})-\gD_{j+1}^{-})^2}
\\
&\le&
\frac{(1+\eps/2)\gD_{j+1}^{-}}
{(\eps\gD_{j+1}^{-}/2)^2}
=O\biggl(\frac{1}{\gD_{j+1}^{-}}\biggr).
\end{eqnarray*}
Say that $j$ is \textit{bad} if $j$ is not good and that $j$ \textit
{fails} if
$j$ is fat or bad.
Then, by stopping at the
first $j$ that fails we see that
\begin{eqnarray*}
{\mathbb P}(\mbox{some $j\le\goxx$ fails})
&\le&
{\mathbb P}(\mbox{some $j\le\goxx$ is fat}) +{\mathbb P}(\mbox{$0$
is bad})
\\
&&{}
+
\sum_{j>0\dvtx\gD_{j}^{-}\ge\goxx}
{\mathbb P}(\mbox{$j$ is bad}\mid\mbox{$j-1$ is good and not fat})
\\
&\le&
{\mathbb P}\bigl(A(n) > n - (1 -\eps/4)\bc\bigr) + o(1)\\
&&{} + \sum
_{j\dvtx\gD_{j}^{-}\ge\goxx} O\biggl(\frac{1}{\gD_{j}^{-}}\biggr)
\\
&=& o(1),
\end{eqnarray*}
since $A(n)<n-(1-\eps/4)\bc$ w.h.p. by Lemma~\ref{LG3} and the
final sum is
$O(1/\goxx)=o(1)$ because the terms $1/\gD_{j}^{-}$ increase geometrically,
so the
sum
is dominated by its largest (and last) term.

We have shown that w.h.p., if $\gD_{j}^{-}\ge\goxx$,
then $\gD_j\ge\gD_{j}^{-}>0$ and thus
$\tauw\ge j$. Hence, by (\ref{manne}) and (\ref{siv}), w.h.p.
%
%
\begin{equation}\label{g5}
\tauw\ge\biggl\lfloor\frac{\log((1-\eps)\bcx/\goxx)}{|{\log}((1-\eps
)p\bc)|}\biggr\rfloor
=\frac{\log\bcx}{|{\log}(p\bc)|}\bigl(1+o(1)\bigr)+O(1).
\end{equation}

Combining the upper bound (\ref{g4}) and the lower bound (\ref{g5}),
we find
%
%
\begin{equation}\label{g6}
\tauw
=\frac{\log\bcx}{|{\log}(p\bc)|}\bigl(1+o(1)\bigr)+\Op(1).
\end{equation}
By (\ref{bc}), $\log(p\bc)=-(np-r\log(np)+O(1))$ and
\[
\log n\ge\log\bcx\ge\log\bc\ge\log n-pn-O(1).
\]
Hence, finally (\ref{g6}) yields
\[
\tauw= \frac{\log n+O(np)}{np-r\log(np)+O(1)}\bigl(1+o(1)\bigr)+\Op(1)
= \frac{\log n}{np}\bigl(1+o(1)\bigr)+\Op(1).
\]
The result now follows from (\ref{cec}).
\end{pf}

\section{\texorpdfstring{Proofs of Theorems \protect\ref{pn=cn}, \protect\ref{sqrtn}, \protect\ref{sqrtnn}}
{Proofs of Theorems 5.2, 5.6, 5.8}}\label{sec11}
\label{Spf+}\label{Slast}

We prove in this section Theorems~\ref{pn=cn},~\ref{sqrtn} and~\ref{sqrtnn}
related to the boundary cases. We consider first the case $p\sim
c/n$.\vadjust{\goodbreak}

\begin{pf*}{Proof of Lemma~\ref{Lf}}
By the implicit function theorem, at least locally, the root $\xo(\gth
)$ is
smooth except at points where
%
%
\begin{equation}\label{critical}
f(x,c,\gth)=\ddx f(x,c,\gth)=0.
\end{equation}
We begin
by studying such \textit{critical points}.

Let $g(y):={\mathbb P}(\Po(y)\le r-1)=1-\psi(y)$; cf. (\ref{tpi}).
Differentiations yield
%
%
\begin{eqnarray}
\label{cg}
g'(y)&=&-{\mathbb P}\bigl(\Po(y)=r-1\bigr)=-\frac{y^{r-1}}{(r-1)!}e^{-y},
\\
\label{cgg}
g''(y)&=&\biggl(\frac{r-1}y-1\biggr)g'(y)=\frac{r-1-y}{y}g'(y).
\end{eqnarray}
We have [see (\ref{f1b1})] $f(x,c,\gth)=1-x-(1-\gth)g(cx)$ and thus
$\ddx f(x,c,\gth)=-1-c(1-\gth)g'(cx)$.
Hence, (\ref{critical}) holds if and only if
\[
\cases{
(1-\gth)g(cx)=1-x,\cr
c(1-\gth)g'(cx)=-1,}
\]
which imply $g(cx)=-c(1-x)g'(cx)$ and thus
%
%
\begin{equation}
\label{ch}
c=cx-\frac{g(cx)}{g'(cx)}.
\end{equation}
Let $h(y):=y-g(y)/g'(y)$, $y>0$, so (\ref{ch}) says $c=h(cx)$.
Then, by (\ref{cgg}),
\[
h'(y)=1-\frac{g'(y)}{g'(y)}+\frac{g(y)g''(y)}{g'(y)^2}
=\frac{r-1-y}{y}\cdot\frac{g(y)}{g'(y)}.
\]
Since $g(y)>0$ and $g'(y)<0$ for $y>0$, $h$ has a global minimum at $y=r-1$,
and the minimum value is
\[
\min_{y>0} h(y)=h(r-1)=r-1-\frac{g(r-1)}{g'(r-1)}
=r+\frac{{\mathbb P}(\Po(r-1)\le r-2)}{{\mathbb P}(\Po(r-1)= r-1)}
=\ccr.
\]
Furthermore, $h(y)>y\to\infty$ as $y\to\infty$, and $h(y)\to\infty
$ as
$y\to0$ too, because then $g(y)\to1$ and $g'(y)\to0$.

Consequently, if $0\le c<\ccr$, then (\ref{ch}) has no solution
$x>0$, and
thus there is no critical point.
If $c=\ccr$, there is exactly one $x>0$ satisfying (\ref{ch}) [viz.,
$x=(r-1)/\ccr$], and if $c>\ccr$, there are two. Since (\ref{ch}) implies
$c>cx$, these roots are in $(0,1)$.

To complete the proof, it is perhaps simplest to rewrite (\ref{f10}) as
$\gth=\vartheta(x)$, with
%
%
\begin{equation}
\label{gthx}
\vartheta(x):=1-(1-x)/g(cx).
\end{equation}
Since $g(y)>0$ for $y\ge0$,
$\gthx$ is a smooth function on $\oi$, with $\gthx(0)=0$ and $\gthx(1)=1$.
Moreover, $f(x,c,\gth)=g(cx)(\gth-\vartheta(x))$, which implies
that
\[
f(x,c,\gth)=\ddx f(x,c,\gth)=0
\quad\iff\quad
\gth=\vartheta(x)\mbox{ and } \gthx'(x)=0.
\]
Consequently, by the results above, if $c<\ccr$, then $\gthx'\neq0$ so
$\gthx'(x)>0$ for $x\ge0$. In this case, $\gthx$ is strictly
increasing and
thus a bijection $\oi\to\oi$, and
$\xo$ is its inverse.

If $c=\ccr$, then $\gthx'=0$ only at a single point, and it follows again
that $\gthx$ is a strictly increasing function and $\xo$ is its inverse.

If $c>\ccr$, then $\gthx'(x)=0$ at two values $x_1$ and $x_2$ with
\mbox{$0<x_1<x_2<1$} and $cx_1<r-1<cx_2$. It can be seen, for example, using
(\ref{cgg}),
that
$\gthx''(x_1)<0<\gthx''(x_2)$, and thus $\gthx$ is \textit
{decreasing} on the
interval $[x_1,x_2]$. The result follows, with $\gthc=\gthx(x_1)$,
$\gthcq=\max(\gthx(x_2),0)$ and
$\xo(\gthc)=x_1$. [Note that $\gthx(x_2)=\min_{x\in\oi}\vartheta(x)<0$
if $c$ is
large enough.]
\end{pf*}
%
%
\begin{remark}
\label{Rgthc}
If $c>\ccr$, then thus $\xo(\gthc)=x_1$ is the smallest root of
$\gthx
'(x)=0$, or
equivalently $x_1=y_1/c$ where $y_1$ is the smallest root of $h(y)=c$;
further, $\gthc=\gthx(x_1)$ while $\xo(\gthc+)$ is the other root of
$\gthx(x)=\gthx(x_1)$.

If $c=\ccr$, we have $y_1=r-1$ and thus $x_1=(r-1)/\ccr$ and, using
(\ref{gthx}) and (\ref{ccr}),
\begin{eqnarray*}
\gthc(\ccr) & = & \gthx\biggl(\frac{r-1}{\ccr}\biggr) =1-\frac{1-(r-1)/\ccr}{g(r-1)}
\\
& = & 1-\frac{1}{r{\mathbb P}(\Po(r-1)=r-1)+{\mathbb P}(\Po(r-1)\le r-2)}.
\end{eqnarray*}
\end{remark}

For $c>\ccr$, the two roots $x_1(c)$ and $x_2(c)$ of $\gthx'(x)=0$ are
smooth functions of~$c$, and thus
\[
\frac{\ddq\gthc}{\ddq c}
= \frac{\ddq}{\ddq c}\gthx(x_1(c))
=\frac{\partial\gthx}{\partial c}(x_1(c))
+\frac{\partial\gthx}{\partial x}(x_1(c))x_1'(c)
=\frac{\partial\gthx}{\partial c}(x_1(c))
<0,
\]
where the last inequality follows from (\ref{gthx}), and similarly
$\ddq\gthcq/\ddq c<0$. Hence, $\gthc(c)$ and $\gthcq(c)$ are decreasing
functions of $c$, as claimed in Remark~\ref{Rgthcc}.
%
%
\begin{lemma}\label{ll10r}
Suppose that $r\ge2$, $p=O(1/n)$ and $tp=o(1)$. Then\break
$S(t)=o_{\mathrm p}(t)$.
\end{lemma}
\begin{pf}
We may assume $1\le t\le1/p$. [Note that $S(t)=0$ for $t<r$.]
Then $\pi(t)=O(t^rp^r)=o(tp)$ by (\ref{pia}), and thus
the expected number of activated vertices is
$\E S(t)=(n-a)\pi(t)=o(npt)=o(t)$.
\end{pf}
\begin{pf*}{Proof of Theorem~\ref{pn=cn}}
First,
in~\ref{p=cnaon} and~\ref{pn=cn=0}, $ap\to\theta c=0$.
Let $\eps>0$.
Taking $t=(1+\eps)a$ in
Lemma~\ref{ll10r}, we find
w.h.p. $S((1+\eps)a)<\eps a$ and thus
\[
A\bigl((1+\eps)a\bigr)=a+S\bigl((1+\eps)a\bigr)
<{(1+\eps)a},
\]
whence $\Ax=T<{(1+\eps)a}$. Consequently, $1\le\Ax/a<1+\eps$ w.h.p.,
proving~\ref{p=cnaon} and~\ref{pn=cn=0}.\vadjust{\goodbreak}

Next,
by (\ref{as}), Lemma~\ref{L1} and (\ref{es}), uniformly for all
$t\ge0$,
\[
A(t)=a+S(t)=a+\E S(t)+o_{\mathrm p}(n)
=(n-a)\pi(t)+a+o_{\mathrm p}(n)
\]
and thus, using also (\ref{dtv}),
\[
n\qw A(t)
=(1-\gth)\pi(t)+\gth+o_{\mathrm p}(1)
=(1-\gth)\tpi(t)+\gth+o_{\mathrm p}(1).
\]
Substituting $t=xn$, we find by (\ref{tpi}), since $tp=xc+o(x)$,
uniformly in
all $x\ge0$,
\begin{eqnarray*}
n\qw A(xn) &=& (1-\gth){\mathbb P}\bigl(\Po(tp) \ge r\bigr)+\gth
+o_{\mathrm p}(1)\\
&=& (1-\gth){\mathbb P}\bigl(\Po(cx) \ge r\bigr)+\gth+o_{\mathrm p}(1)
\end{eqnarray*}
and, recalling (\ref{f1a}), still uniformly in $x\ge0$,
%
%
\begin{equation}\label{f1c}
n\qw\bigl(A(xn)-xn\bigr)
=f(x,c,\gth)+o_{\mathrm p}(1).
\end{equation}
Let $\eps>0$. Since $f(x,c,\gth)>0$ for $x\in[0,x_0(\theta))$, and
thus by
compactness $f(\cdot,c,\gth)$ is bounded from below on $[0,x_0(\theta
)-\eps]$,
(\ref{f1c}) implies that w.h.p.
$A(xn)-xn>0$ on $[0,x_0(\theta)-\eps]$, and thus $T>(x_0(\theta
)-\eps)n$.
Furthermore, both in~\ref{pn=cnsub} and in~\ref{pn=cnsuper} with
$\gth\neq\gthc$, we have $\ddx f(x_0(\theta),c,\gth)\neq0$ and
thus if
$\eps>0$
is small
enough, $f(x_0(\theta)+\eps,c,\gth)<0$, so (\ref{f1c}) implies that w.h.p.
$A((x_0(\theta)+\eps)n) <(x_0(\theta)+\eps)n$ and thus
$T<(x_0(\theta
)+\eps)n$.
\end{pf*}

The proof of Theorem~\ref{cngauss} is very similar to the one of
Theorem~\ref{Tac}. We
first give a more precise estimate of the process $S(t)$, which is the
analog of Lemma~\ref{LG} in the case $p= c/n$. However, in this case,
we get a
Brownian bridge because here we consider a large part of the
distribution of
$Y_i$.
%
%
\begin{lemma}\label{step1}
Suppose $r\ge2$, $p= c/n$ and $a \sim\theta n$ with $c>0$ and \mbox{$0<\theta< 1$}.
Then
%
%
\begin{equation}\label{eqstep1}
Z(x):=
\frac{S(x n)-\E S(x n)}{\sqrt{(1- \theta)n}}
\dto
W_0(\psi(cx))
\end{equation}
in $D[0,1]$, where $W_0$ is a Brownian bridge and $\psi(y):={\mathbb
P}(\Po(y)\ge r)$ as in~(\ref{tpi}).
\end{lemma}
\begin{pf}
Let $\tS(u):=\sum_{i=1}^{n-a} \ett\{U_i\le u\}$, $0\le u\le1$, where
$U_i\in U(0,1)$ are i.i.d.
By (\ref{st}) and (\ref{pi}), we have $S(t)\eqd\tS(\pi(t))$, jointly
for all \mbox{$t\ge0$}. Further, $\frac1{n-a}\tS(u)$, $u\in\oi$, is the empirical
distribution function of $U_1,\ldots,U_{n-a}$, and thus by
\cite{Bill}, Theorem 16.4, in $D\oi$,
\[
\frac{\tS(u)-\E\tS(u)}{\sqrt{n-a}} \dto W_0(u).
\]
Furthermore, by (\ref{dtv}) and (\ref{tpi}),
\[
\pi(xn)=\tpi(xn)+O(1/n)=\psi(xnp)+O(1/n)=\psi(cx)+O(1/n),\vadjust{\goodbreak}
\]
uniformly for $x\ge0$, and it follows, using the continuity of $W_0$, that
\[
\frac{S(xn)-\E S(xn)}{\sqrt{n-a}}
\eqd
\frac{\tS(\pi(xn))-\E\tS(\pi(xn))}{\sqrt{n-a}}
\dto W_0(\psi(cx))
\]
in $D\oi$, which proves the result since $n-a\sim(1-\gth)n$.
\end{pf}
\begin{pf*}{Proof of Theorem~\ref{cngauss}}
It suffices to consider $a$ such that $a\sim\gthc n$.
By (\ref{es}), (\ref{dtv}) and (\ref{tpi}),
%
%
\begin{eqnarray}\label{q2}
\E S(xn)&=&(n-a)\pi(xn)=(n-a)\tpi(xn)+O(1)\nonumber\\[-8pt]\\[-8pt]
&=&(n-a)\psi(cx)+O(1).\nonumber
\end{eqnarray}
By the Skorohod coupling theorem (\cite{Kallenberg}, Theorem 4.30),
we may assume that the processes for different $n$ are coupled such
that the limit (\ref{eqstep1}) in Lemma~\ref{step1} holds a.s., and not
just in
distribution.
Since convergence in $D[0,1]$ to a continuous function is equivalent
to uniform convergence, this means that a.s.
$Z(x)\to W_0(\psi(cx))$
uniformly for $x\in\oi$.
Hence, we have, using (\ref{q2}) and (\ref{f1a}),
%
%
\begin{eqnarray}\label{q1}\quad
&&
A(xn) - xn \nonumber\\
&&\qquad = a + S(xn) - xn\nonumber\\
&&\qquad= a + \E S(xn) + \sqrt{(1-\gthc)n} Z(x)
- xn
\nonumber\\
&&\qquad = a + (n-a)\psi(cx) + \sqrt{(1-\gthc)n} Z(x) - xn+O(1)
\\
&&\qquad = (a-\gthc n)\bigl(1-\psi(cx)\bigr)+nf(x,c,\gthc) + \sqrt{(1-\gthc)n}
Z(x) + O(1)
\nonumber\\
&&\qquad = nf(x,c,\gthc)+(a-\gthc n)\bigl(1-\psi(cx)\bigr) + \sqrt{(1-\gthc)n}
W_0(\psi(cx))
\nonumber\\
&&\qquad\quad{} +o_{\mathrm p}(n^{1/2}),
\nonumber
\end{eqnarray}
uniformly for $x\in\oi$.

We first use (\ref{q1}) to derive the simple estimate
%
%
\begin{equation}\label{q3}
A(xn)-xn=nf(x,c,\gthc)+o_{\mathrm p}(n),
\end{equation}
uniformly for $x\in\oi$.
By Lemma~\ref{Lf}, $f(x,c,\gthc)=0$ for $x=\xo$ or $x=\xox$, with
$f(x,c,\gthc)>0$ for $x\in[0,\xo)\cup(\xo,\xox)$ and $f(x,c,\gthc
)<0$ for
$x\in(\xox,1]$.
Hence, for any fixed small $\eps>0$, (\ref{q3}) implies that w.h.p.
$A(xn)-xn>0$ for $x\in[0,\xo-\eps]\cup[\xo+\eps,\xox-\eps]$
and $A(xn)-xn<0$ for $x\in[\xox+\eps,1]$, and hence
$T\in[\xo-\eps,\xo+\eps]\cup[\xox-\eps,\xox+\eps]$.
It follows by a standard argument that there exists a sequence $\eps
_n\downto0$
such that w.h.p.
\[
\Ax=T\in[\xo-\eps_n,\xo+\eps_n]\cup[\xox-\eps_n,\xox+\eps_n].
\]
Moreover, w.h.p. $T\in[\xo-\eps_n,\xo+\eps_n]$ if and only if
$\inf_{[\xo-\eps_n,\xo+\eps_n]}(A(xn)-xn)<0$.
(We may also
assume that $\eps_n$ is so small that $\eps_n<x_0$ and $2\eps_n<x_1-x_0$.)\vadjust{\goodbreak}

For $x\in[\xo-\eps_n,\xo+\eps_n]$, we have by (\ref{q1}) again,
and the continuity of~$\psi$ and $W_0$,
%
%
\begin{eqnarray}\label{q4}
A(xn) - xn&=& nf(x,c,\gthc)+(a-\gthc n)\bigl(1-\psi(c\xo)+o(1)\bigr)
\nonumber\\[-8pt]\\[-8pt]
&&{}+ \sqrt{(1-\gthc)n} W_0(\psi(c\xo)) +o_{\mathrm p}(n^{1/2}).
\nonumber
\end{eqnarray}
Further, $f(x_0,c,\gthc)=0$ and $f(x,c,\gthc)\ge0$ for
$x\in[\xo-\eps_n,\xo+\eps_n]$, and thus~(\ref{q4}) yields
%
%
\begin{eqnarray}\label{q5}\hspace*{32pt}
\inf_{x\in[\xo-\eps_n,\xo+\eps_n]} \bigl(A(xn) - xn\bigr)&=&(a-\gthc
n)\bigl(1-\psi(c\xo)+o(1)\bigr)
\nonumber\\[-4pt]\\[-12pt]
&&{}+ \sqrt{(1-\gthc)n} W_0(\psi(c\xo)) +o_{\mathrm p}(n^{1/2}).
\nonumber
\end{eqnarray}
The cases~\ref{cngauss-} and~\ref{cngauss+} are easily derived.
We thus focus on~\ref{cngauss0}.
We then have, from (\ref{q5}),
\begin{eqnarray*}
&&n\qqw\inf_{x\in[\xo-\eps_n,\xo+\eps_n]}
\bigl(A(xn) - xn\bigr) \\
&&\qquad= y\bigl(1 - \psi(c\xo)\bigr) +
\sqrt{1-\gthc} W_0(\psi(c\xo)) + o_{\mathrm p}(1)
\end{eqnarray*}
and thus,
since $(1-\psi(c\xo))\qw\sqrt{1-\gthc} W_0(\psi(c\xo))\in
N(0,\gss
)$, where
$\gss=(1-\gthc)\psi(c\xo)/(1-\psi(c\xo))>0$,
\begin{eqnarray*}
&&
{\mathbb P}\Bigl(\inf_{x\in[\xo-\eps_n,\xo+\eps_n]} \bigl(A(xn) - xn\bigr) <0\Bigr)
\\
&&\qquad={\mathbb P}\bigl(y\bigl(1-\psi(c\xo)\bigr)+ \sqrt{1-\gthc} W_0(\psi
(c\xo))<0\bigr) +o_{\mathrm p}(1)
\\
&&\qquad=1-\Phi(y/\gs) +o_{\mathrm p}(1).
\end{eqnarray*}
The result follows.
\end{pf*}

To prove Theorem~\ref{sqrtn} ($p \sim cn^{-1/r}$), we first show using the
previous results that if we can activate $\omega(n) \to\infty$ vertices,
then the activation spreads w.h.p. to the entire graph. It remains to show
that starting with a finite number of active vertices, the process activates
$\omega(n)$ vertices with a probability bounded away from $0$ and $1$. This
will be done using a branching process argument.
%
%
\begin{lemma}\label{Lsqrtn}
Suppose that $p\ge cn^{-1/r}$ for some $c>0$.
If $\go(n)\to\infty$,
then w.h.p. $A(t)>t$ for all $t$ with $\go(n)\le t\le n-1$.
\end{lemma}
\begin{pf}
This is easy to prove directly, but we prefer to view it as a~corollary
of our estimates for smaller $p$.
Thus, let $\xp:=\go(n)^{-1/2r} n^{-1/r}$.
We may assume $\go(n)\le n$ and then $n\qw\ll\xp\ll n^{-1/r}$, so
$\xp<
p$, at least for large~$n$, and we may assume that
$G_{n,\xp}\subseteq G_{n,p}$.
We may consider bootstrap percolation on $G_{n,\xp}$ and $G_{n,p}$
simultaneously, with the same initial set $\aox$ of size $a$; we use the\vadjust{\goodbreak}
description in Section~\ref{Ssetup}, starting with families of i.i.d. random
indicators $I_i'(s)\in\Be(\xp)$ and $I_i(s)\in\Be(p)$ where we
may assume $I_i'(s)\le I_i(s)$. Then, using $'$ to denote variables for
$G_{n,\xp}$, $S'(t)\le S(t)$ and $A'(t)\le A(t)$.

We apply Lemma~\ref{Lbulk} to $G_{n,\xp}$.
The critical time for $G_{n,\xp}$ is [see (\ref{tc})]
\[
\tc'=O\bigl((n(\xp)^r)^{-1/(r-1)}\bigr)=\go(n)^{1/2(r-1)}=o(\go(n)).
\]
Further, $\xp\ge n^{-3/2r}\ge n^{-3/4}$ so, by (\ref{bc}), $\bc'\to0$,
and we may choose ${\bcx}'$ with ${\bcx}'\to0$.
Hence, Lemma~\ref{Lbulk} shows that w.h.p.
$A(t)\ge A'(t)>t$ for $t\in[3\tc',n-{\bcx}']$, and the result
follows since,
for large $n$,
$3\tc'\le\go(n)$ and $n-{\bcx}'>n-1$.~%
\end{pf}

\begin{pf*}{Proof of Theorem~\ref{sqrtn}}
For (ii), we apply Lemma~\ref{Lsqrtn} (if necessary with a smaller
$c$). Taking
$\go(n)=a$, we see that w.h.p. $A(t)>t$ for all $t\in[a,n-1]$. Since also
$A(t)\ge a$, we have $A(t)>t$ for all $t\le n-1$, and thus $\Ax=T=n$.

For (i)
suppose $r \geq2$, $p\sim c n^{-1/r}$ and let $a \geq r$ be some
constant. The
probability that a vertex is activated at a given time $k$ is by (\ref{yik})
%
%
\begin{equation}\label{samu}
{\mathbb P} ( \yix= k) = \pmatrix{k-1\cr r-1}p^r (1-p)^{k-r}
\sim\pmatrix{k-1\cr r-1} \frac{c^r}{n}.
\end{equation}

For any fixed $K$, the random variables
\[
X_k:= A(k) - A(k-1) = S(k) - S(k-1) = \sum_{i\notin\cao} \ett\{
Y_i=k\},
\]
$ k = 1,\ldots, K$, form together with
\[
X_{K+1}:=n-a-A(K)=\sum_{i\notin\cao}\ett\{Y_i>K\}
\]
a random vector with the multinomial distribution
$\Mul(n-a,(p_k)_{k=1}^{K+1})$ with $p_k={\mathbb P}(\yix=k)$,
$k\le K$,
and $p_{K+1}={\mathbb P}(\yix\ge K+1)$.
By (\ref{samu}), $(n-a)p_k\to{k-1\choose r-1}c^r$ for $k\le K$, and it
follows that $X_k$ for $k\le K$ have a joint Poisson limit,
%
%
\begin{equation}\label{sjw}
(X_k)_{k=1}^K\dto(\xi_k)_{k=1}^K
\qquad\mbox{with } \xi_k\in\Po\biggl(\pmatrix{k-1\cr
r-1}c^r\biggr)\quad
\mbox{independent}.\hspace*{-40pt}
\end{equation}
Using the notation of Remark~\ref{RMS} we thus obtain
\[
A(k)\dto a+\sum_{j=1}^k\xi_j =a+k+\tsk\qquad
\mbox{for $k=1,\ldots,t$}\qquad \mbox{jointly}
\]
and thus ${\mathbb P}(T=k)\to{\mathbb P}(\tT=k)$ for $k\le K$
and ${\mathbb P}(T>K)\to{\mathbb P}(\tT>K)$.

Since $K$ is arbitrary, we have shown
${\mathbb P}(\Ax=k)={\mathbb P}(T=k)\to{\mathbb P}(\tT=k)=\zeta(a,c,k)$
for every finite $k\ge1$.
Furthermore, ${\mathbb P}(T>K)-{\mathbb P}(\tT>K)\to0$ for any fixed
$K$, and
a standard argument shows\vspace*{1pt} that there exists a sequence $K_n\to\infty$ such
that ${\mathbb P}(T>K_n)-{\mathbb P}(\tT>K_n)\to0$, and thus
${\mathbb
P}(T>K_n)\to{\mathbb P}(\tT=\infty)$.
On the other hand,\vadjust{\goodbreak} Lemma~\ref{Lsqrtn} with $\go(n)=K_n$ shows that
$P(K_n\le
T<n)\to0$. Consequently,
${\mathbb P}(T=n)={\mathbb P}(T>K_n)+o(1)\to{\mathbb P}(\tT=\infty
)=\zeta(a,c)$.

It is clear that $\zeta(a,c,k)={\mathbb P}(\tT=k)>0$ for every $k\ge a$.
To see that also $\zeta(a,c)={\mathbb P}(\tT=\infty)>0$, note that, see
(\ref{sjw}),
$\E\xi_k={k-1\choose r-1}c^r \to\infty$ as $k\to\infty$. Hence,
there is
some $K_0$ such that $\E\xi_{K_0}>1$. Since $\xi_k$ stochastically dominates
$\xi_{K_0}$ for $k\ge K_0$, it follows that if the process reaches $K_0$
without stopping, the continuation dominates (up to a change of time)
a~Galton--Watson branching process
with offspring distribution $\xi_{K_0}$, which is supercritical and
thus has a positive probability of living forever. Hence, ${\mathbb
P}(\tT=\infty)>0$.
\end{pf*}
\begin{pf*}{Proof of Theorem~\ref{sqrtnn}}
It suffices to consider $a=r$. Thus assume $a=r$, and consider the
vertices activated in the first generation, that is, at time $t=r$. There
are $S(r)\in\Bin(n-r,p^r)$ such vertices.
[Note that, see~(\ref{pi}),
$\pi(r)={\mathbb P}(\Bin(r,p)=r)=p^r$.]
Consequently,
$\E S(r)=(n-r)p^r\to\infty$. Let $\go(n)=\E S(r)/2$, so $\go(n)\to
\infty$.
It follows from Chebyshev's inequality (or Chernoff's) that w.h.p.
$S(r)>\go(n)$. Hence, w.h.p. for all $t\in[r,\go(n)]$,
$A(t)\ge A(r)>S(r)>\go(n)\ge t$. Together with the trivial $A(t)\ge a=r>t$
for $t<r$ and Lemma~\ref{Lsqrtn}, this shows that w.h.p. $A(t)>t$ for all
$t\le
n-1$, and thus $\Ax=T=n$.
\end{pf*}

\section*{Acknowledgments}

The authors gratefully acknowledge the hospitality and the stimulating
environment of Institut Mittag-Leffler where the majority of this work
was carried out during the program ``Discrete Probability,'' 2009. The
authors thank the referee for helpful suggestions.


%

\printaddresses

\end{document}